\documentclass{amsart}

\newif\iflics
\licsfalse

\usepackage{amsmath,amsthm,amsfonts}
\usepackage{cite}
\usepackage[utf8]{inputenc}
\usepackage{tikz-cd}
\usepackage{xspace}
\usepackage[all,cmtip]{xy}
\usepackage{mathtools}
\usepackage{color}
\definecolor{myurlcolor}{rgb}{0.4,0,0}
\definecolor{mycitecolor}{rgb}{0,0,0.8}
\definecolor{myrefcolor}{rgb}{0,0,0.8}
\usepackage[pagebackref]{hyperref}
\usepackage{hyperref}
\hypersetup{colorlinks,linkcolor=myrefcolor,citecolor=mycitecolor,urlcolor=myurlcolor}

\usepackage{tikz}
\usepackage{tikz-cd}
  \usetikzlibrary{backgrounds,circuits,circuits.ee.IEC,shapes,petri,fit,matrix}
\usepackage{adjustbox} 

\usepackage[status=draft,mode=multiuser]{fixme}
\fxusetheme{colorsig}

\usepackage{cleveref,aliascnt}

\newcommand{\defin}[1]{{\bf \boldmath{#1}}}

\newcommand{\id}[1]{\mathrm{Id}_{#1}}
\newcommand{\op}{^\mathrm{op}}
\newcommand{\PN}{\ensuremath{\mathcal{P}}\xspace}
\newcommand{\N}{\mathbb{N}}
\newcommand{\SigmaNet}{\ensuremath{\Sigma\text{-net}}\xspace}

\newcommand{\CatDef}[1]{\ensuremath{\mathsf{#1}}}

\newcommand{\Set}{\CatDef{Set}}
\newcommand{\SSMC}{\CatDef{SSMC}}
\newcommand{\StrMC}{\CatDef{StrMC}}
\newcommand{\CMC}{\CatDef{CMC}}
\newcommand{\PreNet}{\CatDef{PreNet}}
\newcommand{\SigmaNetCat}{\CatDef{\Sigma{\text{-}net}}}
\newcommand{\Petri}{\CatDef{Petri}}
\newcommand{\WGPet}{\CatDef{WGPet}}
\newcommand{\Cat}{\CatDef{Cat}}
\newcommand{\CommDuoid}{\CatDef{CommDuoid}}
\newcommand{\PROP}{\CatDef{PROP}}
\newcommand{\Open}{\mathbb{O}\mathbf{pen}}
\newcommand{\pre}{\mathrm{pre}}
\newcommand{\pet}{\mathrm{pet}}

\newcommand{\X}{\CatDef{X}}
\newcommand{\C}{\CatDef{C}}
\newcommand{\D}{\CatDef{D}}
\newcommand{\E}{\CatDef{E}}
\newcommand{\Comma}[2]{(#1\downarrow #2)}

\newcommand{\maps}{\colon}

\newcommand{\simRightarrow}{\xRightarrow{\raisebox{-3pt}[0pt][0pt]{\ensuremath{\sim}}}}
\newcommand\NEarrow{\mathrel{\rotatebox[origin=c]{45}{$\Rightarrow$}}}
\newcommand\Narrow{\mathrel{\rotatebox[origin=c]{90}{$\Rightarrow$}}}

\newtheorem{thm}{Theorem}[section]
\crefname{thm}{Theorem}{Theorems}
\newaliascnt{prop}{thm}
\newtheorem{prop}[prop]{Proposition}
\aliascntresetthe{prop}
\crefname{prop}{Proposition}{Propositions}
\newaliascnt{lem}{thm}
\newtheorem{lem}[lem]{Lemma}
\aliascntresetthe{lem}
\crefname{lem}{Lemma}{Lemmas}
\theoremstyle{definition}
\newaliascnt{defn}{thm}
\newtheorem{defn}[defn]{Definition}
\aliascntresetthe{defn}
\crefname{defn}{Definition}{Definitions}
\theoremstyle{remark}
\newaliascnt{ex}{thm}
\newtheorem{ex}[ex]{Example}
\aliascntresetthe{ex}
\crefname{ex}{Example}{Examples}
\newaliascnt{rmk}{thm}
\newtheorem{rmk}[rmk]{Remark}
\aliascntresetthe{rmk}
\crefname{rmk}{Remark}{Remarks}

\crefformat{section}{Section~#2#1#3}

\tikzstyle{simple}=[-,line width=2.000]
\tikzstyle{arrow}=[-,postaction={decorate},decoration={markings,mark=at position .5 with {\arrow{>}}},line width=1.100]
\pgfdeclarelayer{edgelayer}
\pgfdeclarelayer{nodelayer}
\pgfsetlayers{edgelayer,nodelayer,main}

\tikzstyle{none}=[inner sep=0pt]

\tikzstyle{place}=
[circle,thick,draw=blue!75,fill=blue!20,minimum size=6mm]

\tikzstyle{transition}=
[rectangle,thick,draw=black!75,fill=black!20,minimum width=1mm, minimum height=4mm]

\def\xcoord{-1.2}
\def\ycoord{0.5}

\def\transitionbgrd{gray!50}
\def\transitionop{0.3}

\tikzstyle{inarrow}=[->, >=stealth, shorten >=.03cm,line width=1.5]
\tikzstyle{empty}=[circle,fill=none, draw=none]
\tikzstyle{inputdot}=[circle,fill=purple,draw=purple, scale=.25]
\tikzstyle{inputarrow}=[->,draw=purple, shorten >=.05cm]
\tikzstyle{simple}=[-,draw=purple,line width=1.000]

\makeatletter
\def\slashedarrowfill@#1#2#3#4#5{%
  $\m@th\thickmuskip0mu\medmuskip\thickmuskip\thinmuskip\thickmuskip
   \relax#5#1\mkern-7mu%
   \cleaders\hbox{$#5\mkern-2mu#2\mkern-2mu$}\hfill
   \mathclap{#3}\mathclap{#2}%
   \cleaders\hbox{$#5\mkern-2mu#2\mkern-2mu$}\hfill
   \mkern-7mu#4$%
}
\def\rightslashedarrowfill@{%
  \slashedarrowfill@\relbar\relbar\mapstochar\rightarrow}
\newcommand{\xslashedrightarrow}[2][]{%
  \ext@arrow 0055{\rightslashedarrowfill@}{#1}{#2}}

\FXRegisterAuthor{ms}{anms}{Mike}
\FXRegisterAuthor{jb}{anjb}{John}
\FXRegisterAuthor{jm}{anjm}{Jade}
\FXRegisterAuthor{fg}{anfg}{Fab}

\title{Categories of Nets}
\date{\today}

\iflics
\author{
        \IEEEauthorblockN{
            John C.\ Baez\IEEEauthorrefmark{1}, 
            Fabrizio Genovese\IEEEauthorrefmark{2}\IEEEauthorrefmark{3}, 
            Jade Master\IEEEauthorrefmark{4}
            and Michael Shulman\IEEEauthorrefmark{5}
        }

    \IEEEauthorblockA{
        \IEEEauthorrefmark{1}Department of Mathematics, U.\ C.\ Riverside, Riverside, California 92521\\
        Email: baez@math.ucr.edu\\
    }
    \IEEEauthorblockA{
        \IEEEauthorrefmark{2}Department of Informatics, University of Pisa, Pisa, Italia 56127\\
    }
    \IEEEauthorblockA{
        \IEEEauthorrefmark{3}Statebox, Statebox.IO B.V., Leiden, Netherlands 2314ED\\
        Email: fabrizio.romano.genovese@gmail.com\\
    }
    \IEEEauthorblockA{
        \IEEEauthorrefmark{4}Department of Mathematics, U.\ C.\ Riverside, Riverside, California 92521\\
        Email: jmast003@ucr.edu\\
    }
    \IEEEauthorblockA{
        \IEEEauthorrefmark{5}Department of Mathematics, University of San Diego, San Diego, California 92110\\
        Email: shulman@sandiego.edu
    }
  }

\else

\author{John C.\ Baez}
\address{Department of Mathematics, U.\ C.\ Riverside, Riverside, California 92521}
\email{baez@math.ucr.edu}

\author{Fabrizio Genovese}
\address{Department of Informatics, University of Pisa, Pisa, Italia 56127}
\email{fabrizio.romano.genovese@gmail.com}

\author{Jade Master}
\address{Department of Mathematics, U.\ C.\ Riverside, Riverside, California 92521}
\email{jmast003@ucr.edu}

\author{Michael Shulman}
\address{Department of Mathematics, University of San Diego, San Diego, California 92110}
\email{shulman@sandiego.edu}

\fi

\usepackage{etoolbox}
\apptocmd{\sloppy}{\hbadness 10000\relax}{}{}

\begin{document}

\iflics
\IEEEoverridecommandlockouts
\IEEEpubid{\makebox[\columnwidth]{978-1-6654-4895-6/21/\$31.00~
\copyright2021 IEEE \hfill} \hspace{\columnsep}\makebox[\columnwidth]{ }}
\fi

\iflics\maketitle\fi


\begin{abstract} 
We present a unified framework for Petri nets and various variants, such as pre-nets and Kock's whole-grain Petri nets.  Our framework is based on a less well-studied notion that we call $\Sigma$-nets, which allow fine-grained control over whether each transition behaves according to the collective or individual token philosophy.  We describe three forms of execution semantics in which pre-nets generate strict monoidal categories, $\Sigma$-nets (including whole-grain Petri nets) generate symmetric strict monoidal categories, and Petri nets generate commutative monoidal categories, all by left adjoint functors.  We also construct adjunctions relating these categories of nets to each other, in particular showing that all kinds of net can be embedded in the unifying category of $\Sigma$-nets, in a way that commutes coherently with their execution semantics.
\end{abstract}

\iflics\else\maketitle\fi

\section{Introduction}
\label{sec:introduction}

A Petri net is a seemingly simple thing: 
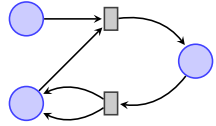
\[
  \scalebox{0.75}{
\begin{tikzpicture}
  \begin{scope}
    \begin{pgfonlayer}{nodelayer}
        \node [place,tokens=0] (1a) at (-1.5,0.75) {};
        \node [place,tokens=0] (1b) at (-1.5,-0.75) {};      			
        \node [place,tokens=0] (3a) at (1.5,0) {};
        \node[transition] (2a) at (0,0.75) {};
        \node[transition] (2b) at (0,-0.75) {};
    \end{pgfonlayer}
    \begin{pgfonlayer}{edgelayer}
        \draw[style=inarrow, thick] (1a) to (2a);
        \draw[style=inarrow, thick] (1b) to (2a);
        \draw[style=inarrow, thick, bend right] (2b) to (1b);
        \draw[style=inarrow, thick, bend left] (2b) to (1b);
    
        \draw[style=inarrow, thick, bend left] (2a) to (3a);
        \draw[style=inarrow, thick, bend left] (3a) to (2b);
    \end{pgfonlayer}
  \end{scope}
\end{tikzpicture}
}
\]
\noindent It consists of ``places'' (drawn as circles) and ``transitions'' (drawn as boxes), with directed edges called ``arcs'' from places to transitions and from transitions to places.  The idea is that when we use a Petri net, we place dots called ``tokens'' in the places, and then move them around using the transitions:
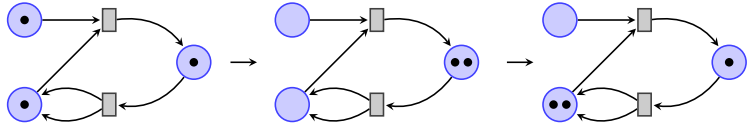
\[
  \scalebox{0.75}{  
\begin{tikzpicture}
  \begin{scope}
    \begin{pgfonlayer}{nodelayer}
        \node [place,tokens=1] (1a) at (-1.5,0.75) {};
        \node [place,tokens=1] (1b) at (-1.5,-0.75) {};      			
        \node [place,tokens=1] (3a) at (1.5,0) {};
        \node[transition] (2a) at (0,0.75) {};
        \node[transition] (2b) at (0,-0.75) {};
    \end{pgfonlayer}
    \begin{pgfonlayer}{edgelayer}
        \draw[style=inarrow, thick] (1a) to (2a);
        \draw[style=inarrow, thick] (1b) to (2a);
        \draw[style=inarrow, thick, bend right] (2b) to (1b);
        \draw[style=inarrow, thick, bend left] (2b) to (1b);
    
        \draw[style=inarrow, thick, bend left] (2a) to (3a);
        \draw[style=inarrow, thick, bend left] (3a) to (2b);
    \end{pgfonlayer}
  \end{scope}

  \begin{scope}[xshift=135]
    \begin{pgfonlayer}{nodelayer}
        \node [place,tokens=0] (1a) at (-1.5,0.75) {};
        \node [place,tokens=0] (1b) at (-1.5,-0.75) {};      			
        \node [place,tokens=2] (3a) at (1.5,0) {};
        \node[transition] (2a) at (0,0.75) {};
        \node[transition] (2b) at (0,-0.75) {};
    \end{pgfonlayer}
    \begin{pgfonlayer}{edgelayer}
        \draw[style=inarrow, thick] (1a) to (2a);
        \draw[style=inarrow, thick] (1b) to (2a);
        \draw[style=inarrow, thick, bend right] (2b) to (1b);
        \draw[style=inarrow, thick, bend left] (2b) to (1b);
    
        \draw[style=inarrow, thick, bend left] (2a) to (3a);
        \draw[style=inarrow, thick, bend left] (3a) to (2b);
    \end{pgfonlayer}
  \end{scope}

    \draw[style=inarrow, thick] (2.15,0) -- (2.65,0);

  \iflics\else
  \begin{scope}[xshift=270]
    \begin{pgfonlayer}{nodelayer}
        \node [place,tokens=0] (1a) at (-1.5,0.75) {};
        \node [place,tokens=2] (1b) at (-1.5,-0.75) {};      			
        \node [place,tokens=1] (3a) at (1.5,0) {};
        \node[transition] (2a) at (0,0.75) {};
        \node[transition] (2b) at (0,-0.75) {};
    \end{pgfonlayer}
    \begin{pgfonlayer}{edgelayer}
        \draw[style=inarrow, thick] (1a) to (2a);
        \draw[style=inarrow, thick] (1b) to (2a);
        \draw[style=inarrow, thick, bend right] (2b) to (1b);
        \draw[style=inarrow, thick, bend left] (2b) to (1b);
    
        \draw[style=inarrow, thick, bend left] (2a) to (3a);
        \draw[style=inarrow, thick, bend left] (3a) to (2b);
    \end{pgfonlayer}
  \end{scope}
    \draw[style=inarrow, thick] (7.05,0) -- (7.55,0);
  \fi

\end{tikzpicture}
}
\]
Thanks in part to their simplicity, Petri nets are widely used in computer science, chemistry, biology and other fields to model systems where entities interact and change state \cite{GiraultValk,Peterson}.

Ever since the work of Meseguer and Montanari \cite{MM}, parallels have been drawn between Petri nets and symmetric strict monoidal categories (SSMCs). Intuitively, a Petri net can be interpreted as a \emph{presentation} of such a category, by using its places to generate a commutative monoid of objects, and its transitions to generate the morphisms.   An object in the SSMC represents a ``marking'' of the net---a given placement of tokens in it---while a morphism represents a ``firing sequence'': a sequence of transitions that carry markings to other markings. One of the advantages of this ``execution semantics'' for a net is that it can be compositionally interfaced to other structures using monoidal functors. 

However, the apparent simplicity of Petri nets hides many subtleties.  There are various ways to make the definition of Petri net precise.  For example: is there a finite \emph{set} of arcs from a given place to a given transition (and the other way around), or merely a natural number?  If there is a finite set, is this set equipped with an ordering?  Furthermore, what is a morphism between Petri nets?  A wide variety of answers to these questions have been explored in the literature.  

Different answers are good for different purposes.  In the ``individual token philosophy'', we allow a finite set of tokens in each place, and tokens have their own individual identity.   In the ``collective token philosophy'', we merely allow a natural number of tokens in each place, so it means nothing to switch two tokens in the same place~\cite{GlabbeekPlotkin}.

Moreover, the idea of using SSMCs to represent net semantics, albeit intuitive, presents subtleties of its own.  There has been a great deal of work on this subject \cite{BM, BMMS2001, DMM, GH, GGHPPV, Master, SassoneStrong, SassoneCategory, SassoneAxiomatization}.   Nevertheless, we still lack a general answer describing the relations between nets and SSMCs.

\section{Dramatis person\ae}
  \label{sec: dramatis personae}
Our goal is to bring some order to this menagerie.  Our attitude is that though there may be multiple kinds of Petri net, each should freely generate a monoidal category of an appropriate sort, and these processes should be left adjoint functors.  
We thus consider three kinds of nets, and three corresponding kinds of monoidal categories:
\[
  \begin{tikzcd}[row sep=large,column sep=large]
  \StrMC \ar[r,shift left] \ar[dd,shift left] &
  \SSMC \ar[l,shift left] \ar[r,shift left] \ar[dd, shift left] &
  \CMC \ar[dd,shift left] \ar[l,shift left] \\
  \\
  \PreNet 
  \ar[uu,shift left] \ar[r,shift left=2] \ar[r,shift right=2] 
  &\SigmaNetCat \ar[uu,shift left] \ar[l] \ar[r,shift left]
  &\Petri \ar[uu,shift left] \ar[l,shift left]
\end{tikzcd}
\]
On the top row we have:
\begin{itemize}
  \item $\StrMC$, with strict monoidal categories as objects and strict monoidal functors as morphisms.
  \item $\SSMC$, with symmetric strict monoidal categories as objects and strict symmetric monoidal functors as their morphisms.  A \defin{symmetric strict monoidal category} is a symmetric monoidal category whose monoidal structure is strictly associative and unital; its symmetry may not be the identity.
  \item $\CMC$, with commutative monoidal categories as objects and strict symmetric monoidal functors as morphisms.  A \defin{commutative monoidal category} is a symmetric strict monoidal category where the symmetry is the identity. 
\end{itemize}
Monoidal categories of these three kinds are freely generated by three kinds of nets, on the bottom row of the diagram:
\begin{itemize}
  \item $\PreNet$, with pre-nets as objects.   A \defin{pre-net} consists of a set $S$ of \defin{places}, a set $T$ of \defin{transitions}, and functions $T \xrightarrow{s,t} S^*\times S^*$, where $S^*$ is the underlying set of the free monoid on $S$.  We describe the category $\PreNet$, and its adjunction with $\StrMC$,  in \cref{sec:pre-nets}.  These ideas are due to Bruni, Meseguer, Montanari and Sassone~\cite{BMMS2001}.
  \item $\SigmaNetCat$, with \SigmaNet{}s as objects.  A \defin{\SigmaNet{}} consists of a set $S$ and a discrete opfibration $T \to PS \times PS\op$, where $PS$ is the free symmetric strict monoidal category generated by a set of objects $S$ and no generating morphisms.  We describe the category $\SigmaNetCat$ in \cref{sec:sigma nets}, and its adjunction with $\SSMC$ in \cref{thm:pre-net - sigma-net adjunction}.  
  \item $\Petri$, with Petri nets as objects.  A \defin{Petri net}, as we will use the term, consists of a set $S$, a set $T$, and functions $T \xrightarrow{s,t} \N[S] \times \N[S]$, where $\N[S]$ is the free commutative monoid on $S$.   We describe the category $\Petri$, and its adjunction with $\CMC$, in \cref{sec:petri nets}.  This material can be found in~\cite{BM,Master}.
 \end{itemize}

These three notions of net obviously have a similar flavor. Their parallel relationships to the three notions of monoidal category is made even clearer when we note that regarded as discrete categories, $S^*$ and $\N[S]$ are respectively the free monoidal category and the free commutative monoidal category on the set $S$.

Besides the three adjunctions between the categories on the top row and those on the bottom row, in which the left adjoints point upward, there are also adjunctions running horizontally across the diagram: adjoint pairs in the top row and bottom right and an adjoint triple in the bottom left, with left adjoints drawn above their right adjoints. In \cref{sec:adjunctions}, we examine these adjunctions in detail.

Of particular importance are the right adjoint mapping Petri nets to \SigmaNet{}s, and the left adjoint mapping pre-nets to \SigmaNet{}s.
We think of these as ``embedding'' the collective token world (Petri nets) and the individual token world (pre-nets) into the unifying context of \SigmaNet{}s.  In the case of Petri nets, the functor is literally an embedding (i.e., fully faithful), and since it is a right adjoint it preserves all limits (though not all colimits).  In the case of pre-nets, the functor is faithful but not full, but it is an equivalence onto a slice category of \SigmaNetCat{}, and preserves all colimits and all connected limits (such as pullbacks).  These embeddings also respect the most common categorical semantics: in \cref{sec:adjunctions} we will show that the left adjoints $\PreNet \to \SSMC$ and $\Petri\to \CMC$ both factor through $\SigmaNetCat$.

The images of pre-nets and Petri nets in \SigmaNet{}s have a large intersection, consisting of those nets in which no places are ever duplicated in the inputs or outputs of any transition.  These are the nets for which there is no difference between the individual and collective token philosophies.  As we shall see, general \SigmaNet{}s allow more fine-grained control than either pre-nets or Petri nets: for example, some transitions may obey the individual token philosophy while others obey the collective token philosophy.

 Our work is closely related to that of Kock~\cite{Kock}.
  He refers to \SigmaNet{}s as ``digraphical species'', and sketches a proof, different from ours, that there is an adjunction relating them to $\SSMC$.  But his focus is on a fourth notion of net: ``whole-grain Petri nets''.  He sketches a proof that these are the image of pre-nets inside $\SigmaNet$, which we detail in \cref{sec:whole-grain} (so that in particular, whole-grain Petri nets also generate symmetric strict monoidal categories); but says nothing about their relationship to Petri nets as traditionally conceived.
  
\vfill

\section{Petri Nets}
\label{sec:petri nets}
Symmetric monoidal categories are a general algebraic framework to represent processes that can be performed in sequence and in parallel. Because Petri nets represent schematics for such processes, we expect them to freely generate symmetric monoidal categories.   In fact they generate a special sort of
symmetric monoidal categories: \emph{commutative} ones.

\begin{defn}
  Let $\Petri$ be the category where:
  \begin{itemize}
    \item An object is a \defin{Petri net}: a pair of functions $T \xrightarrow{s,t} \N[S]$, where $\N[S]$ denotes the underlying set of the free commutative monoid on $S$.
    \item A morphism from $T_1 \xrightarrow{s_1,t_1} \N[S_1]$ to $T_2 \xrightarrow{s_2,t_2} \N[S_2]$ is a pair of functions $f \maps S_1 \to S_2, g \maps T_1 \to T_2$ such that the following diagram commutes:
    \[
      \begin{tikzcd}
        \N[S_1]\ar[d, "{\N[f]}"'] & T_1 \ar[d, "g"']  \ar[l, "s_1"'] \ar[r, "t_1"] & \N[S_1] \ar[d, "{\N[f]}"] \\
        \N[S_2] & T_2 \ar[l, "s_2"] \ar[r, "t_2"'] & \N[S_2]
      \end{tikzcd}  
    \]  
    where $\N[f]$ denotes the unique monoid homomorphism extending $f$.
  \end{itemize}
\end{defn}
Our concept of Petri net morphism is more restrictive than that of Meseguer and Montanari \cite{MM}: where we have $\N[f] \maps \N[S_1] \to \N[S_2]$; they allow an arbitrary monoid homomorphism from $\N[S_1]$ to $\N[S_2]$. 
\begin{defn}
  A \defin{commutative monoidal category} is a commutative monoid object in $\Cat$. Equivalently, it is a strict monoidal category $(C,\otimes,I)$ such that for all objects $a$ and $b$ and morphisms $f$ and $g$ in $C$ 
  \[a \otimes b = b \otimes a \text{ and } f \otimes g = g \otimes f.\]
  A morphism of commutative monoidal categories is a strict monoidal functor.  We write $\CMC$ for the category of commutative monoidal categories and such morphisms between them.
\end{defn}
A commutative monoidal category can be seen as a particularly strict sort of symmetric monoidal category. Ordinarily, symmetric monoidal categories are equipped with ``symmetry'' isomorphisms
\[\sigma_{x,y} \maps x \otimes y \xrightarrow{\sim} y \otimes x \]
for every pair of objects $x$ and $y$.
In a commutative monoidal category $x \otimes y$ is \textit{equal} to $y\otimes x$, so we can --- and henceforth will --- make it symmetric by choosing $\sigma_{x,y}$ to be the identity for all $x$ and $y$.  Any morphism of commutative monoidal categories then becomes a strict symmetric monoidal functor.   

The following adjunction shows that Petri nets are the right sort of generating data for commutative monoidal categories.
\begin{prop}\label{prop: adjunction Petri - CMC}
  There is an adjunction 
  \[
    \begin{tikzcd}
      \Petri \ar[r,bend left,"F_{\Petri}"]  \ar[r, phantom, "\bot"] & \CMC \ar[l, bend left,"U_{\Petri}"] 
    \end{tikzcd} 
  \]
  whose left adjoint sends a Petri net $P$ to the commutative monoidal category $F_{\Petri}(P)$ where:
  \begin{itemize}
      \item Objects are markings of $P$, i.e., elements of the free commutative monoid on its set of places.
      \item Morphisms are generated inductively by the following rules:
      \begin{itemize}
      \item for each place $s$ there is an identity $1_s \maps s \to s$
      \item for each transition $\tau$ of $P$, there is a morphism going from its source to its target
      \item for every pair of morphisms $f \maps x \to y$ and $f' \maps x' \to y'$, there is a morphism $f \otimes f' \maps x \otimes x' \to y \otimes y'$
      \item for every pair of composable morphisms $f \maps x \to y$ and $g \maps y \to z$, there is a morphism $g \circ f \maps x \to z$
      \end{itemize}
      and quotiented to satisfy the axioms of a commutative monoidal category. 
  \end{itemize}
\end{prop}
\begin{proof}
  This is a special case of \cite[Theorem 5.1]{Master}.  See also \cite[Lemma 9]{BM}.
\end{proof}
As noted in \cref{sec:introduction}, we view this construction as associating to each net a monoidal category of its ``executions''.   The objects of this category are markings that accord with the collective token philosophy: for instance, if $p$ and $q$ are places, the object $2p+3q$ has two tokens on $p$ and three tokens on $q$, but no way to distinguish between the former two tokens or between the latter three. Similarly, the morphisms in this category are equivalence classes of firing sequences.  This interpretation is particularly captivating if we represent morphisms in a monoidal category using string diagrams.

However, the equivalence relation on firing sequences that determines when two define the same morphism is very coarse when we take the commutative monoidal category freely generated by a Petri net.  Indeed, if $f,g \maps x \to x$ are morphisms in a commutative monoidal category, the following sequence of equations holds:
\begin{equation*}
  \scalebox{\iflics 0.6\else 0.85\fi}{
  \begin{tikzpicture}[baseline=0.45cm]
    \node[draw,thick,minimum width=0.6cm,minimum height=0.6cm, rounded corners=3pt] (f) at (1,1) {$f$};
    \node[draw,thick,minimum width=0.6cm,minimum height=0.6cm, rounded corners=3pt] (g) at (1,0) {$g$};
    \draw (0,1) -- node[above] {$x$} (f);
    \draw (f) -- node[above] {$x$}(3,1);
    \draw (0,0) -- node[below] {$x$} (g);
    \draw (g) -- node[below] {$x$} (3,0);
  \end{tikzpicture}
  =
  \begin{tikzpicture}[baseline=0.45cm]
    \node[draw,thick,minimum width=0.6cm,minimum height=0.6cm, rounded corners=3pt] (f) at (1,1) {$f$};
    \node[draw,thick,minimum width=0.6cm,minimum height=0.6cm, rounded corners=3pt] (g) at (1,0) {$g$};
    \draw (0,1) -- node[above] {$x$} (f);
    \draw (f) -- (1.5,1);
    \draw (0,0) -- node[below] {$x$} (g);
    \draw (g) -- (1.5,0);
    \draw[out=0, in=180] (1.5,0) to node[above,near end] {$x$} (3,1);
    \draw[out=0, in=180] (1.5,1) to node[below,near end] {$x$} (3,0);
  \end{tikzpicture}
  =
  \begin{tikzpicture}[baseline=0.45cm]
    \node[draw,thick,minimum width=0.6cm,minimum height=0.6cm, rounded corners=3pt] (f) at (2.25,0) {$f$};
    \node[draw,thick,minimum width=0.6cm,minimum height=0.6cm, rounded corners=3pt] (g) at (0.75,0) {$g$};
    \draw[out=0, in=180] (0,1) to node[near start, above] {$x$} (f);
    \draw (0,0) -- node[below] {$x$} (g);
    \draw [out=0, in=180] (g) to node[near end, above] {$x$} (3,1);
    \draw (f) --  node[below] {$x$} (3,0);
  \end{tikzpicture}
  =
  \begin{tikzpicture}[baseline=0.45cm]
    \node[draw,thick,minimum width=0.6cm,minimum height=0.6cm, rounded corners=3pt] (f) at (2.25,0) {$f$};
    \node[draw,thick,minimum width=0.6cm,minimum height=0.6cm, rounded corners=3pt] (g) at (0.75,0) {$g$};
    \draw (0,1) -- node[above] {$x$} (3,1);
    \draw (0,0) -- node[below] {$x$} (g);
    \draw (g) -- node[below] {$x$} (f) ;
    \draw (f) --  node[below] {$x$} (3,0);
  \end{tikzpicture}}
\end{equation*}
These equations imply that given any two firing sequences $f$ and $g$ that start and end at some marking $x$, the commutative monoidal category cannot distinguish whether they act independently or whether $f$ acts on the tokens already processed by $g$.  When $x$ is the tensor unit, the equations above hold in \emph{any} symmetric monoidal category.  But in a commutative monoidal category, the above equations hold for any object~$x$.
\section{Pre-nets}
\label{sec:pre-nets}
The shortcomings of commutative monoidal categories we presented in the last section are overcome by using symmetric monoidal categories where the symmetries are not necessarily identity morphisms.  One approach first builds monoidal categories and then freely adds symmetries.  In 1991 Joyal and Street~\cite{JoyalStreet} introduced  ``tensor schemes'', which can be used to describe free strict monoidal categories.  In 2001, essentially the same idea was introduced by Bruni, Meseguer, Montanari and Sassone~\cite{BMMS2001} under the name ``pre-nets''. However, for these authors, the use of pre-nets to describe free strict monoidal categories was just the first stage of a procedure to obtain free \emph{symmetric} strict monoidal categories.  We recall this procedure now.
\begin{defn}\label{defn:prenet}
  Let $\PreNet$ be the category where:
  \begin{itemize}
    \item An object is a \defin{pre-net}: a pair of functions $T \xrightarrow{s,t} S^*$, where $S^*$ is the underlying set of the free monoid on $S$.
    \item A morphism from $T_1 \xrightarrow{s_1,t_1} S_1^*$ to $T_2 \xrightarrow{s_2,t_2} S_2^*$ is a pair of functions $f \maps S_1 \to S_2$, $g \maps T_1 \to T_2$ such that the following diagram commutes,
    where $f^{*}$ denotes the unique monoid homomorphism extending $f$:
    \[
      \begin{tikzcd}
        S_1^{*}\ar[d, "f^{*}"'] & T_1 \ar[d, "g"']  \ar[l, "s_1"'] \ar[r, "t_1"] & S_1^{*} \ar[d, "f^{*}"] \\
        S_2^{*} & T_2 \ar[l, "s_2"] \ar[r, "t_2"'] & S_2^{*}
      \end{tikzcd}  
    \]  
  \end{itemize}
\end{defn}
\noindent Graphically, a pre-net looks very similar to a Petri net, and we follow the convention of~\cite{BBM} by decorating arcs with numbers to indicate their input/output order in a transition, as in:
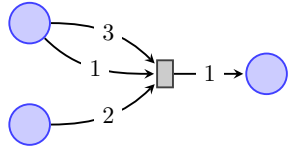
\[
    \scalebox{0.9}{  
\begin{tikzpicture}
  \begin{scope}
    \begin{pgfonlayer}{nodelayer}
        \node [place,tokens=0] (1a) at (-2,0.75) {};
        \node [place,tokens=0] (1b) at (-2,-0.75) {};      			
        \node [place,tokens=0] (3a) at (1.5,0) {};
        \node[transition] (2a) at (0,0) {};
    \end{pgfonlayer}
    \begin{pgfonlayer}{edgelayer}
        \draw[style=inarrow, out=0, in=135, thick] (1a) to node[fill=white, midway] {$3$} (2a);
        \draw[style=inarrow, out=-45, in=180, thick] (1a) to node[fill=white, midway] {$1$} (2a);
        \draw[style=inarrow, out=0, in=225, thick] (1b) to node[fill=white,midway] {$2$} (2a);
        \draw[style=inarrow, thick] (2a) to node[fill=white, midway] {$1$} (3a);
    \end{pgfonlayer}
  \end{scope}
\end{tikzpicture}
}
\]
This denotes that the place in the top left is used as the first and third input of the transition, while the place in the bottom left is the second input. 

Pre-nets give rise to strict monoidal categories as follows.  

%
%
\begin{prop}\label{prop: adjunction PreNet - StrMC}
  There is an adjunction
  \[
    \begin{tikzcd}
      \PreNet \ar[r,bend left,"F_\PreNet"]  \ar[r, phantom, "\bot"] & \StrMC \ar[l,bend left,"U_\PreNet"] 
    \end{tikzcd}
  \]
  whose left adjoint sends a pre-net $Q$ to the strict monoidal category where:
  \begin{itemize}
    \item Objects are elements of the free monoid on the set of places.
    \item   Morphisms are generated inductively by the following rules:
      \begin{itemize}
      \item for each place $s$ there is an identity $1_s \maps s \to s$
      \item for each transition $\tau$ of $Q$, there is a morphism going from its source to its target
      \item for every pair of morphisms $f \maps x \to y$ and $f' \maps x' \to y'$, there is a morphism $f \otimes f' \maps x \otimes x' \to y \otimes y'$
      \item for every pair of composable morphisms $f \maps x \to y$ and $g \maps y \to z$, there is a morphism $g \circ f \maps x \to z$
    \end{itemize}
    and quotiented to satisfy the axioms of a strict monoidal category.
  \end{itemize}
\end{prop}
\begin{proof}
This is \cite[Prop.\ 6.1]{Master}.
\end{proof}
The above adjunction can be composed with one defined in \cref{prop: free symmetries} to obtain an adjunction between pre-nets and strict symmetric monoidal categories:
 \[\begin{tikzcd} \PreNet \ar[r,bend left,"F_{\PreNet}"]  \ar[r, phantom, "\bot"]
 & \StrMC \ar[r,bend left,"F_{\StrMC}"] \ar[r, phantom, "\bot"] \ar[l,bend left,"U_{\PreNet}"]
 & \ar[l,bend left,"U_{\StrMC}"] \SSMC. \end{tikzcd} \]
 The composite adjunction
 \[\begin{tikzcd} \PreNet \ar[r,bend left,"F_{\ast}"]  \ar[r, phantom, "\bot"]
 & \ar[l,bend left,"U_{\ast}"] \SSMC \end{tikzcd} \]
 is used to obtain the categorical operational semantics of pre-nets under the individual token philosophy. This adjunction was first presented by Bruni, Meseguer, Montanari, and Sassone \cite{BMMS2001} with codomain a full subcategory of $\SSMC$ and was later refined to the above form in \cite{Master}.

 This composite adjunction has also been used to give a categorical semantics for Petri nets \cite{BMMS1999,BMMS2001,SassoneStrong}.  For this, given a Petri net $P$, one first chooses a pre-net $Q$ having $P$ as its underlying Petri net, and then forms the symmetric strict monoidal category $F_*(Q)$.  However this semantics is not functorial, due to the arbitrary choice involved.

The category $\PreNet$ is better behaved than $\Petri$.  The latter is not even cartesian closed, for essentially the same reasons described in~\cite{Cheng,MakkaiZawadowski}, but the former is cartesian closed, and even a topos:

\begin{prop}
  \label{prop:topos of pre-nets}
  The category $\PreNet$ is equivalent to a presheaf category.
\end{prop}
\begin{proof}
It suffices to construct a category $C$ so that functors from $C$ to $\Set$ are the same as pre-nets. Let $C$ have an object $p$, and for every pair of natural numbers $(n,m)$, let $C$ contain an object $t(m,n)$.  Here $p$ stands for `places', while $t(m,n)$ stands for `transitions with $m$ inputs and $n$ outputs'. Besides identity morphisms, $C$ contains $m$ morphisms $s_i \maps t(m,n) \to p$ representing the source maps and $n$ morphisms $t_j \maps t(m,n) \to p$ representing the target maps.  Composition in $C$ is trivial.

A pre-net $T \xrightarrow{(s,t)} S^* \times S^*$ can be identified with the functor $C \to \Set$ that sends $p$ to the set of places $S$, sends the object $t(m,n)$ to the subset of $T$ consisting of transitions with $m$ inputs and $n$ outputs, and sends the morphisms $s_i,t_j \maps t(m,n)\to p$ to the functions that map each transition to its $i$-th input and $j$-th output. A morphism of pre-nets $(f,g)$ can then be identified with a natural transformation between such functors, with the $p$-component given by $g$ and with the $t(m,n)$-components given by the restrictions of $f$ to the set of transitions with $m$ inputs and $n$ outputs. Naturality follows from the commutative diagrams in \cref{defn:prenet}.
\end{proof}
One downside of pre-nets is that ordering the inputs and outputs of transitions seems artificial in many applications where Petri nets are heavily used~\cite{GGHPPV,Statebox}. This ordering also greatly restricts the available morphisms between pre-nets. For example, there is no morphism between the following pre-nets:
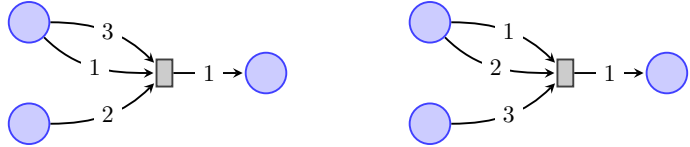
\[
  \scalebox{0.9}{  
\begin{tikzpicture}
  \begin{scope}
    \begin{pgfonlayer}{nodelayer}
        \node [place,tokens=0] (1a) at (-2,0.75) {};
        \node [place,tokens=0] (1b) at (-2,-0.75) {};      			
        \node [place,tokens=0] (3a) at (1.5,0) {};
        \node[transition] (2a) at (0,0) {};
    \end{pgfonlayer}
    \begin{pgfonlayer}{edgelayer}
        \draw[style=inarrow, out=0, in=135, thick] (1a) to node[fill=white, midway] {$3$} (2a);
        \draw[style=inarrow, out=-45, in=180, thick] (1a) to node[fill=white, midway] {$1$} (2a);
        \draw[style=inarrow, out=0, in=225, thick] (1b) to node[fill=white,midway] {$2$} (2a);
        \draw[style=inarrow, thick] (2a) to node[fill=white, midway] {$1$} (3a);
    \end{pgfonlayer}
  \end{scope}
\end{tikzpicture}
} \qquad \iflics\else\qquad\fi \scalebox{0.9}{  
\begin{tikzpicture}
  \begin{scope}
    \begin{pgfonlayer}{nodelayer}
        \node [place,tokens=0] (1a) at (-2,0.75) {};
        \node [place,tokens=0] (1b) at (-2,-0.75) {};      			
        \node [place,tokens=0] (3a) at (1.5,0) {};
        \node[transition] (2a) at (0,0) {};
    \end{pgfonlayer}
    \begin{pgfonlayer}{edgelayer}
        \draw[style=inarrow, out=0, in=135, thick] (1a) to node[fill=white, midway] {$1$} (2a);
        \draw[style=inarrow, out=-45, in=180, thick] (1a) to node[fill=white, midway] {$2$} (2a);
        \draw[style=inarrow, out=0, in=225, thick] (1b) to node[fill=white,midway] {$3$} (2a);
        \draw[style=inarrow, thick] (2a) to node[fill=white, midway] {$1$} (3a);
    \end{pgfonlayer}
  \end{scope}
\end{tikzpicture}
}
\]
even though there is a morphism between their underlying Petri nets, in which the ordering information has been forgotten.

Furthermore, in the symmetric strict monoidal category $F_\ast(Q)$ coming from a pre-net $Q$, \emph{none} of the symmetries $\sigma_{x,y} \maps x \otimes y \to y \otimes x$ are identities, except when $x$ or $y$ is the unit object.   Perhaps more importantly, if a transition in the pre-net $Q$ gives a morphism
 \[    t \maps x_1 \otimes \cdots \otimes x_m \to y_1 \otimes \cdots \otimes y_n ,\]
 where $x_i,y_j$ are objects coming from places of $Q$, the composite of $t$ with symmetries that permute the inputs $x_i$ and outputs $y_j$ is only equal to $f$ if both these permutations are the identity.
 
 Thus, the SSMCs obtained from pre-nets exemplify an extreme version of the individual token philosophy.  Not only does each token have its own individual identity, switching two tokens before or after executing the morphism corresponding to a transition \emph{always} gives a different morphism.

\section{\texorpdfstring{$\SigmaNet$s}{Sigma-nets}}
\label{sec:sigma nets}
We have seen that Petri nets generate symmetric monoidal categories naturally suited to the collective token philosophy, where the tokens have no individual identity, so it makes no sense to speak of switching them.  On the other hand, we have just seen that pre-nets generate symmmetric  monoidal categories suited to an extreme version of the individual token philosophy, in which switching tokens \emph{always} has an effect.

Now we introduce a new kind of nets, called $\Sigma$-nets, which in some sense lie between these two extremes.  In a \SigmaNet{}, one has control over which permutations of the input or output of a transition alter the morphism it defines and which do not.  This finer ability to control the action of permutations allow $\Sigma$-nets to behave either like Petri nets or pre-nets---or in a mixed way.
\begin{lem} 
  There is a forgetful functor
  \[Q \colon \SSMC \to \Set  \]
  that sends a symmetric strict monoidal category to its set of objects and sends a strict symmetric monoidal functor to its underlying function on objects.
  $Q$ has a left adjoint
  \[P \colon \Set \to \SSMC \]
  that sends a set $S$ to the symmetric strict monoidal category $PS$ having (possibly empty) words in $S$ as objects, and permutations as morphisms.  
\end{lem}

\begin{proof} 
See Sassone \cite[Sec.\ 3]{SassoneStrong} or Gambino and Joyal \cite[Sec.\ 3.1]{GJ}.
\end{proof}
\begin{defn}\label{def: sigma net}
  A \SigmaNet{} is a set $S$ together with a functor
  \[
    N \colon PS \times PS\op \to \Set.
  \]
  A morphism between \SigmaNet{}s $PS_1 \times PS_1\op \xrightarrow{N} \Set$ and $PS_2 \times PS_2\op \xrightarrow{M} \Set$ is a pair $(g,\alpha)$ where $g \maps S_1 \to S_2$ is a function and $\alpha$ is a natural transformation filling the following diagram:
  \[
    \begin{tikzcd}
      PS_1 \times PS_1\op 
      \ar[dd, "Pg \times Pg\op", swap] 
      \ar[dr, " "{name=U, pos=0.2}, "N"]
      &\\&  
      \Set.
      \\
      PS_2\times PS_2\op
      \ar[ur, " "{name=D, pos=0.15}, "M", swap]
      & 
      \arrow[Rightarrow, from = U, to = D, shorten <=3ex, shorten >=3ex, "\alpha", swap]
    \end{tikzcd}
  \]
  This defines the category $\SigmaNetCat$.
\end{defn}
The definition of \SigmaNet{} may seem unintuitive, but it is easily explained.  Suppose $N \maps PS \times PS\op \to \Set$ is a \SigmaNet.  We call $S$ its set of \defin{places}.  Objects of $PS$ are words of places.   Given $m,m' \in PS$, we call an element of $N(m,m')$ a \defin{transition} with \defin{source} $m$ and \defin{target} $m'$.   In \cref{thm:SigmaNet_SSMC_adjunction} we describe how to freely generate a symmetric strict monoidal category $C$ from the \SigmaNet $N$.  In this construction, the transitions of $N$ give morphisms that generate all the morphisms in $C$.

More precisely: the objects of $C$ are words of places.  The tensor product of objects is given by concatenation of words, and the symmetry in $C$ acts by permuting places in a word.  Each transition $t \in N(m,m')$ gives a morphism $t \maps m \to m'$ in $C$, and all other morphisms are generated by composition, tensor product and symmetries.

Given a transition $t \in N(m,m')$, the action of the functor $N \maps PS \times PS\op \to \Set$ on morphisms describes what happens to the corresponding morphism $t \maps m \to m'$ when we permute the places in its source and target: it gets sent to some other transition (possibly the same one).
We say two transitions $t \in N(m,m')$ and $u \in N(n,n')$ are in the same \defin{transition class} if and only if there exists a morphism $\sigma \maps (m,m') \to (n,n')$ in $PS \times PS\op$ such that 
\[    N(\sigma)(t) = u. \]
\begin{ex}\label{ex: no dup sigma net}
  There is a $\SigmaNet$ $N$ with just two places, say $p$ and $q$, and just two transitions, $t_1 \in N(pq,\epsilon)$ and $t_2 \in N(qp,\epsilon)$, where $\epsilon$ stands for the empty word. There are two morphisms in $PS\times PS\op$ with domain $(pq,\epsilon)$, namely the identity and the swap $(pq,\epsilon) \to (qp, \epsilon)$.
  Since there is a unique function between any two singleton sets, $N$ of this swap must map $t_1$ to $t_2$.  Thus, both $t_1$ and $t_2$ lie in the same transition class, and this $\SigmaNet$ has just one transition class.
\end{ex}
\begin{ex}\label{ex: no dup two trans sigma net}
  Now consider a \SigmaNet $M$ with just two places $p$ and $q$ and exactly four transitions, with $M(pq,\epsilon) = \{t_1,u_1\}$ and $M(qp,\epsilon) = \{t_2,u_2\}$.  We set $M$ of the swap $(pq,\epsilon) \to (qp,\epsilon)$ to act as the function $\{t_1,u_1\} \to \{t_2,u_2\}$ sending $t_1$ to $t_2$ and $u_1$ to $u_2$. This \SigmaNet has exactly two transition classes: $t_1$ and $t_2$ represent one transition class, and $u_1$ and $u_2$ represent the other.  Note that if we instead define the action of $M$ on the swap to send $t_1$ to $u_2$ and $t_2$ to $u_1$, then we would still have two transition classes; in fact this would be an isomorphic \SigmaNet.
\end{ex}
\cref{ex: no dup sigma net,ex: no dup two trans sigma net} illustrate the situation when no place occurs more than once in the source or target of any transition.  In this case, any two values of the functor $N$ will be related by \emph{at most one} morphism in $PS\times PS\op$, and the functorial action of $N$ on such a morphism provides a way to canonically identify their values.
\begin{ex}\label{ex: petri style sigma net}
  Next consider a $\Sigma$-net $O$ with one place $p$ and one transition, namely $t \in O(pp,\epsilon)$.    There are still two morphisms in $PS\times PS\op$ with domain $(pp,\epsilon)$, the identity and the swap, but now both have $(pp,\epsilon)$ as codomain as well.  There is still only one transition class, but now $t$ is mapped to itself by \emph{both} morphisms $(pp,\epsilon) \to (pp,\epsilon)$, the identity and the swap.
\end{ex}
Given a group $G$ acting on a set $X$, the \defin{isotropy group} of $x \in X$ is the subgroup of $G$ consisting of elements that map $x$ to itself.
Thus, in \cref{ex: petri style sigma net}, unlike \cref{ex: no dup two trans sigma net}, we are seeing a transition with a nontrivial isotropy group.  In fact, because permutations act \emph{trivially} on all transitions, the \SigmaNet of \cref{ex: petri style sigma net} belongs to the image of $\Petri$ under the functor $G_{\pet} \maps \Petri \to \SigmaNetCat$ described in \cref{thm:petri-net-sigma-net-adjunction}.
\begin{ex}\label{ex: prenet style sigma net}
   Next consider a \SigmaNet $Q$ with one place $p$ and precisely two transitions with $Q(pp,\epsilon) = \{t_1,t_2\}$. Suppose that $Q$ of the identity $(pp,\epsilon) \to (pp,\epsilon)$ acts as the identity function (as it must), while $Q$ of the swap acts by $t_1 \mapsto t_2$ and $t_2 \mapsto t_1$.  Then $t_1$ and $t_2$ represent the same transition class, so there is once again only one transition class.  The isotropy groups of $t_1$ and $t_2$ are trivial.   In fact, because permutations act \emph{freely} on the transitions in every transition class, this \SigmaNet belongs to the image of $\PreNet$ under the functor $F_\pre \maps \PreNet \to \SigmaNetCat$ described in \cref{thm:pre-net - sigma-net adjunction}.
\end{ex}
\begin{ex}\label{ex: mixed sigma net}
   Now let us give an example blending features from Examples~\ref{ex: petri style sigma net} and~\ref{ex: prenet style sigma net}.   For this, we create a \SigmaNet $R$ that has one place $p$ and three transitions $t_1,t_2,u \in R(pp,\epsilon)$, such that $t_1$ and $t_2$ are order-sensitive while $u$ is not.  This \SigmaNet maps $(pp,\epsilon)$ to $\{ t_1,t_2,u \}$ and everything else to the empty set.  The action of the swap automorphism of $(pp,\epsilon)$ switches $t_1, t_2$ and fixes $u$. As a result, this \SigmaNet has two transition classes: $t_1, t_2$ are both representatives of one transition class, while $u$  represents the other.    This \SigmaNet is not in the image of $G_{\pet} \maps \Petri \to \SigmaNetCat$ or $F_\pre \maps \PreNet \to \SigmaNetCat$; it mixes the two worlds.
\end{ex}
We now introduce some graphical representations for $\Sigma$-nets.
The first depicts a transition class as a three-dimensional tank containing pictures of the permutations that act trivially on an arbitrarily chosen transition in this class; for instance:
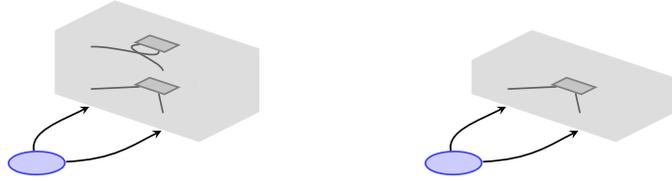
\[
    \scalebox{0.8}{
\begin{tikzpicture}
  \begin{scope}[
      yshift=0,every node/.append style={
      yslant=\ycoord,xslant=\xcoord},yslant=\ycoord,xslant=\xcoord
      ]
    \begin{pgfonlayer}{nodelayer}
      \node [place,tokens=0] (1a) at (-2,0) {};
    \end{pgfonlayer}
    \begin{pgfonlayer}{edgelayer}
      \draw[style=inarrow, thick, out=45, in=180] (1a) to (-0.5,0.5);
      \draw[style=inarrow, thick, out=-45, in=180] (1a) to  (-0.5,-0.5);
    \end{pgfonlayer}
    
      \fill[\transitionbgrd, opacity=\transitionop]  (-0.5,-1) rectangle (0.5,1);
      \node (11) at (-0.5,1) {};
      \node (21) at (-0.5,-1) {};
      \node (31) at (0.5,-1) {};
      \node (41) at (0.5,1) {};
  \end{scope}

  \begin{scope}[
      yshift=8,every node/.append style={
      yslant=\ycoord,xslant=\xcoord},yslant=\ycoord,xslant=\xcoord
      ]
    \begin{pgfonlayer}{nodelayer}
      \node [transition] (2b) []  at (0,0) {};  
    \end{pgfonlayer}
    \begin{pgfonlayer}{edgelayer}
      \draw[thick, in=135, out=-45] (-0.5,0.5) to  (2b);
      \draw[thick,  in=225, out=45] (-0.5,-0.5) to  (2b);
    \end{pgfonlayer}
  \end{scope}

  \begin{scope}[
      yshift=28,every node/.append style={
      yslant=\ycoord,xslant=\xcoord},yslant=\ycoord,xslant=\xcoord
      ]
    \begin{pgfonlayer}{nodelayer}
      \node [transition] (2b) []  at (0,0) {};  
    \end{pgfonlayer} 
    \begin{pgfonlayer}{edgelayer}
      \draw[thick, in=225, out=-45] (-0.5,0.5) to  (2b);
      \draw[thick,  in=135, out=45] (-0.5,-0.5) to  (2b);
    \end{pgfonlayer}
  \end{scope}

  \begin{scope}[
      yshift=30,every node/.append style={
      yslant=\ycoord,xslant=\xcoord},yslant=\ycoord,xslant=\xcoord
      ]

      \fill[\transitionbgrd, opacity=\transitionop]  (-0.5,-1) rectangle (0.5,1);
      \node (14) at (-0.5,1) {};
      \node (24) at (-0.5,-1) {};
      \node (34) at (0.5,-1) {};
      \node (44) at (0.5,1) {};
  \end{scope}

    \fill[\transitionbgrd, opacity=\transitionop]   (11.center) -- (14.center) -- (24.center) -- (21.center) -- (11.center);
    \fill[\transitionbgrd, opacity=\transitionop]   (21.center) -- (24.center) -- (34.center) -- (31.center) -- (21.center);
    \fill[\transitionbgrd, opacity=\transitionop]   (31.center) -- (34.center) -- (44.center) -- (41.center) -- (31.center);
    \fill[\transitionbgrd, opacity=\transitionop]   (41.center) -- (44.center) -- (14.center) -- (11.center) -- (41.center);

\end{tikzpicture}
} \iflics\hspace{3em}\else\hspace{5em}\fi \scalebox{0.8}{
\begin{tikzpicture}
  \begin{scope}[
      yshift=0,every node/.append style={
      yslant=\ycoord,xslant=\xcoord},yslant=\ycoord,xslant=\xcoord
      ]
    \begin{pgfonlayer}{nodelayer}
      \node [place,tokens=0] (1a) at (-2,0) {};
    \end{pgfonlayer}
    \begin{pgfonlayer}{edgelayer}
      \draw[style=inarrow, thick, out=45, in=180] (1a) to (-0.5,0.5);
      \draw[style=inarrow, thick, out=-45, in=180] (1a) to  (-0.5,-0.5);
    \end{pgfonlayer}
    
      \fill[\transitionbgrd, opacity=\transitionop]  (-0.5,-1) rectangle (0.5,1);
      \node (11) at (-0.5,1) {};
      \node (21) at (-0.5,-1) {};
      \node (31) at (0.5,-1) {};
      \node (41) at (0.5,1) {};
  \end{scope}

  \begin{scope}[
      yshift=8,every node/.append style={
      yslant=\ycoord,xslant=\xcoord},yslant=\ycoord,xslant=\xcoord
      ]
    \begin{pgfonlayer}{nodelayer}
      \node [transition] (2b) []  at (0,0) {};  
    \end{pgfonlayer}
    \begin{pgfonlayer}{edgelayer}
      \draw[thick, in=135, out=-45] (-0.5,0.5) to  (2b);
      \draw[thick,  in=225, out=45] (-0.5,-0.5) to  (2b);
    \end{pgfonlayer}
  \end{scope}


  \begin{scope}[
      yshift=16,every node/.append style={
      yslant=\ycoord,xslant=\xcoord},yslant=\ycoord,xslant=\xcoord
      ]

      \fill[\transitionbgrd, opacity=\transitionop]  (-0.5,-1) rectangle (0.5,1);
      \node (14) at (-0.5,1) {};
      \node (24) at (-0.5,-1) {};
      \node (34) at (0.5,-1) {};
      \node (44) at (0.5,1) {};
  \end{scope}

    \fill[\transitionbgrd, opacity=\transitionop]   (11.center) -- (14.center) -- (24.center) -- (21.center) -- (11.center);
    \fill[\transitionbgrd, opacity=\transitionop]   (21.center) -- (24.center) -- (34.center) -- (31.center) -- (21.center);
    \fill[\transitionbgrd, opacity=\transitionop]   (31.center) -- (34.center) -- (44.center) -- (41.center) -- (31.center);
    \fill[\transitionbgrd, opacity=\transitionop]   (41.center) -- (44.center) -- (14.center) -- (11.center) -- (41.center);

\end{tikzpicture}
}
\]
On the left is Example~\ref{ex: petri style sigma net}, where the identity and the swap both act trivially.  On the right is Example~\ref{ex: prenet style sigma net}, where only the identity acts trivially.  \cref{ex: no dup sigma net,ex: mixed sigma net} are instead:
\[
  \scalebox{\iflics 0.6\else 0.75\fi}{
\begin{tikzpicture}
  \begin{scope}[
      yshift=0,every node/.append style={
      yslant=\ycoord,xslant=\xcoord},yslant=\ycoord,xslant=\xcoord
      ]
    \begin{pgfonlayer}{nodelayer}
      \node [place,tokens=0] (1a) at (-2,0.75) {};
      \node [place,tokens=0] (1b) at (-2,-0.75) {};      			
    \end{pgfonlayer}
    \begin{pgfonlayer}{edgelayer}
      \draw[style=inarrow, thick] (1a) to  (-0.5,0.75);
      \draw[style=inarrow, thick] (1b) to  (-0.5,-0.75);
    \end{pgfonlayer}
    
      \fill[\transitionbgrd, opacity=\transitionop]  (-0.5,-1) rectangle (0.5,1);
      \node (11) at (-0.5,1) {};
      \node (21) at (-0.5,-1) {};
      \node (31) at (0.5,-1) {};
      \node (41) at (0.5,1) {};
  \end{scope}

  \begin{scope}[
      yshift=8,every node/.append style={
      yslant=\ycoord,xslant=\xcoord},yslant=\ycoord,xslant=\xcoord
      ]
    \begin{pgfonlayer}{nodelayer}
      \node [transition] (2b) []  at (0,0) {};  
    \end{pgfonlayer}
    \begin{pgfonlayer}{edgelayer}
      \draw[thick, in=135, out=-45] (-0.5,0.5) to  (2b);
      \draw[thick,  in=225, out=45] (-0.5,-0.5) to  (2b);
    \end{pgfonlayer}
  \end{scope}


  \begin{scope}[
      yshift=16,every node/.append style={
      yslant=\ycoord,xslant=\xcoord},yslant=\ycoord,xslant=\xcoord
      ]

      \fill[\transitionbgrd, opacity=\transitionop]  (-0.5,-1) rectangle (0.5,1);
      \node (14) at (-0.5,1) {};
      \node (24) at (-0.5,-1) {};
      \node (34) at (0.5,-1) {};
      \node (44) at (0.5,1) {};
  \end{scope}

    \fill[\transitionbgrd, opacity=\transitionop]   (11.center) -- (14.center) -- (24.center) -- (21.center) -- (11.center);
    \fill[\transitionbgrd, opacity=\transitionop]   (21.center) -- (24.center) -- (34.center) -- (31.center) -- (21.center);
    \fill[\transitionbgrd, opacity=\transitionop]   (31.center) -- (34.center) -- (44.center) -- (41.center) -- (31.center);
    \fill[\transitionbgrd, opacity=\transitionop]   (41.center) -- (44.center) -- (14.center) -- (11.center) -- (41.center);

\end{tikzpicture}
} \iflics\hspace{2em}\else\hspace{5em}\fi \scalebox{\iflics 0.6\else 0.75\fi}{
\begin{tikzpicture}
  \begin{scope}[
      yshift=0,every node/.append style={
      yslant=\ycoord,xslant=\xcoord},yslant=\ycoord,xslant=\xcoord
      ]
    \begin{pgfonlayer}{nodelayer}
      \node [place,tokens=0] (1a) at (-2,1.5) {};
    \end{pgfonlayer}
    \begin{pgfonlayer}{edgelayer}
      \draw[style=inarrow, thick, out=-30, in=180] (1a) to (-0.5,0.5);
      \draw[style=inarrow, thick, out=-60, in=180] (1a) to  (-0.5,-0.5);
      \draw[style=inarrow, thick, out=60, in=180] (1a) to (-0.5,3.5);
      \draw[style=inarrow, thick, out=30, in=180] (1a) to  (-0.5,2.5);
    \end{pgfonlayer}

      \fill[\transitionbgrd, opacity=\transitionop]  (-0.5,2) rectangle (0.5,4);
      \node (13) at (-0.5,4) {};
      \node (23) at (-0.5,2) {};
      \node (33) at (0.5,2) {};
      \node (43) at (0.5,4) {};
      
      \fill[\transitionbgrd, opacity=\transitionop]  (-0.5,-1) rectangle (0.5,1);
      \node (11) at (-0.5,1) {};
      \node (21) at (-0.5,-1) {};
      \node (31) at (0.5,-1) {};
      \node (41) at (0.5,1) {};
  \end{scope}

  \begin{scope}[
      yshift=8,every node/.append style={
      yslant=\ycoord,xslant=\xcoord},yslant=\ycoord,xslant=\xcoord
      ]
      
    \begin{pgfonlayer}{nodelayer}
      \node [transition] (2a) []  at (0,3) {};  
      \node [transition] (2b) []  at (0,0) {};  
    \end{pgfonlayer}
    \begin{pgfonlayer}{edgelayer}
      \draw[thick, in=135, out=-45] (-0.5,0.5) to  (2b);
      \draw[thick,  in=225, out=45] (-0.5,-0.5) to  (2b);

      \draw[thick, in=135, out=-45] (-0.5,3.5) to  (2a);
      \draw[thick,  in=225, out=45] (-0.5,2.5) to  (2a);
    \end{pgfonlayer}
  \end{scope}
  
  \begin{scope}[
      yshift=28,every node/.append style={
      yslant=\ycoord,xslant=\xcoord},yslant=\ycoord,xslant=\xcoord
      ]
      
    \begin{pgfonlayer}{nodelayer}
      \node [transition] (2a) []  at (0,3) {};  
    \end{pgfonlayer}
    \begin{pgfonlayer}{edgelayer}

      \draw[thick, in=225, out=-45] (-0.5,3.5) to  (2a);
      \draw[thick,  in=135, out=45] (-0.5,2.5) to  (2a);
    \end{pgfonlayer}
  \end{scope}

  \begin{scope}[
      yshift=30,every node/.append style={
      yslant=\ycoord,xslant=\xcoord},yslant=\ycoord,xslant=\xcoord
      ]
      
      \fill[\transitionbgrd, opacity=\transitionop]  (-0.5,2) rectangle (0.5,4);

      \node (14) at (-0.5,4) {};
      \node (24) at (-0.5,2) {};
      \node (34) at (0.5,2) {};
      \node (44) at (0.5,4) {};
  \end{scope}
  
  \begin{scope}[
      yshift=10,every node/.append style={
      yslant=\ycoord,xslant=\xcoord},yslant=\ycoord,xslant=\xcoord
      ]

      \fill[\transitionbgrd, opacity=\transitionop]  (-0.5,-1) rectangle (0.5,1);
      \node (12) at (-0.5,1) {};
      \node (22) at (-0.5,-1) {};
      \node (32) at (0.5,-1) {};
      \node (42) at (0.5,1) {};
  \end{scope}

    \fill[\transitionbgrd, opacity=\transitionop]   (11.center) -- (12.center) -- (22.center) -- (21.center) -- (11.center);
    \fill[\transitionbgrd, opacity=\transitionop]   (21.center) -- (22.center) -- (32.center) -- (31.center) -- (21.center);
    \fill[\transitionbgrd, opacity=\transitionop]   (31.center) -- (32.center) -- (42.center) -- (41.center) -- (31.center);
    \fill[\transitionbgrd, opacity=\transitionop]   (41.center) -- (42.center) -- (12.center) -- (11.center) -- (41.center);

    \fill[\transitionbgrd, opacity=\transitionop]   (13.center) -- (14.center) -- (24.center) -- (23.center) -- (13.center);
    \fill[\transitionbgrd, opacity=\transitionop]   (23.center) -- (24.center) -- (34.center) -- (33.center) -- (23.center);
    \fill[\transitionbgrd, opacity=\transitionop]   (33.center) -- (34.center) -- (44.center) -- (43.center) -- (33.center);
    \fill[\transitionbgrd, opacity=\transitionop]   (43.center) -- (44.center) -- (14.center) -- (13.center) -- (43.center);

\end{tikzpicture}
}
\]
Although each transition class is a set of transitions, and each tank represents a single transition class, the pictures inside that tank do \emph{not} represent the transitions in that class, but rather the \emph{isotropy group} of a \emph{single} transition in the class.  For a given transition class, with $m$ inputs and $n$ outputs, say, the 
number of transitions and the size of the isotropy group are inversely related: their product is the cardinality $m!n!$ of the total symmetry group $S_m\times S_n$.

Of course, this formalism quickly becomes disadvantageous for large nets.  Another approach is to draw a \SigmaNet{} using the usual two-dimensional representation of a Petri net, with one node for each transition class, but decorated by the relevant isotropy group. When this group is trivial, we can omit it.
In this style of drawing, the nets of \cref{ex: petri style sigma net,ex: prenet style sigma net} look as follows:
\[
  \scalebox{0.9}{  
\begin{tikzpicture}
  \begin{scope}
    \begin{pgfonlayer}{nodelayer}
        \node [place,tokens=0] (1a) at (-2,0) {};
        \node[transition, label=above:$S_2$] (2a) at (0,0) {};
    \end{pgfonlayer}
    \begin{pgfonlayer}{edgelayer}
      \draw[style=inarrow, thick, bend left] (1a) to (2a);
      \draw[style=inarrow, thick, bend right] (1a) to  (2a);
      \end{pgfonlayer}
  \end{scope}
\end{tikzpicture}
} \hspace{5em} \scalebox{0.9}{  
\begin{tikzpicture}
  \begin{scope}
    \begin{pgfonlayer}{nodelayer}
        \node [place,tokens=0] (1a) at (-2,0) {};
        \node[transition] (2a) at (0,0) {};
    \end{pgfonlayer}
    \begin{pgfonlayer}{edgelayer}
      \draw[style=inarrow, thick, bend left] (1a) to (2a);
      \draw[style=inarrow, thick, bend right] (1a) to  (2a);
    \end{pgfonlayer}
  \end{scope}
\end{tikzpicture}
}
\]
In the rest of the paper, we will use the three-dimensional representation.
However, both representations are potentially misleading, in different ways.  The two-dimensional representation does not show how the isotropy group sits inside the ambient group $S_m\times S_n$.  The three-dimensional representation indicates this, but it is sensitive to the arbitrary choice of one transition in each class.
For instance, the following three pictures all represent the same \SigmaNet{}, with one transition class containing three transitions, each with isotropy group isomorphic to $S_2$.
But these three copies of $S_2$ sit inside $S_3$ differently, yielding different three-dimensional pictures.
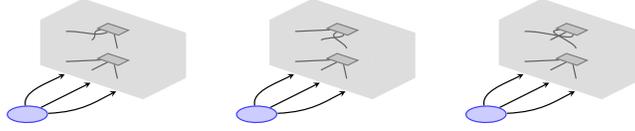
\[
  \scalebox{0.57}{
\begin{tikzpicture}
  \begin{scope}[
      yshift=0,every node/.append style={
      yslant=\ycoord,xslant=\xcoord},yslant=\ycoord,xslant=\xcoord
      ]
    \begin{pgfonlayer}{nodelayer}
      \node [place,tokens=0] (1a) at (-2,0) {};
    \end{pgfonlayer}
    \begin{pgfonlayer}{edgelayer}
      \draw[style=inarrow, thick, out=45, in=180] (1a) to (-0.5,0.5);
      \draw[style=inarrow, thick] (1a) to  (-0.5,0);
      \draw[style=inarrow, thick, out=-45, in=180] (1a) to  (-0.5,-0.5);
    \end{pgfonlayer}
    
      \fill[\transitionbgrd, opacity=\transitionop]  (-0.5,-1) rectangle (0.5,1);
      \node (11) at (-0.5,1) {};
      \node (21) at (-0.5,-1) {};
      \node (31) at (0.5,-1) {};
      \node (41) at (0.5,1) {};
  \end{scope}

  \begin{scope}[
      yshift=8,every node/.append style={
      yslant=\ycoord,xslant=\xcoord},yslant=\ycoord,xslant=\xcoord
      ]
    \begin{pgfonlayer}{nodelayer}
      \node [transition] (2b) []  at (0,0) {};  
    \end{pgfonlayer}
    \begin{pgfonlayer}{edgelayer}
      \draw[thick, in=135, out=-45] (-0.5,0.5) to  (2b);
      \draw[thick] (-0.5,0) to  (2b);
      \draw[thick,  in=225, out=45] (-0.5,-0.5) to  (2b);
    \end{pgfonlayer}
  \end{scope}

  \begin{scope}[
      yshift=28,every node/.append style={
      yslant=\ycoord,xslant=\xcoord},yslant=\ycoord,xslant=\xcoord
      ]
    \begin{pgfonlayer}{nodelayer}
      \node [transition] (2b) []  at (0,0) {};  
    \end{pgfonlayer} 
    \begin{pgfonlayer}{edgelayer}
      \draw[thick, in=180, out=-45] (-0.5,0.5) to  (2b);
      \draw[thick, out=0, in=135] (-0.5,0) to  (2b);
      \draw[thick,  in=225, out=45] (-0.5,-0.5) to  (2b);
    \end{pgfonlayer}
  \end{scope}

  \begin{scope}[
      yshift=30,every node/.append style={
      yslant=\ycoord,xslant=\xcoord},yslant=\ycoord,xslant=\xcoord
      ]

      \fill[\transitionbgrd, opacity=\transitionop]  (-0.5,-1) rectangle (0.5,1);
      \node (14) at (-0.5,1) {};
      \node (24) at (-0.5,-1) {};
      \node (34) at (0.5,-1) {};
      \node (44) at (0.5,1) {};
  \end{scope}

    \fill[\transitionbgrd, opacity=\transitionop]   (11.center) -- (14.center) -- (24.center) -- (21.center) -- (11.center);
    \fill[\transitionbgrd, opacity=\transitionop]   (21.center) -- (24.center) -- (34.center) -- (31.center) -- (21.center);
    \fill[\transitionbgrd, opacity=\transitionop]   (31.center) -- (34.center) -- (44.center) -- (41.center) -- (31.center);
    \fill[\transitionbgrd, opacity=\transitionop]   (41.center) -- (44.center) -- (14.center) -- (11.center) -- (41.center);

\end{tikzpicture}
} \quad \scalebox{0.57}{
\begin{tikzpicture}
  \begin{scope}[
      yshift=0,every node/.append style={
      yslant=\ycoord,xslant=\xcoord},yslant=\ycoord,xslant=\xcoord
      ]
    \begin{pgfonlayer}{nodelayer}
      \node [place,tokens=0] (1a) at (-2,0) {};
    \end{pgfonlayer}
    \begin{pgfonlayer}{edgelayer}
      \draw[style=inarrow, thick, out=45, in=180] (1a) to (-0.5,0.5);
      \draw[style=inarrow, thick] (1a) to  (-0.5,0);
      \draw[style=inarrow, thick, out=-45, in=180] (1a) to  (-0.5,-0.5);
    \end{pgfonlayer}
    
      \fill[\transitionbgrd, opacity=\transitionop]  (-0.5,-1) rectangle (0.5,1);
      \node (11) at (-0.5,1) {};
      \node (21) at (-0.5,-1) {};
      \node (31) at (0.5,-1) {};
      \node (41) at (0.5,1) {};
  \end{scope}

  \begin{scope}[
      yshift=8,every node/.append style={
      yslant=\ycoord,xslant=\xcoord},yslant=\ycoord,xslant=\xcoord
      ]
    \begin{pgfonlayer}{nodelayer}
      \node [transition] (2b) []  at (0,0) {};  
    \end{pgfonlayer}
    \begin{pgfonlayer}{edgelayer}
      \draw[thick, in=135, out=-45] (-0.5,0.5) to  (2b);
      \draw[thick] (-0.5,0) to  (2b);
      \draw[thick,  in=225, out=45] (-0.5,-0.5) to  (2b);
    \end{pgfonlayer}
  \end{scope}

  \begin{scope}[
      yshift=28,every node/.append style={
      yslant=\ycoord,xslant=\xcoord},yslant=\ycoord,xslant=\xcoord
      ]
    \begin{pgfonlayer}{nodelayer}
      \node [transition] (2b) []  at (0,0) {};  
    \end{pgfonlayer} 
    \begin{pgfonlayer}{edgelayer}
      \draw[thick, in=135, out=-45] (-0.5,0.5) to  (2b);
      \draw[thick, out=0, in=225] (-0.5,0) to  (2b);
      \draw[thick,  in=180, out=45] (-0.5,-0.5) to  (2b);
    \end{pgfonlayer}
  \end{scope}

  \begin{scope}[
      yshift=30,every node/.append style={
      yslant=\ycoord,xslant=\xcoord},yslant=\ycoord,xslant=\xcoord
      ]

      \fill[\transitionbgrd, opacity=\transitionop]  (-0.5,-1) rectangle (0.5,1);
      \node (14) at (-0.5,1) {};
      \node (24) at (-0.5,-1) {};
      \node (34) at (0.5,-1) {};
      \node (44) at (0.5,1) {};
  \end{scope}

    \fill[\transitionbgrd, opacity=\transitionop]   (11.center) -- (14.center) -- (24.center) -- (21.center) -- (11.center);
    \fill[\transitionbgrd, opacity=\transitionop]   (21.center) -- (24.center) -- (34.center) -- (31.center) -- (21.center);
    \fill[\transitionbgrd, opacity=\transitionop]   (31.center) -- (34.center) -- (44.center) -- (41.center) -- (31.center);
    \fill[\transitionbgrd, opacity=\transitionop]   (41.center) -- (44.center) -- (14.center) -- (11.center) -- (41.center);

\end{tikzpicture}
} \quad \scalebox{0.57}{
\begin{tikzpicture}
  \begin{scope}[
      yshift=0,every node/.append style={
      yslant=\ycoord,xslant=\xcoord},yslant=\ycoord,xslant=\xcoord
      ]
    \begin{pgfonlayer}{nodelayer}
      \node [place,tokens=0] (1a) at (-2,0) {};
    \end{pgfonlayer}
    \begin{pgfonlayer}{edgelayer}
      \draw[style=inarrow, thick, out=45, in=180] (1a) to (-0.5,0.5);
      \draw[style=inarrow, thick] (1a) to  (-0.5,0);
      \draw[style=inarrow, thick, out=-45, in=180] (1a) to  (-0.5,-0.5);
    \end{pgfonlayer}
    
      \fill[\transitionbgrd, opacity=\transitionop]  (-0.5,-1) rectangle (0.5,1);
      \node (11) at (-0.5,1) {};
      \node (21) at (-0.5,-1) {};
      \node (31) at (0.5,-1) {};
      \node (41) at (0.5,1) {};
  \end{scope}

  \begin{scope}[
      yshift=8,every node/.append style={
      yslant=\ycoord,xslant=\xcoord},yslant=\ycoord,xslant=\xcoord
      ]
    \begin{pgfonlayer}{nodelayer}
      \node [transition] (2b) []  at (0,0) {};  
    \end{pgfonlayer}
    \begin{pgfonlayer}{edgelayer}
      \draw[thick, in=135, out=-45] (-0.5,0.5) to  (2b);
      \draw[thick] (-0.5,0) to  (2b);
      \draw[thick,  in=225, out=45] (-0.5,-0.5) to  (2b);
    \end{pgfonlayer}
  \end{scope}

  \begin{scope}[
      yshift=28,every node/.append style={
      yslant=\ycoord,xslant=\xcoord},yslant=\ycoord,xslant=\xcoord
      ]
    \begin{pgfonlayer}{nodelayer}
      \node [transition] (2b) []  at (0,0) {};  
    \end{pgfonlayer} 
    \begin{pgfonlayer}{edgelayer}
      \draw[thick, in=225, out=-45] (-0.5,0.5) to  (2b);
      \draw[thick] (-0.5,0) to  (2b);
      \draw[thick,  in=135, out=45] (-0.5,-0.5) to  (2b);
    \end{pgfonlayer}
  \end{scope}

  \begin{scope}[
      yshift=30,every node/.append style={
      yslant=\ycoord,xslant=\xcoord},yslant=\ycoord,xslant=\xcoord
      ]

      \fill[\transitionbgrd, opacity=\transitionop]  (-0.5,-1) rectangle (0.5,1);
      \node (14) at (-0.5,1) {};
      \node (24) at (-0.5,-1) {};
      \node (34) at (0.5,-1) {};
      \node (44) at (0.5,1) {};
  \end{scope}

    \fill[\transitionbgrd, opacity=\transitionop]   (11.center) -- (14.center) -- (24.center) -- (21.center) -- (11.center);
    \fill[\transitionbgrd, opacity=\transitionop]   (21.center) -- (24.center) -- (34.center) -- (31.center) -- (21.center);
    \fill[\transitionbgrd, opacity=\transitionop]   (31.center) -- (34.center) -- (44.center) -- (41.center) -- (31.center);
    \fill[\transitionbgrd, opacity=\transitionop]   (41.center) -- (44.center) -- (14.center) -- (11.center) -- (41.center);

\end{tikzpicture}
}
\]
We do not know whether there is a graphical representation of \SigmaNet{}s that avoids both of these problems.

\section{Perspectives on the category \texorpdfstring{$\SigmaNet$}{Sigma-net}}
\label{sec:fib}

The definition of \SigmaNet{}s in \cref{sec:sigma nets} gives what may be called a ``profunctorial'' perspective: a \SigmaNet{} is a functor from
$PS \times PS\op$ to $\Set$, which is the same as a profunctor from $PS$ to itself.  This perspective will be useful in constructing the adjunction between $\SigmaNetCat$ and $\SSMC$ in \cref{thm:SigmaNet_SSMC_adjunction}.
However, there are other perspectives on \SigmaNet{}s, leading to two alternative descriptions of $\SigmaNetCat$, useful for other purposes.

\subsection{The presheaf perspective}
\label{sec:presheaf-perspective}

\begin{thm}\label{thm: sigma nets is a topos}
  \SigmaNetCat{} is equivalent to a presheaf category.
\end{thm}

\begin{proof} 
We construct a category $D$ so that functors from $D$ to $\Set$ can be identified with $\Sigma$-nets. To construct $D$, we take the category $C$ from \cref{prop:topos of pre-nets} and throw in extra automorphisms of each object $t(m,n)$, making its automorphism group $S_m \times S_n$.   For a source map $s_i \maps t(m,n) \to p$ and an automorphism $(\sigma,\tau) \in S_m \times S_n$,  we set the composite $s_i \circ (\sigma,\tau)$ equal to $s_{\sigma(i)}$.  Similarly, for a target map $t_j \maps t(m,n) \to p$, we set the composite $t_j \circ (\sigma,\tau)$ equal to $t_{\tau(j)}$.  Then, for each $\Sigma$-net $N \maps PS \times PS\op \to \Set$, there is a corresponding functor $\nu \maps D \to \Set$ defined as follows.  It sends the object $p \in D$ to the set of places of $N$.  It sends each object $t(m,n) \in D$ to the disjoint union of the sets $N(a,b)$ over all $a \in PS$ with length $m$ and $b \in PS$ with length $n$.  It sends the morphisms $s_i,t_j \maps t(m,n) \to p$ to the functions that map any transition to its $i$th input and $j$th output.  Finally, this functor $\nu$ sends the permutations $(\sigma,\tau)$ to the natural actions of the symmetric group on the transitions of $N$. For a morphism of $\Sigma$-nets $(g,\alpha) \maps N \to N'$, there is a natural transformation between their functors whose $p$-component is given by $g$ and whose $t(m,n)$-components are given by disjoint unions of the components of $\alpha$.  One can check that the resulting functor from $\SigmaNetCat{}$ to $\Set^D$ is an equivalence.
\end{proof}

\cref{thm: sigma nets is a topos} has a lot of nice consequences, such as:
\begin{itemize}
  \item $\SigmaNetCat$ is complete and cocomplete. This is particularly important since many compositional approaches to Petri nets rely on colimits; for example, composition of open Petri nets is done using pushouts\iflics\else\ (see \cref{sec:open-nets})\fi, while tensoring them is done using coproducts~\cite{BM}.
  \item $\SigmaNetCat$ is a topos, and thus an adhesive category \cite{LS2}, so it admits a theory of double pushout rewriting \cite{LS1}.   This is relevant as double pushout rewriting  is a widely used technique to transform graph-like structures in the literature \cite{Ehrig}.  The internal logic of toposes is very rich, and understanding its implications for $\Sigma$-nets is an interesting direction for future work.
\end{itemize}
Note that $\Petri$ is not a presheaf category, whereas the category of directed graphs is.  Indeed, as noted in~\cite{Kock2}, graphs are functors $C_{1,1}\to\Set$, where $C_{1,1}$ is the full subcategory of the above $C$ (or $D$) on the objects $p$ and $t(1,1)$.

The category $D$ in \cref{thm: sigma nets is a topos} is equivalent to the opposite of Kock's category of ``elementary graphs''~\cite[1.5.4]{Kock2}.  Thus, $\SigmaNetCat$ is equivalent to his category of ``digraphical species''~\cite[2.1]{Kock2}.  Similarly, $C$ is the opposite of Kock's ``elementary planar graphs''.
\subsection{The groupoidal perspective}
The profunctorial and presheaf perspectives highlight the transitions of a \SigmaNet{} over its transition classes.
Sometimes, however, we want to work directly with the transition classes; we now describe a third perspective that permits this.

Firstly, it is well-known \cite[Theorem 2.1.2]{LR} that a functor $N \maps PS\times PS\op\to \Set$ is equivalent to a discrete opfibration $T \to PS\times PS\op$ for some category $T$.  Since $PS\times PS\op$ is a groupoid, $T$ is as well.  In addition to this, a morphism of \SigmaNet{}s is equivalently a commutative square
\[
  \begin{tikzcd}[column sep=1.5cm]
    T_1 \ar[r,"g"] \ar[d] & T_2 \ar[d] \\
    PS_1 \times PS_1\op \ar[r,"Pf \times Pf\op"'] & PS_2\times PS_2\op
  \end{tikzcd}
\]
Note that this looks much more similar to the definitions of the categories $\PreNet$ and $\Petri$.
The set of objects of $T$ here is the disjoint union of the sets $N(p,p')$, i.e., the transitions rather than the transition classes.
The transition classes are the \emph{isomorphism classes} of the groupoid $T$.   To contract these down to single objects, we can replace $T$ by an equivalent groupoid that is \defin{skeletal}, i.e., there are no morphisms $x \to y$  for objects $x \neq y$, or equivalently each isomorphism class contains exactly one object.   After such a replacement the functor $T\to PS\times PS\op$ is no longer a discrete opfibration, but it is still faithful.
To compensate for this replacement of discrete opfibrations by faithful functors with skeletal domain, when defining the morphisms of \SigmaNet{}s we have to allow the squares to commute up to isomorphism rather than strictly.
\begin{thm}\label{thm: fibrational representation}
  The category of \SigmaNet{}s is equivalent to the following category:
  \begin{itemize}
  \item Its objects are faithful functors $T \to PS\times PS\op$, where $S$ is a set and $T$ is a skeletal groupoid.
  \item Its morphisms are squares that commute up to specified isomorphism
    \[
      \begin{tikzcd}[column sep=1.5cm]
        T_1 \ar[r,"g"] \ar[d] \ar[dr,phantom,near end,"\cong\scriptstyle\theta"] & T_2 \ar[d] \\
        PS_1 \times PS_1\op \ar[r,"Pf \times Pf\op"'] & PS_2\times PS_2\op
      \end{tikzcd}
    \]
    modulo the equivalence relation that two such morphisms $(f,g,\theta)$ and $(f',g',\theta')$ are considered equal if $f=f'$ and there is a natural isomorphism $\phi \colon g\simRightarrow g'$ such that
     \[ 
      \scalebox{\iflics 0.6\else 0.8\fi}{
      \begin{tikzpicture}
        \begin{pgfonlayer}{nodelayer}
          \node (1a) at (0,2) {$T_1$};
          \node (2a) at (4,2) {$T_2$};
          \node (3a) at (5.5,1) {$=$};
          \node (4a) at (7,2) {$T_1$};
          \node (5a) at (11,2) {$T_2$}; 
          \node(nat1) at (2,0.75) {$\cong_{\theta'}$};
          \node(nat1) at (9,1) {$\cong_{\theta}$};
          \node[rotate=-90, label=above:$\phi$](phi) at (2,2) {$\Rightarrow$};

          \node (1b) at (0,0) {$PS_1 \times PS_1\op$};
          \node (2b) at (4,0) {$PS_2 \times PS_2\op$};
          \node (4b) at (7,0) {$PS_1 \times PS_1\op$};
          \node (5b) at (11,0) {$PS_2 \times PS_2\op$};
        \end{pgfonlayer}
        \begin{pgfonlayer}{edgelayer}
          \draw[->](1a.north east) to node[above]{$g$} (2a.north west);
          \draw[->] (1a.south east) to node[below]{$g'$} (2a.south west);
          \draw[->] (4a) to node[above]{$g$} (5a);
          \draw[->](1a) to (1b);
          \draw[->](2a) to (2b);
          \draw[->](4a) to (4b);
          \draw[->](5a) to (5b);
          \draw[->] (1b) to node[below]{$Pf \times Pf\op$} (2b);
          \draw[->] (4b) to node[below]{$Pf \times Pf\op$} (5b);
        \end{pgfonlayer}
      \end{tikzpicture}}
    \]
  \end{itemize}
\end{thm}
Note that since $T_2 \to PS_2\times PS_2\op$ is faithful, such a $\phi$ is unique if it exists.
Thus, the category described in the theorem is in fact \emph{equivalent} to the evident 2-category having as 2-morphisms natural isomorphisms $\phi$ as above.
\begin{proof}
  Define a category \PN as in the theorem, but where the objects allow $T$ to be any groupoid.
  Then there is a functor $\SigmaNetCat\to \PN$, since discrete opfibrations are faithful and strictly commutative squares also commute up to isomorphism.  This functor is faithful, since if $\theta$ and $\theta'$ are identities so is $\phi$, by the faithfulness of $T_2 \to PS_2\times PS_2\op$.   
  Moreover, any morphism in \PN whose target $T_2 \to PS_2\times PS_2\op$ is a discrete opfibration has a representative that commutes strictly, since we can lift the isomorphism $\theta$ to an isomorphism $\phi$ with $\theta'$ an identity.  Thus, the functor $\SigmaNetCat\to \PN$ is also full.
  
  Let the \defin{pseudo slice 2-category} over a groupoid $B$ be the 2-category with groupoids over $B$ as objects, triangles commuting up to natural isomorphism
  \[
      \begin{tikzcd}[column sep=0.5cm,row sep=0.3cm]
        A \ar[ddr] \ar[rr,"g"]  & &  A' \ar[ddl] \\ 
         & \raisebox{1 em}{$\cong\scriptstyle\theta$}\\
         & B
      \end{tikzcd}
    \]
  as morphisms, and the evident 2-morphisms \cite[Definition 3.2]{Palmgren}.  Any groupoid over $B$ is equivalent, in the pseudo slice 2-category of $B$, to a fibration \cite[Theorem 6.7]{Strickland}, which in the groupoid case is the same as an opfibration.   If $f \maps A\to B$ is faithful then this opfibration will be as well, so $f$ is equivalent to a discrete opfibration.    Since equivalences in the pseudo slice 2-category yield isomorphisms in \PN, the functor $\SigmaNetCat\to \PN$ is also essentially surjective, and hence an equivalence.

  The category described in the theorem is a full subcategory of \PN, so it suffices to show that every object of \PN is isomorphic to one where $T$ is skeletal.  But any groupoid is equivalent to a skeletal one, and such an equivalence preserves faithfulness and yields an isomorphism in \PN.
\end{proof}
Note that the construction in the final paragraph taking a groupoid to a skeletal one preserves connected components, while in the output each connected component has exactly one object.  Thus, in the representation described in \cref{thm: fibrational representation} the objects of the groupoid $T$ really are precisely the transition classes.  Since the transition classes of a \SigmaNet{} correspond to the transitions of a Petri net, we can think of a \SigmaNet as a Petri net together with, for each transition, (1) a lifting of its source and target multisets to words, and (2) an \emph{isotropy group} that acts faithfully on those words, i.e., maps injectively to the subgroup of $S_m\times S_n$ that fixes both of those words.
\begin{ex}
  If we start from a transition in a Petri net $t \maps 3a+2b \to 4c$, then we could lift it to a transition in a \SigmaNet{} by defining $t \maps aaabb \to cccc$ and equipping it with any subgroup of $S_3 \times S_2 \times S_4$, which describes the ``degree of collectivization'' of $t$.
  If the isotropy group is trivial, then this particular transition behaves like one in a pre-net---tokens are ``fully individualized''---whereas if it is as large as possible then it behaves as a Petri net---tokens are ``fully collectivized''. This idea is heavily used in the next section to describe the adjunctions between $\Petri$, $\PreNet$ and $\SigmaNetCat$.
\end{ex} 
\section{Description of the adjunctions}
\label{sec:adjunctions}
Now we describe in detail all the adjunctions between the categories in play. We again include the diagram of \cref{sec: dramatis personae}, but now with most of the functors labeled.
\begin{equation}
\label{eq:adjunctions}
\begin{tikzcd}[row sep=large,column sep=large]
  \StrMC \ar[r,shift left, "F_{\StrMC}"] \ar[dd,shift left, "U_{\PreNet}"] 
  &
  \SSMC \ar[l,shift left, "U_{\StrMC}"] \ar[r,shift left, "F_{\SSMC}"] \ar[dd, shift left, "U_{\SigmaNetCat}"] &
  \CMC \ar[dd,shift left, "U_{\Petri}"] \ar[l,shift left, "U_{\SSMC}"] \\
  \\
  \PreNet 
  \ar[uu,shift left, "F_{\PreNet}"] \ar[r,shift left=2,"F_{\pre}"] \ar[r,shift right=2,swap,"H_{\pre}"] 
  &
  \SigmaNetCat \ar[uu,shift left, "F_{\SigmaNetCat}"] \ar[l] \ar[r,shift left, "F_{\pet}"]
  &
  \Petri \ar[uu,shift left, "F_{\Petri}"] \ar[l,shift left, "G_{\pet}"]
\end{tikzcd}
\end{equation}
The adjunctions in the top row can be constructed using standard tools, such as the adjoint functor theorem or the adjoint lifting theorem.
\begin{prop}\label{prop: free symmetries}
  There is an adjunction
  \[
    \begin{tikzcd}[row sep=large,column sep=large]
      \StrMC \ar[r, bend left, "F_{\StrMC}"] \ar[r, phantom, "\bot"]  &
      \SSMC. \ar[l, bend left, "U_{\StrMC}"]
    \end{tikzcd}  
  \]
\end{prop}
\noindent Here, $U_{\StrMC}$ freely adds symmetries to a strict monoidal category, while $U_{\StrMC}$ sends any symmetric strict monoidal category to its underlying strict monoidal category.
\begin{prop}
  There is an adjunction
  \[
    \begin{tikzcd}[row sep=large,column sep=large]
      \SSMC \ar[r,bend left, "F_{\SSMC}"] \ar[r, phantom, "\bot"]  &
      \CMC. \ar[l,bend left, "U_{\SSMC}"]
    \end{tikzcd}  
  \]
\end{prop}
\noindent $F_{\SSMC}$ takes a symmetric strict monoidal category and imposes a law saying that all symmetries are identity morphism, while $U_{\SSMC}$ sends any commutative monoidal category to its underlying symmetric strict monoidal category.

The adjunction between $\Petri$ and $\CMC$ was recalled in \cref{prop: adjunction Petri - CMC}, while that between $\PreNet$ and $\StrMC$ was recalled in \cref{prop: adjunction PreNet - StrMC}. We now cover the middle column and bottom row of the diagram, which are new.
\begin{thm}
\label{thm:pre-net - sigma-net adjunction}
  There is a triple of adjoint functors
  \[
    \begin{tikzcd}[row sep=large,column sep=1.5cm]
      \PreNet \ar[r, bend left, shift left=3, "F_{\pre}"]
      \ar[r, shift left=5, phantom, "\bot"] 
      \ar[r, shift right=5, phantom, "\bot"] 
      \ar[r, bend right, shift right=3, "H_{\pre}"'] &
      \SigmaNetCat. \ar[l, "G_{\pre}" description]
    \end{tikzcd}  
  \]
\end{thm}
\begin{proof}
For this proof it is most convenient to work with the presheaf perspective.
In the proof of \cref{prop:topos of pre-nets} we described a category $\C$ such that 
$\PreNet \cong [\C,\Set]$
and in the proof of \cref{thm: sigma nets is a topos} we described a category $\D$ such that 
$\SigmaNetCat \cong [\D,\Set]$.  Recall that $\D$ is built by starting with the objects and morphisms of $\C$ and adding new morphisms and equations. The inclusion gives a functor $i \maps \C \to \D$ which induces a functor 
\[\SigmaNetCat \cong [\D,\Set] \xrightarrow{(-) \circ i} [\C,\Set] \cong \PreNet \]
given by precomposition with $i$. The composite functor above is the forgetful functor $G_{\pre}$. Therefore $G_{\pre}$ has a left adjoint $F_{\pre} \maps \PreNet \to \SigmaNetCat$ given by left Kan extension along $i$ and a right adjoint $H_{\pre} \maps \PreNet \to \SigmaNetCat$ given by right Kan extension along $i$. 
\end{proof}

Let us spell out what the functors $F_{\pre}$, $G_{\pre}$, $H_{\pre}$ do in detail, using our three-dimensional graphical representation.
\begin{description}
\item[$F_{\pre}$\iflics)\else~constructs tanks\fi]
  For this functor we work in the groupoid representation of \SigmaNet{}s.
  A pre-net $T \xrightarrow{(s,t)} S^{*}\times S^{*}$ is sent to the \SigmaNet{} $T \xrightarrow{F_{\pre}(s,t)} PS \times PS\op$, where $T$ denotes the \emph{discrete} groupoid having $T$ as underlying set of objects. Since $T$ is discrete, the functor $F_{\pre}(s,t)$ only needs to be defined on objects, which we do by taking the composite
  \[ T \xrightarrow{(s,t)} S^{*}\times S^{*} \to PS \times PS\op\]
  using the fact that $S^{*}$ is the set of objects of $PS$.
  A morphism of pre-nets $(f,g) \maps (s_1,t_1) \to (s_2,t_2)$ induces a morphism of \SigmaNet{}s: $g \maps T_1 \to T_2$ lifts to a morphism between discrete groupoids and $f \maps S_1 \to S_2$ lifts to a functor $PS_1 \times PS_1\op\to PS_2 \times PS_2\op$. The relevant square as in \cref{thm: fibrational representation} commutes strictly.

  $F_{\pre}$ takes a pre-net and builds from it a \SigmaNet{} with trivial isotropy groups. Graphically, this amounts to enclosing every transition of the given pre-net in a tank:
  \[
    \scalebox{\iflics 0.55\else 0.8\fi}
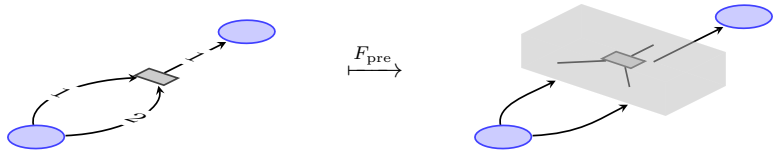
{  
\begin{tikzpicture}[baseline=0cm]
  \begin{scope}[
    yshift=0,every node/.append style={
    yslant=\ycoord,xslant=\xcoord},yslant=\ycoord,xslant=\xcoord
    ]    \begin{pgfonlayer}{nodelayer}
        \node [place,tokens=0] (1a) at (-2,0) {};
        \node [place,tokens=0] (3a) at (1.5,0) {};
        \node[transition] (2a) at (0,0) {};
    \end{pgfonlayer}
    \begin{pgfonlayer}{edgelayer}
        \draw[style=inarrow, out=45, in=135, thick] (1a) to node[rotate=-90, fill=white, midway] {$1$} (2a);
        \draw[style=inarrow, out=-45, in=225, thick] (1a) to node[rotate=-90, fill=white,midway] {$2$} (2a);
        \draw[style=inarrow, thick] (2a) to node[rotate=-90, fill=white, midway] {$1$} (3a);
    \end{pgfonlayer}
  \end{scope}
\end{tikzpicture}
} \iflics\quad\else\qquad\fi \xmapsto{F_\pre} \iflics\quad\else\qquad\fi \scalebox{\iflics 0.55\else 0.8\fi}{
\begin{tikzpicture}[baseline=0cm]
  \begin{scope}[
      yshift=0,every node/.append style={
      yslant=\ycoord,xslant=\xcoord},yslant=\ycoord,xslant=\xcoord
      ]
    \begin{pgfonlayer}{nodelayer}
      \node [place,tokens=0] (1a) at (-2,0) {};
      \node [place,tokens=0] (3a) at (2,0) {};
    \end{pgfonlayer}
    \begin{pgfonlayer}{edgelayer}
      \draw[style=inarrow, thick, out=45, in=180] (1a) to (-0.5,0.5);
      \draw[style=inarrow, thick, out=-45, in=180] (1a) to  (-0.5,-0.5);
      \draw[style=inarrow, thick] (0.5,0) to (3a);
    \end{pgfonlayer}
    
      \fill[\transitionbgrd, opacity=\transitionop]  (-0.5,-1) rectangle (0.5,1);
      \node (11) at (-0.5,1) {};
      \node (21) at (-0.5,-1) {};
      \node (31) at (0.5,-1) {};
      \node (41) at (0.5,1) {};
  \end{scope}

  \begin{scope}[
      yshift=8,every node/.append style={
      yslant=\ycoord,xslant=\xcoord},yslant=\ycoord,xslant=\xcoord
      ]
    \begin{pgfonlayer}{nodelayer}
      \node [transition] (2b) []  at (0,0) {};  
    \end{pgfonlayer}
    \begin{pgfonlayer}{edgelayer}
      \draw[thick, in=135, out=-45] (-0.5,0.5) to  (2b);
      \draw[thick,  in=225, out=45] (-0.5,-0.5) to  (2b);
      \draw[thick] (2b) to  (0.5,0);
    \end{pgfonlayer}
  \end{scope}


  \begin{scope}[
      yshift=16,every node/.append style={
      yslant=\ycoord,xslant=\xcoord},yslant=\ycoord,xslant=\xcoord
      ]

      \fill[\transitionbgrd, opacity=\transitionop]  (-0.5,-1) rectangle (0.5,1);
      \node (14) at (-0.5,1) {};
      \node (24) at (-0.5,-1) {};
      \node (34) at (0.5,-1) {};
      \node (44) at (0.5,1) {};
  \end{scope}

    \fill[\transitionbgrd, opacity=\transitionop]   (11.center) -- (14.center) -- (24.center) -- (21.center) -- (11.center);
    \fill[\transitionbgrd, opacity=\transitionop]   (21.center) -- (24.center) -- (34.center) -- (31.center) -- (21.center);
    \fill[\transitionbgrd, opacity=\transitionop]   (31.center) -- (34.center) -- (44.center) -- (41.center) -- (31.center);
    \fill[\transitionbgrd, opacity=\transitionop]   (41.center) -- (44.center) -- (14.center) -- (11.center) -- (41.center);

\end{tikzpicture}
}
  \]
  In particular, the \emph{transition classes} of $F_\pre(N)$ are the transitions of $N$. 

\item[$G_{\pre}$\iflics)\else~explodes tanks.\fi]
  For this functor we work in the profunctor representation.
  A \SigmaNet{} $PS \times PS\op \xrightarrow{N} \Set$ is sent to the pre-net having $S$ as its set of places and the disjoint union of all sets $N(a,b)$, for any $a,b$ objects of $PS$, as its set of transitions. For each transition, input  and output places are defined using the inverse image of $N$.
  That is, the transitions of $G_{\pre}N$ are the \emph{transitions} of $N$, with their grouping into classes and their isotropy groups forgotten.

  We can give a different interpretation of this using the groupoid perspective.  Suppose $T \xrightarrow{N} PS \times PS\op$ is a \SigmaNet{}.  Then for each object $t$ of $T$ such that $N(t)$ is a pair of strings of length $m$ and $n$ there will be $S_m \times S_n/\hom_{T}(t,t)$ transitions in $G_{\pre}N$, where $S_n$ denotes the group of permutations over a string of $n$ elements.  Graphically, this is represented by ``exploding'' a tank with $m$ inputs and $n$ outputs and introducing $m!n!/k$ pre-net transitions, where $k$ is the number of elements in the tank.
  \[
    \scalebox{\iflics 0.55\else 0.8\fi}
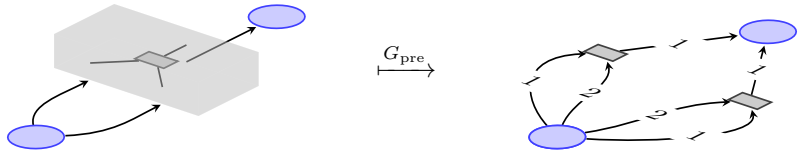
{
\begin{tikzpicture}[baseline=0cm]
  \begin{scope}[
      yshift=0,every node/.append style={
      yslant=\ycoord,xslant=\xcoord},yslant=\ycoord,xslant=\xcoord
      ]
    \begin{pgfonlayer}{nodelayer}
      \node [place,tokens=0] (1a) at (-2,0) {};
      \node [place,tokens=0] (3a) at (2,0) {};
    \end{pgfonlayer}
    \begin{pgfonlayer}{edgelayer}
      \draw[style=inarrow, thick, out=45, in=180] (1a) to (-0.5,0.5);
      \draw[style=inarrow, thick, out=-45, in=180] (1a) to  (-0.5,-0.5);
      \draw[style=inarrow, thick] (0.5,0) to (3a);
    \end{pgfonlayer}
    
      \fill[\transitionbgrd, opacity=\transitionop]  (-0.5,-1) rectangle (0.5,1);
      \node (11) at (-0.5,1) {};
      \node (21) at (-0.5,-1) {};
      \node (31) at (0.5,-1) {};
      \node (41) at (0.5,1) {};
  \end{scope}

  \begin{scope}[
      yshift=8,every node/.append style={
      yslant=\ycoord,xslant=\xcoord},yslant=\ycoord,xslant=\xcoord
      ]
    \begin{pgfonlayer}{nodelayer}
      \node [transition] (2b) []  at (0,0) {};  
    \end{pgfonlayer}
    \begin{pgfonlayer}{edgelayer}
      \draw[thick, in=135, out=-45] (-0.5,0.5) to  (2b);
      \draw[thick,  in=225, out=45] (-0.5,-0.5) to  (2b);
      \draw[thick] (2b) to  (0.5,0);
    \end{pgfonlayer}
  \end{scope}


  \begin{scope}[
      yshift=16,every node/.append style={
      yslant=\ycoord,xslant=\xcoord},yslant=\ycoord,xslant=\xcoord
      ]

      \fill[\transitionbgrd, opacity=\transitionop]  (-0.5,-1) rectangle (0.5,1);
      \node (14) at (-0.5,1) {};
      \node (24) at (-0.5,-1) {};
      \node (34) at (0.5,-1) {};
      \node (44) at (0.5,1) {};
  \end{scope}

    \fill[\transitionbgrd, opacity=\transitionop]   (11.center) -- (14.center) -- (24.center) -- (21.center) -- (11.center);
    \fill[\transitionbgrd, opacity=\transitionop]   (21.center) -- (24.center) -- (34.center) -- (31.center) -- (21.center);
    \fill[\transitionbgrd, opacity=\transitionop]   (31.center) -- (34.center) -- (44.center) -- (41.center) -- (31.center);
    \fill[\transitionbgrd, opacity=\transitionop]   (41.center) -- (44.center) -- (14.center) -- (11.center) -- (41.center);

\end{tikzpicture}
} \iflics\quad\else\qquad\fi \xmapsto{G_{\pre}} \iflics\quad\else\qquad\fi \scalebox{\iflics 0.55\else 0.8\fi}{  
\begin{tikzpicture}[baseline=0cm]
  \begin{scope}[
    yshift=0,every node/.append style={
    yslant=\ycoord,xslant=\xcoord},yslant=\ycoord,xslant=\xcoord
    ]
    \begin{pgfonlayer}{nodelayer}
        \node [place,tokens=0] (1a) at (-2,0) {};
        \node [place,tokens=0] (1b) at (-2,-0) {};
        \node [place,tokens=0] (3a) at (1.5,0) {};
        \node[transition] (2a) at (0,1) {};
        \node[transition] (2b) at (0,-1) {};

    \end{pgfonlayer}
    \begin{pgfonlayer}{edgelayer}
        \draw[style=inarrow, out=60, in=135, thick] (1a) to node[rotate=-90, fill=white, midway] {$1$} (2a);
        \draw[style=inarrow, out=30, in=225, thick] (1b) to node[rotate=-90, fill=white,midway] {$2$} (2a);
        \draw[style=inarrow, out=-30, in=135, thick] (1a) to node[rotate=-90, fill=white, midway] {$2$} (2b);
        \draw[style=inarrow, out=-60, in=225, thick] (1b) to node[rotate=-90, fill=white,midway] {$1$} (2b);
        \draw[style=inarrow, thick] (2a) to node[rotate=-90, fill=white, midway] {$1$} (3a);
        \draw[style=inarrow, thick] (2b) to node[rotate=-90, fill=white, midway] {$1$} (3a);
    \end{pgfonlayer}
  \end{scope}
\end{tikzpicture}
}
  \]
  In the image above, we see the behavior of $G_{\pre}$ on a $\Sigma$-net having a transition with trivial isotropy group, while in the image below $G_{\pre}$ is used on a \SigmaNet{} having a transition with 2-element isotropy group.
  \[
    \scalebox{\iflics 0.55\else 0.8\fi}
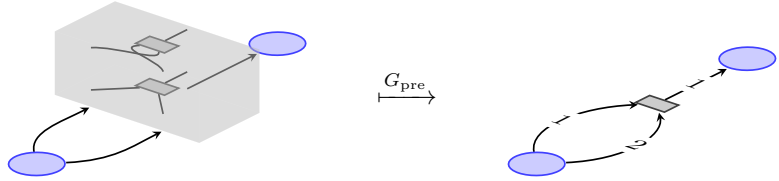
{
\begin{tikzpicture}[baseline=0cm]
  \begin{scope}[
      yshift=0,every node/.append style={
      yslant=\ycoord,xslant=\xcoord},yslant=\ycoord,xslant=\xcoord
      ]
    \begin{pgfonlayer}{nodelayer}
      \node [place,tokens=0] (1a) at (-2,0) {};
      \node [place,tokens=0] (3a) at (2,0) {};
    \end{pgfonlayer}
    \begin{pgfonlayer}{edgelayer}
      \draw[style=inarrow, thick, out=45, in=180] (1a) to (-0.5,0.5);
      \draw[style=inarrow, thick, out=-45, in=180] (1a) to  (-0.5,-0.5);
      \draw[style=inarrow, thick] (0.5,0) to (3a);
    \end{pgfonlayer}
    
      \fill[\transitionbgrd, opacity=\transitionop]  (-0.5,-1) rectangle (0.5,1);
      \node (11) at (-0.5,1) {};
      \node (21) at (-0.5,-1) {};
      \node (31) at (0.5,-1) {};
      \node (41) at (0.5,1) {};
  \end{scope}

  \begin{scope}[
      yshift=8,every node/.append style={
      yslant=\ycoord,xslant=\xcoord},yslant=\ycoord,xslant=\xcoord
      ]
    \begin{pgfonlayer}{nodelayer}
      \node [transition] (2b) []  at (0,0) {};  
    \end{pgfonlayer}
    \begin{pgfonlayer}{edgelayer}
      \draw[thick, in=135, out=-45] (-0.5,0.5) to  (2b);
      \draw[thick,  in=225, out=45] (-0.5,-0.5) to  (2b);
      \draw[thick] (2b) to  (0.5,0);
    \end{pgfonlayer}
  \end{scope}

  \begin{scope}[
      yshift=28,every node/.append style={
      yslant=\ycoord,xslant=\xcoord},yslant=\ycoord,xslant=\xcoord
      ]
    \begin{pgfonlayer}{nodelayer}
      \node [transition] (2b) []  at (0,0) {};  
    \end{pgfonlayer} 
    \begin{pgfonlayer}{edgelayer}
      \draw[thick, in=225, out=-45] (-0.5,0.5) to  (2b);
      \draw[thick,  in=135, out=45] (-0.5,-0.5) to  (2b);
      \draw[thick] (2b) to  (0.5,0);
    \end{pgfonlayer}
  \end{scope}

  \begin{scope}[
      yshift=30,every node/.append style={
      yslant=\ycoord,xslant=\xcoord},yslant=\ycoord,xslant=\xcoord
      ]

      \fill[\transitionbgrd, opacity=\transitionop]  (-0.5,-1) rectangle (0.5,1);
      \node (14) at (-0.5,1) {};
      \node (24) at (-0.5,-1) {};
      \node (34) at (0.5,-1) {};
      \node (44) at (0.5,1) {};
  \end{scope}

    \fill[\transitionbgrd, opacity=\transitionop]   (11.center) -- (14.center) -- (24.center) -- (21.center) -- (11.center);
    \fill[\transitionbgrd, opacity=\transitionop]   (21.center) -- (24.center) -- (34.center) -- (31.center) -- (21.center);
    \fill[\transitionbgrd, opacity=\transitionop]   (31.center) -- (34.center) -- (44.center) -- (41.center) -- (31.center);
    \fill[\transitionbgrd, opacity=\transitionop]   (41.center) -- (44.center) -- (14.center) -- (11.center) -- (41.center);

\end{tikzpicture}
} \iflics\quad\else\qquad\fi \xmapsto{G_{\pre}} \iflics\quad\else\qquad\fi \scalebox{\iflics 0.55\else 0.8\fi}{  
\begin{tikzpicture}[baseline=0cm]
  \begin{scope}[
    yshift=0,every node/.append style={
    yslant=\ycoord,xslant=\xcoord},yslant=\ycoord,xslant=\xcoord
    ]    \begin{pgfonlayer}{nodelayer}
        \node [place,tokens=0] (1a) at (-2,0) {};
        \node [place,tokens=0] (3a) at (1.5,0) {};
        \node[transition] (2a) at (0,0) {};
    \end{pgfonlayer}
    \begin{pgfonlayer}{edgelayer}
        \draw[style=inarrow, out=45, in=135, thick] (1a) to node[rotate=-90, fill=white, midway] {$1$} (2a);
        \draw[style=inarrow, out=-45, in=225, thick] (1a) to node[rotate=-90, fill=white,midway] {$2$} (2a);
        \draw[style=inarrow, thick] (2a) to node[rotate=-90, fill=white, midway] {$1$} (3a);
    \end{pgfonlayer}
  \end{scope}
\end{tikzpicture}
}
  \]
  \item[$H_{\pre}$\iflics)\else~matches transitions.\fi] While $F_{\pre}$ builds as many tanks as we can get from a pre-net's transitions, $H_{\pre}$ bundles pre-net transitions sharing the same inputs/outputs modulo permutations, whenever they complete their corresponding symmetry groups. For instance, in the figure below transitions $x$ and $y$ complete the permutation group $S_2 \times S_1$, and hence they give rise to the tank denoted with $\langle x, y \rangle$. The same
  happens for transitions $x$ and $z$, giving rise to tank $\langle x , z \rangle$.
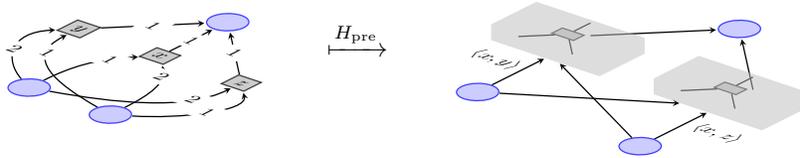
  \[
    \scalebox{\iflics 0.43\else 0.6\fi}{  
\begin{tikzpicture}[baseline=0cm]
  \begin{scope}[
    yshift=0,every node/.append style={
    yslant=\ycoord,xslant=\xcoord},yslant=\ycoord,xslant=\xcoord
    ]
    \begin{pgfonlayer}{nodelayer}
        \node [place,tokens=0] (1a) at (-2,0.75) {};
        \node [place,tokens=0] (1b) at (-2,-0.75) {};
        \node [place,tokens=0] (3a) at (1.5,0) {};
        \node[transition] (2a) at (0,0) {\rotatebox{-90}{$x$}};
        \node[transition] (2b) at (0,1.5) {\rotatebox{-90}{$y$}};
        \node[transition] (2c) at (0,-1.5) {\rotatebox{-90}{$z$}};

    \end{pgfonlayer}
    \begin{pgfonlayer}{edgelayer}
        \draw[style=inarrow, out=0, in=135, thick] (1a) to node[rotate=-90, fill=white, near end] {$1$} (2a);
        \draw[style=inarrow, out=0, in=-135, thick] (1b) to node[rotate=-90, fill=white,near end] {$2$} (2a);
        \draw[style=inarrow, out=60, in=135, thick] (1a) to node[rotate=-90, fill=white, midway] {$2$} (2b);
        \draw[style=inarrow, out=75, in=180, thick] (1b) to node[rotate=-90, fill=white,near end] {$1$} (2b);
        \draw[style=inarrow, out=-75, in=180, thick] (1a) to node[rotate=-90, fill=white, near end] {$2$} (2c);
        \draw[style=inarrow, out=-60, in=-135, thick] (1b) to node[rotate=-90, fill=white, midway] {$1$} (2c);
        \draw[style=inarrow, thick] (2a) to node[rotate=-90, fill=white, midway] {$1$} (3a);
        \draw[style=inarrow, thick] (2b) to node[rotate=-90, fill=white, midway] {$1$} (3a);
        \draw[style=inarrow, thick] (2c) to node[rotate=-90, fill=white, midway] {$1$} (3a);
    \end{pgfonlayer}
  \end{scope}
\end{tikzpicture}
} \iflics\,\else\qquad\fi \xmapsto{H_{\pre}} \iflics\,\else\qquad\fi \scalebox{\iflics 0.43\else 0.6\fi}{
\begin{tikzpicture}[baseline=1cm]
  \begin{scope}[
      yshift=0,every node/.append style={
      yslant=\ycoord,xslant=\xcoord},yslant=\ycoord,xslant=\xcoord
      ]
    \begin{pgfonlayer}{nodelayer}
      \node [place,tokens=0] (1a) at (-2,3) {};
      \node [place,tokens=0] (1b) at (-2,0) {};      			
      \node [place,tokens=0] (3a) at (2,1.5) {};
    \end{pgfonlayer}
    \begin{pgfonlayer}{edgelayer}
      \draw[style=inarrow, thick] (1a) to node[right, rotate=-90]{} (-0.5,0.5);
      \draw[style=inarrow, thick] (1b) to  (-0.5,0);
      \draw[style=inarrow, thick] (0.5,0) to (3a);
      \draw[style=inarrow, thick] (1a) to node[left, rotate=-90]{} (-0.5,3);
      \draw[style=inarrow, thick] (1b) to  (-0.5,2.75);
      \draw[style=inarrow, thick] (0.5,3) to (3a);
    \end{pgfonlayer}

      \fill[\transitionbgrd, opacity=\transitionop]  (-0.5,2) rectangle (0.5,4);
      \node (13) at (-0.5,4) {};
      \node (23) at (-0.5,2) {};
      \node (33) at (0.5,2) {};
      \node (43) at (0.5,4) {};
      
      \fill[\transitionbgrd, opacity=\transitionop]  (-0.5,-1) rectangle (0.5,1);
      \node (11) at (-0.5,1) {};
      \node (21) at (-0.5,-1) {};
      \node (31) at (0.5,-1) {};
      \node (41) at (0.5,1) {};
  \end{scope}

  \begin{scope}[
      yshift=8,every node/.append style={
      yslant=\ycoord,xslant=\xcoord},yslant=\ycoord,xslant=\xcoord
      ]
      
    \begin{pgfonlayer}{nodelayer}
      \node[transition, label={[xshift=-1.0cm, yshift=-0.5cm]\rotatebox{-90}{$\langle x,y \rangle$}}] (2a) []  at (0,3) {};  
      \node[transition, label={[xshift=-1.01cm, yshift=-1.6cm]\rotatebox{-90}{$\langle x,z \rangle$}}] (2b) []  at (0,0) {};  
    \end{pgfonlayer}
    \begin{pgfonlayer}{edgelayer}
      \draw[thick, in=135, out=-45] (-0.5,0.5) to  (2b);
      \draw[thick,  in=225, out=45] (-0.5,-0.5) to  (2b);
      \draw[thick] (2b) to  (0.5,0);

      \draw[thick, in=135, out=-45] (-0.5,3.5) to  (2a);
      \draw[thick,  in=225, out=45] (-0.5,2.5) to  (2a);
      \draw[thick] (2a) to  (0.5,3);
    \end{pgfonlayer}
  \end{scope}
  
  \begin{scope}[
      yshift=10,every node/.append style={
      yslant=\ycoord,xslant=\xcoord},yslant=\ycoord,xslant=\xcoord
      ]
      \fill[\transitionbgrd, opacity=\transitionop]  (-0.5,2) rectangle (0.5,4);

      \node (14) at (-0.5,4) {};
      \node (24) at (-0.5,2) {};
      \node (34) at (0.5,2) {};
      \node (44) at (0.5,4) {};

      \fill[\transitionbgrd, opacity=\transitionop]  (-0.5,-1) rectangle (0.5,1);
      \node (12) at (-0.5,1) {};
      \node (22) at (-0.5,-1) {};
      \node (32) at (0.5,-1) {};
      \node (42) at (0.5,1) {};
  \end{scope}

    \fill[\transitionbgrd, opacity=\transitionop]   (11.center) -- (12.center) -- (22.center) -- (21.center) -- (11.center);
    \fill[\transitionbgrd, opacity=\transitionop]   (21.center) -- (22.center) -- (32.center) -- (31.center) -- (21.center);
    \fill[\transitionbgrd, opacity=\transitionop]   (31.center) -- (32.center) -- (42.center) -- (41.center) -- (31.center);
    \fill[\transitionbgrd, opacity=\transitionop]   (41.center) -- (42.center) -- (12.center) -- (11.center) -- (41.center);

    \fill[\transitionbgrd, opacity=\transitionop]   (13.center) -- (14.center) -- (24.center) -- (23.center) -- (13.center);
    \fill[\transitionbgrd, opacity=\transitionop]   (23.center) -- (24.center) -- (34.center) -- (33.center) -- (23.center);
    \fill[\transitionbgrd, opacity=\transitionop]   (33.center) -- (34.center) -- (44.center) -- (43.center) -- (33.center);
    \fill[\transitionbgrd, opacity=\transitionop]   (43.center) -- (44.center) -- (14.center) -- (13.center) -- (43.center);

\end{tikzpicture}
}
  \]
  The following pre-net does not have enough transitions to complete the symmetry group of its inputs/outputs. As such, $H_{\pre}$ cannot match this transition with anything, and does not produce any tank.
  \[
    \scalebox{\iflics 0.45\else 0.7\fi}
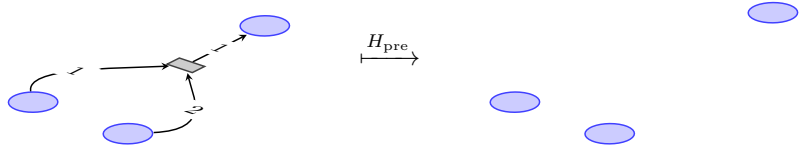
{  
\begin{tikzpicture}[baseline=0cm]
  \begin{scope}[
    yshift=0,every node/.append style={
    yslant=\ycoord,xslant=\xcoord},yslant=\ycoord,xslant=\xcoord
    ]    \begin{pgfonlayer}{nodelayer}
        \node [place,tokens=0] (1a) at (-2,0.75) {};
        \node [place,tokens=0] (1b) at (-2,-0.75) {}; 
        \node [place,tokens=0] (3a) at (1.5,0) {};
        \node[transition] (2a) at (0,0) {};
    \end{pgfonlayer}
    \begin{pgfonlayer}{edgelayer}
      \draw[style=inarrow, out=45, in=135, thick] (1a) to node[rotate=-90, fill=white, midway] {$1$} (2a);
      \draw[style=inarrow, out=-45, in=225, thick] (1b) to node[rotate=-90, fill=white,midway] {$2$} (2a);
      \draw[style=inarrow, thick] (2a) to node[rotate=-90, fill=white, midway] {$1$} (3a);
    \end{pgfonlayer}
  \end{scope}
\end{tikzpicture}
} \iflics\quad\else\qquad\fi \xmapsto{H_{\pre}} \iflics\quad\else\qquad\fi \scalebox{\iflics 0.45\else 0.7\fi}{
\begin{tikzpicture}[baseline=0cm]
  \begin{scope}[
      yshift=0,every node/.append style={
      yslant=\ycoord,xslant=\xcoord},yslant=\ycoord,xslant=\xcoord
      ]
    \begin{pgfonlayer}{nodelayer}
      \node [place,tokens=0] (1a) at (-2,-0.75) {};
      \node [place,tokens=0] (1b) at (-2,0.75) {};
      \node [place,tokens=0] (3a) at (2,0) {};
    \end{pgfonlayer}
    
  \end{scope}

\end{tikzpicture}
}
  \]
  In the following case, the pre-net has a repeated input. $H_{\pre}$ is then able to match the transtion with itself, producing a maximally commutative tank.
  \[
    \scalebox{\iflics 0.57\else 0.8\fi}
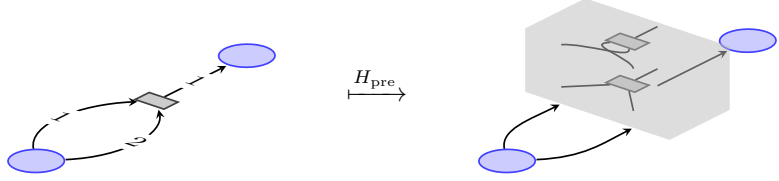
{  
\begin{tikzpicture}[baseline=0cm]
  \begin{scope}[
    yshift=0,every node/.append style={
    yslant=\ycoord,xslant=\xcoord},yslant=\ycoord,xslant=\xcoord
    ]    \begin{pgfonlayer}{nodelayer}
        \node [place,tokens=0] (1a) at (-2,0) {};
        \node [place,tokens=0] (3a) at (1.5,0) {};
        \node[transition] (2a) at (0,0) {};
    \end{pgfonlayer}
    \begin{pgfonlayer}{edgelayer}
        \draw[style=inarrow, out=45, in=135, thick] (1a) to node[rotate=-90, fill=white, midway] {$1$} (2a);
        \draw[style=inarrow, out=-45, in=225, thick] (1a) to node[rotate=-90, fill=white,midway] {$2$} (2a);
        \draw[style=inarrow, thick] (2a) to node[rotate=-90, fill=white, midway] {$1$} (3a);
    \end{pgfonlayer}
  \end{scope}
\end{tikzpicture}
} \iflics\quad\else\qquad\fi \xmapsto{H_{\pre}} \iflics\quad\else\qquad\fi \scalebox{\iflics 0.57\else 0.8\fi}{
\begin{tikzpicture}[baseline=0cm]
  \begin{scope}[
      yshift=0,every node/.append style={
      yslant=\ycoord,xslant=\xcoord},yslant=\ycoord,xslant=\xcoord
      ]
    \begin{pgfonlayer}{nodelayer}
      \node [place,tokens=0] (1a) at (-2,0) {};
      \node [place,tokens=0] (3a) at (2,0) {};
    \end{pgfonlayer}
    \begin{pgfonlayer}{edgelayer}
      \draw[style=inarrow, thick, out=45, in=180] (1a) to (-0.5,0.5);
      \draw[style=inarrow, thick, out=-45, in=180] (1a) to  (-0.5,-0.5);
      \draw[style=inarrow, thick] (0.5,0) to (3a);
    \end{pgfonlayer}
    
      \fill[\transitionbgrd, opacity=\transitionop]  (-0.5,-1) rectangle (0.5,1);
      \node (11) at (-0.5,1) {};
      \node (21) at (-0.5,-1) {};
      \node (31) at (0.5,-1) {};
      \node (41) at (0.5,1) {};
  \end{scope}

  \begin{scope}[
      yshift=8,every node/.append style={
      yslant=\ycoord,xslant=\xcoord},yslant=\ycoord,xslant=\xcoord
      ]
    \begin{pgfonlayer}{nodelayer}
      \node [transition] (2b) []  at (0,0) {};  
    \end{pgfonlayer}
    \begin{pgfonlayer}{edgelayer}
      \draw[thick, in=135, out=-45] (-0.5,0.5) to  (2b);
      \draw[thick,  in=225, out=45] (-0.5,-0.5) to  (2b);
      \draw[thick] (2b) to  (0.5,0);
    \end{pgfonlayer}
  \end{scope}

  \begin{scope}[
      yshift=28,every node/.append style={
      yslant=\ycoord,xslant=\xcoord},yslant=\ycoord,xslant=\xcoord
      ]
    \begin{pgfonlayer}{nodelayer}
      \node [transition] (2b) []  at (0,0) {};  
    \end{pgfonlayer} 
    \begin{pgfonlayer}{edgelayer}
      \draw[thick, in=225, out=-45] (-0.5,0.5) to  (2b);
      \draw[thick,  in=135, out=45] (-0.5,-0.5) to  (2b);
      \draw[thick] (2b) to  (0.5,0);
    \end{pgfonlayer}
  \end{scope}

  \begin{scope}[
      yshift=30,every node/.append style={
      yslant=\ycoord,xslant=\xcoord},yslant=\ycoord,xslant=\xcoord
      ]

      \fill[\transitionbgrd, opacity=\transitionop]  (-0.5,-1) rectangle (0.5,1);
      \node (14) at (-0.5,1) {};
      \node (24) at (-0.5,-1) {};
      \node (34) at (0.5,-1) {};
      \node (44) at (0.5,1) {};
  \end{scope}

    \fill[\transitionbgrd, opacity=\transitionop]   (11.center) -- (14.center) -- (24.center) -- (21.center) -- (11.center);
    \fill[\transitionbgrd, opacity=\transitionop]   (21.center) -- (24.center) -- (34.center) -- (31.center) -- (21.center);
    \fill[\transitionbgrd, opacity=\transitionop]   (31.center) -- (34.center) -- (44.center) -- (41.center) -- (31.center);
    \fill[\transitionbgrd, opacity=\transitionop]   (41.center) -- (44.center) -- (14.center) -- (11.center) -- (41.center);

\end{tikzpicture}
}
  \]
  Looking at these examples, we see that in general the transitions of $N$ do not correspond directly to either the transitions of $H_\pre(N)$ or the transition classes of $H_\pre(N)$.
\end{description}
\begin{prop}\label{thm:petri-net-sigma-net-adjunction}
  There is an adjunction
  \[
    \begin{tikzcd}[row sep=large,column sep=1.5cm]
      \SigmaNetCat \ar[r,bend left, "F_{\pet}"] \ar[r, phantom, "\bot"] &
      \Petri. \ar[l, bend left, "G_{\pet}"]
    \end{tikzcd}  
  \]
\end{prop}

\begin{proof}
  Note first that \Petri{} is by definition precisely the comma category $\Comma{\Set}{(\N[-] \times \N[-])}$.
  Similarly, if we identify a \SigmaNet{} with a functor $N \maps PS\times PS\op \to \Set$ and thereby with a discrete opfibration $N \to PS\times PS\op$, then \SigmaNetCat{} becomes identified with the full subcategory of the comma category $\Comma{\Cat}{(P(-)\times P(-)\op)}$ consisting of the discrete opfibrations.

  Now note that $\Set$ is a reflective full subcategory of $\Cat$, with reflector $\pi_0$ that takes the set of connected components of a category.
  Moreover, we have $\pi_0(PS\times PS\op) \cong \N[S] \times \N[S]$.
  Thus, \cref{lem:comma-adj}, proven below (and applied with $\D = \Cat$, $\C = \E = \Set$, and $K = P(-) \times P(-)\op$), shows that we have an adjunction
  \[
    \adjustbox{scale=\iflics 0.89\else 1.0\fi}{
    \begin{tikzcd}
      \Comma{\Cat}{(P(-)\times P(-)\op)} \ar[r,shift left,"F"] &
      \Comma{\Set}{(\N[-] \times \N[-])} = \Petri \ar[l, shift left,"G"]
    \end{tikzcd}}
  \]
  in which the left adjoint $F$ applies $\pi_0$ to both domain and codomain, and the right adjoint $G$ pulls back along the unit $PS\times PS\op \to \N[S]\times \N[S]$.
  Therefore, it suffices to observe that this right adjoint takes values in discrete opfibrations, hence in \SigmaNetCat.
\end{proof}

\begin{lem}\label{lem:comma-adj}
  Let $\E$ be a reflective subcategory of $\D$, with reflector $\pi \maps \D \to \E$, and let $K \maps \C\to\D$ be a functor where $\D$ has pullbacks.  Then there is an adjunction
  \[
    \begin{tikzcd}
      \Comma{\D}{K} \ar[r,shift left,"F"] &
      \Comma{\E}{(\pi\circ K)}. \ar[l, shift left,"G"]
    \end{tikzcd}  
  \]
\end{lem}
\begin{proof}
  Let $\eta_X \maps X\to \pi X$ denote the unit of the reflection.  Then for any $f \maps S_1\to S_2$ in \C{}, we have $\eta_{K S_2}\circ K f = \pi K f \circ \eta_{K S_1}$ by naturality; we denote this common map by $\eta_f$.  Now there is a profunctor between $\Comma{\D}{K}$ and $\Comma{\E}{(\pi\circ K)}$ defined to take $T_1 \to K S_1$ and $T_2 \to \pi K S_2$ (where $T_1\in\D$ and $T_2\in\E$) to the set of pairs $(f,g)$ where $f:S_1\to S_2$ in \C{} and $g:T_1 \to T_2$ in \D{} make the following square commute:
  \[
    \begin{tikzcd}
      T_1 \ar[r,"g"] \ar[d] & T_2 \ar[d] \\
      K S_1 \ar[r,"\eta_f"'] & \pi K S_2.
    \end{tikzcd}
  \]
  This profunctor is representable on both sides, because any such square factors uniquely in both of the following ways:
  \[
    \adjustbox{scale=\iflics 0.82\else 1.0 \fi}{
    \begin{tikzcd}
      T_1 \ar[r,"\eta_{T_1}"] \ar[d] & \pi T_1 \ar[d] \ar[r,dashed] & T_2 \ar[d] \\
      K S_1 \ar[r,"\eta_{K S_1}"'] & \pi K S_1 \ar[r,"\pi K f"'] & \pi K S_2
    \end{tikzcd}
    \iflics\quad\else\qquad\fi
    \begin{tikzcd}
      T_1 \ar[r,dashed] \ar[d] & \bullet \ar[d] \ar[r] \ar[dr,phantom,near start,"\lrcorner"] & T_2 \ar[d] \\
      K S_1 \ar[r,"K f"'] & K S_2 \ar[r,"\eta_{K S_2}"'] & \pi K S_2.
    \end{tikzcd}}
  \]
  On the left, the factorization is by the universal property of $\eta_{T_1}$, while on the right it is by the universal property of the pullback.
  Therefore, there is an adjunction $F\dashv G$ as desired, where $F$ takes $T_1 \to K S_1$ to $\pi T_1 \to \pi K S_1$, and $G$ takes $T_2 \to \pi K S_2$ to the pullback of $T_2$ to $K S_2$.
\end{proof}

Note that by construction, this adjunction is a reflection, i.e., the right adjoint $G_{\pet}$ is fully faithful.
We can illustrate the action of $F_{\pet}$ and $G_{\pet}$ with examples.
\begin{description}
  \item[$F_{\pet}$\iflics)\else~deflates tanks.\fi] In the groupoid perspective, this functor takes a \SigmaNet $T \xrightarrow{N} PS \times PS\op$ and maps it to the Petri net having the underlying set of objects of $T$ as transitions, $S$ as places, and input/output functions induced by the mapping on objects of $N$. The action on a morphism $(g,f)$ is obtained by restricting the functor $g$ to its mapping on objects. Graphically, $F_{\pet}$ deflates tanks, replacing each tank by a single transition:
  \begin{gather*}
    \scalebox{\iflics 0.55\else 0.8\fi}
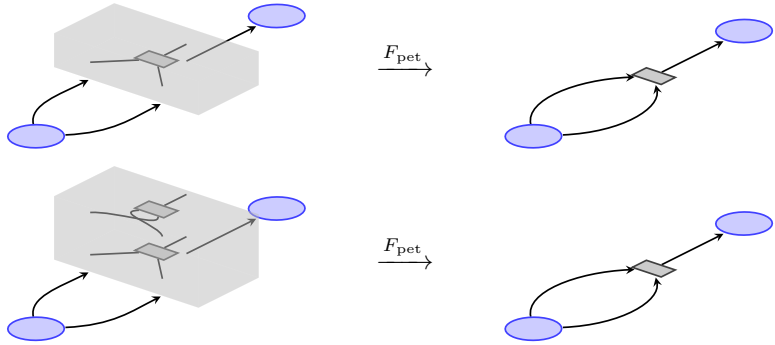
{
\begin{tikzpicture}[baseline=0cm]
  \begin{scope}[
      yshift=0,every node/.append style={
      yslant=\ycoord,xslant=\xcoord},yslant=\ycoord,xslant=\xcoord
      ]
    \begin{pgfonlayer}{nodelayer}
      \node [place,tokens=0] (1a) at (-2,0) {};
      \node [place,tokens=0] (3a) at (2,0) {};
    \end{pgfonlayer}
    \begin{pgfonlayer}{edgelayer}
      \draw[style=inarrow, thick, out=45, in=180] (1a) to (-0.5,0.5);
      \draw[style=inarrow, thick, out=-45, in=180] (1a) to  (-0.5,-0.5);
      \draw[style=inarrow, thick] (0.5,0) to (3a);
    \end{pgfonlayer}
    
      \fill[\transitionbgrd, opacity=\transitionop]  (-0.5,-1) rectangle (0.5,1);
      \node (11) at (-0.5,1) {};
      \node (21) at (-0.5,-1) {};
      \node (31) at (0.5,-1) {};
      \node (41) at (0.5,1) {};
  \end{scope}

  \begin{scope}[
      yshift=8,every node/.append style={
      yslant=\ycoord,xslant=\xcoord},yslant=\ycoord,xslant=\xcoord
      ]
    \begin{pgfonlayer}{nodelayer}
      \node [transition] (2b) []  at (0,0) {};  
    \end{pgfonlayer}
    \begin{pgfonlayer}{edgelayer}
      \draw[thick, in=135, out=-45] (-0.5,0.5) to  (2b);
      \draw[thick,  in=225, out=45] (-0.5,-0.5) to  (2b);
      \draw[thick] (2b) to  (0.5,0);
    \end{pgfonlayer}
  \end{scope}


  \begin{scope}[
      yshift=16,every node/.append style={
      yslant=\ycoord,xslant=\xcoord},yslant=\ycoord,xslant=\xcoord
      ]

      \fill[\transitionbgrd, opacity=\transitionop]  (-0.5,-1) rectangle (0.5,1);
      \node (14) at (-0.5,1) {};
      \node (24) at (-0.5,-1) {};
      \node (34) at (0.5,-1) {};
      \node (44) at (0.5,1) {};
  \end{scope}

    \fill[\transitionbgrd, opacity=\transitionop]   (11.center) -- (14.center) -- (24.center) -- (21.center) -- (11.center);
    \fill[\transitionbgrd, opacity=\transitionop]   (21.center) -- (24.center) -- (34.center) -- (31.center) -- (21.center);
    \fill[\transitionbgrd, opacity=\transitionop]   (31.center) -- (34.center) -- (44.center) -- (41.center) -- (31.center);
    \fill[\transitionbgrd, opacity=\transitionop]   (41.center) -- (44.center) -- (14.center) -- (11.center) -- (41.center);

\end{tikzpicture}
} \iflics\quad\else\qquad\fi \xrightarrow{F_{\pet}} \iflics\quad\else\qquad\fi \scalebox{\iflics 0.55\else 0.8\fi}{  
\begin{tikzpicture}[baseline=0cm]
    \begin{scope}[
      yshift=0,every node/.append style={
      yslant=\ycoord,xslant=\xcoord},yslant=\ycoord,xslant=\xcoord
      ]
    \begin{pgfonlayer}{nodelayer}
        \node [place,tokens=0] (1a) at (-2,0) {};
        \node [place,tokens=0] (3a) at (1.5,0) {};
        \node[transition] (2a) at (0,0) {};
    \end{pgfonlayer}
    \begin{pgfonlayer}{edgelayer}
      \draw[style=inarrow, out=45, in=135, thick] (1a) to (2a);
      \draw[style=inarrow, out=-45, in=225, thick] (1a) to (2a);
      \draw[style=inarrow, thick] (2a) to (3a);
    \end{pgfonlayer}
  \end{scope}
\end{tikzpicture}
}\\
    \scalebox{\iflics 0.55\else 0.8\fi}{
\begin{tikzpicture}[baseline=0cm]
  \begin{scope}[
      yshift=0,every node/.append style={
      yslant=\ycoord,xslant=\xcoord},yslant=\ycoord,xslant=\xcoord
      ]
    \begin{pgfonlayer}{nodelayer}
      \node [place,tokens=0] (1a) at (-2,0) {};
      \node [place,tokens=0] (3a) at (2,0) {};
    \end{pgfonlayer}
    \begin{pgfonlayer}{edgelayer}
      \draw[style=inarrow, thick, out=45, in=180] (1a) to (-0.5,0.5);
      \draw[style=inarrow, thick, out=-45, in=180] (1a) to  (-0.5,-0.5);
      \draw[style=inarrow, thick] (0.5,0) to (3a);
    \end{pgfonlayer}
    
      \fill[\transitionbgrd, opacity=\transitionop]  (-0.5,-1) rectangle (0.5,1);
      \node (11) at (-0.5,1) {};
      \node (21) at (-0.5,-1) {};
      \node (31) at (0.5,-1) {};
      \node (41) at (0.5,1) {};
  \end{scope}

  \begin{scope}[
      yshift=8,every node/.append style={
      yslant=\ycoord,xslant=\xcoord},yslant=\ycoord,xslant=\xcoord
      ]
    \begin{pgfonlayer}{nodelayer}
      \node [transition] (2b) []  at (0,0) {};  
    \end{pgfonlayer}
    \begin{pgfonlayer}{edgelayer}
      \draw[thick, in=135, out=-45] (-0.5,0.5) to  (2b);
      \draw[thick,  in=225, out=45] (-0.5,-0.5) to  (2b);
      \draw[thick] (2b) to  (0.5,0);
    \end{pgfonlayer}
  \end{scope}

  \begin{scope}[
      yshift=28,every node/.append style={
      yslant=\ycoord,xslant=\xcoord},yslant=\ycoord,xslant=\xcoord
      ]
    \begin{pgfonlayer}{nodelayer}
      \node [transition] (2b) []  at (0,0) {};  
    \end{pgfonlayer} 
    \begin{pgfonlayer}{edgelayer}
      \draw[thick, in=225, out=-45] (-0.5,0.5) to  (2b);
      \draw[thick,  in=135, out=45] (-0.5,-0.5) to  (2b);
      \draw[thick] (2b) to  (0.5,0);
    \end{pgfonlayer}
  \end{scope}

  \begin{scope}[
      yshift=30,every node/.append style={
      yslant=\ycoord,xslant=\xcoord},yslant=\ycoord,xslant=\xcoord
      ]

      \fill[\transitionbgrd, opacity=\transitionop]  (-0.5,-1) rectangle (0.5,1);
      \node (14) at (-0.5,1) {};
      \node (24) at (-0.5,-1) {};
      \node (34) at (0.5,-1) {};
      \node (44) at (0.5,1) {};
  \end{scope}

    \fill[\transitionbgrd, opacity=\transitionop]   (11.center) -- (14.center) -- (24.center) -- (21.center) -- (11.center);
    \fill[\transitionbgrd, opacity=\transitionop]   (21.center) -- (24.center) -- (34.center) -- (31.center) -- (21.center);
    \fill[\transitionbgrd, opacity=\transitionop]   (31.center) -- (34.center) -- (44.center) -- (41.center) -- (31.center);
    \fill[\transitionbgrd, opacity=\transitionop]   (41.center) -- (44.center) -- (14.center) -- (11.center) -- (41.center);

\end{tikzpicture}
} \iflics\quad\else\qquad\fi \xrightarrow{F_{\pet}} \iflics\quad\else\qquad\fi \scalebox{\iflics 0.55\else 0.8\fi}{  
\begin{tikzpicture}[baseline=0cm]
    \begin{scope}[
      yshift=0,every node/.append style={
      yslant=\ycoord,xslant=\xcoord},yslant=\ycoord,xslant=\xcoord
      ]
    \begin{pgfonlayer}{nodelayer}
        \node [place,tokens=0] (1a) at (-2,0) {};
        \node [place,tokens=0] (3a) at (1.5,0) {};
        \node[transition] (2a) at (0,0) {};
    \end{pgfonlayer}
    \begin{pgfonlayer}{edgelayer}
      \draw[style=inarrow, out=45, in=135, thick] (1a) to (2a);
      \draw[style=inarrow, out=-45, in=225, thick] (1a) to (2a);
      \draw[style=inarrow, thick] (2a) to (3a);
    \end{pgfonlayer}
  \end{scope}
\end{tikzpicture}
}
  \end{gather*}
  In particular, the transitions of $F_\pet(N)$ are the \emph{transition classes} of $N$.
  \item[$G_{\pet}$\iflics)\else~builds tanks as big as possible.\fi] Petri nets are mapped under $G_\pet$ to corresponding \SigmaNet{}s that have the largest isotropy groups possible. Consider a transition $t$ in a Petri net $N$. Its inputs and outputs will be a couple of unordered strings of length $n,m$, respectively. Pick any ordering for these strings, and call them $a,b$, respectively. Finally, let $G_t$ be the subgroup of $S_n \times S_m$ that fixes the pair of strings $(a,b)$.
  
  $N$ is mapped to a \SigmaNet $T \xrightarrow{G_{\pet} N} PS \times PS\op$ whose groupoid $T$ has transitions of $N$ as objects and, for each $t$ in $T$, $G_t$ as its group of automorphisms. $N$ maps $t$ to the ordering $(a,b)$ we have chosen before. It can be seen that picking different orderings of the input/output of each transition gives isomorphic results.

  Graphically, out of each Petri net we build a corresponding \SigmaNet{} that has its tanks as full as possible:
  \[
    \scalebox{\iflics 0.55\else 0.8\fi}
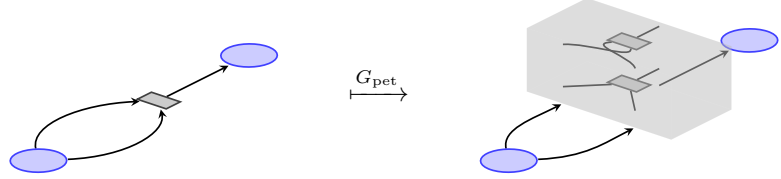
{  
\begin{tikzpicture}[baseline=0cm]
    \begin{scope}[
      yshift=0,every node/.append style={
      yslant=\ycoord,xslant=\xcoord},yslant=\ycoord,xslant=\xcoord
      ]
    \begin{pgfonlayer}{nodelayer}
        \node [place,tokens=0] (1a) at (-2,0) {};
        \node [place,tokens=0] (3a) at (1.5,0) {};
        \node[transition] (2a) at (0,0) {};
    \end{pgfonlayer}
    \begin{pgfonlayer}{edgelayer}
      \draw[style=inarrow, out=45, in=135, thick] (1a) to (2a);
      \draw[style=inarrow, out=-45, in=225, thick] (1a) to (2a);
      \draw[style=inarrow, thick] (2a) to (3a);
    \end{pgfonlayer}
  \end{scope}
\end{tikzpicture}
} \iflics\quad\else\qquad\fi \xmapsto{G_{\pet}} \iflics\quad\else\qquad\fi \scalebox{\iflics 0.55\else 0.8\fi}{
\begin{tikzpicture}[baseline=0cm]
  \begin{scope}[
      yshift=0,every node/.append style={
      yslant=\ycoord,xslant=\xcoord},yslant=\ycoord,xslant=\xcoord
      ]
    \begin{pgfonlayer}{nodelayer}
      \node [place,tokens=0] (1a) at (-2,0) {};
      \node [place,tokens=0] (3a) at (2,0) {};
    \end{pgfonlayer}
    \begin{pgfonlayer}{edgelayer}
      \draw[style=inarrow, thick, out=45, in=180] (1a) to (-0.5,0.5);
      \draw[style=inarrow, thick, out=-45, in=180] (1a) to  (-0.5,-0.5);
      \draw[style=inarrow, thick] (0.5,0) to (3a);
    \end{pgfonlayer}
    
      \fill[\transitionbgrd, opacity=\transitionop]  (-0.5,-1) rectangle (0.5,1);
      \node (11) at (-0.5,1) {};
      \node (21) at (-0.5,-1) {};
      \node (31) at (0.5,-1) {};
      \node (41) at (0.5,1) {};
  \end{scope}

  \begin{scope}[
      yshift=8,every node/.append style={
      yslant=\ycoord,xslant=\xcoord},yslant=\ycoord,xslant=\xcoord
      ]
    \begin{pgfonlayer}{nodelayer}
      \node [transition] (2b) []  at (0,0) {};  
    \end{pgfonlayer}
    \begin{pgfonlayer}{edgelayer}
      \draw[thick, in=135, out=-45] (-0.5,0.5) to  (2b);
      \draw[thick,  in=225, out=45] (-0.5,-0.5) to  (2b);
      \draw[thick] (2b) to  (0.5,0);
    \end{pgfonlayer}
  \end{scope}

  \begin{scope}[
      yshift=28,every node/.append style={
      yslant=\ycoord,xslant=\xcoord},yslant=\ycoord,xslant=\xcoord
      ]
    \begin{pgfonlayer}{nodelayer}
      \node [transition] (2b) []  at (0,0) {};  
    \end{pgfonlayer} 
    \begin{pgfonlayer}{edgelayer}
      \draw[thick, in=225, out=-45] (-0.5,0.5) to  (2b);
      \draw[thick,  in=135, out=45] (-0.5,-0.5) to  (2b);
      \draw[thick] (2b) to  (0.5,0);
    \end{pgfonlayer}
  \end{scope}

  \begin{scope}[
      yshift=30,every node/.append style={
      yslant=\ycoord,xslant=\xcoord},yslant=\ycoord,xslant=\xcoord
      ]

      \fill[\transitionbgrd, opacity=\transitionop]  (-0.5,-1) rectangle (0.5,1);
      \node (14) at (-0.5,1) {};
      \node (24) at (-0.5,-1) {};
      \node (34) at (0.5,-1) {};
      \node (44) at (0.5,1) {};
  \end{scope}

    \fill[\transitionbgrd, opacity=\transitionop]   (11.center) -- (14.center) -- (24.center) -- (21.center) -- (11.center);
    \fill[\transitionbgrd, opacity=\transitionop]   (21.center) -- (24.center) -- (34.center) -- (31.center) -- (21.center);
    \fill[\transitionbgrd, opacity=\transitionop]   (31.center) -- (34.center) -- (44.center) -- (41.center) -- (31.center);
    \fill[\transitionbgrd, opacity=\transitionop]   (41.center) -- (44.center) -- (14.center) -- (11.center) -- (41.center);

\end{tikzpicture}
}  
  \]
  Thus, the \emph{transition classes} of $G_\pet(N)$ are the transitions of $N$.
\end{description}

\begin{rmk}
  Note that $F_\pre$ and $G_\pet$ both build a \SigmaNet{} whose transition classes are the transitions of a pre-net or Petri net.
  On the other hand, $G_\pre$ and $F_\pet$ are ``dual'', in that they build a pre-net or Petri net whose transitions are, respectively, the transitions or the transition classes of a \SigmaNet.
  In particular, the composite $F_\pet \circ F_\pre$ preserves transitions: it is the functor $\PreNet \to \Petri$ that simply forgets the ordering of inputs and outputs.
  Its right adjoint $G_\pre \circ G_\pet$ explodes each transition of a Petri net into as many transitions of a pre-net as possible, giving its inputs and outputs all possible orderings.
\end{rmk}

The last adjunction to construct is the one in the middle column:
\begin{thm}
\label{thm:SigmaNet_SSMC_adjunction}
  There is an adjunction
  \[
    \begin{tikzcd}[row sep=large,column sep=large]
      \SigmaNetCat \ar[r,bend left, "F_{\SigmaNetCat}"] \ar[r, phantom, "\bot"]  &
      \SSMC. \ar[l,bend left, "U_{\SigmaNetCat}"]
    \end{tikzcd}  
  \]
\end{thm}

Two proofs of \cref{thm:SigmaNet_SSMC_adjunction} were sketched by Kock~\cite[\S\S6--7]{Kock}.  Our proof is more similar to the proof of~\cite[Theorem 5.1]{Master}, which is a generalization of \cref{prop: adjunction Petri - CMC,prop: adjunction PreNet - StrMC} involving a Lawvere theory $\mathsf{Q}$: these two propositions follow by taking $\mathsf{Q}$ to be the theory of commutative monoids and the theory of monoids, respectively.   In that proof, an adjunction between $\mathsf{Q}$-nets and $\mathsf{Q}$-categories (i.e., $\mathsf{Q}$-algebras in $\mathsf{Cat}$) was obtained as the composite of two adjunctions where the intermediate category consists of $\mathsf{Q}$-graphs: graphs internal to the category of $\mathsf{Q}$-algebras.
Note that these $\mathsf{Q}$-graphs have operations coming from the Lawvere theory $\mathsf{Q}$, which act both on vertices and edges, but they lack the ability to compose edges (i.e., morphisms) that one has in a $\mathsf{Q}$-category.

Our desired adjunction here is not a special case of~\cite[Theorem 5.1]{Master}, since the symmetries in a symmetric monoidal category cannot be represented by a structure on the object set alone.
However, we can perform a similar factorization through a category containing only the monoidal operations.  We begin by reducing the problem from strict symmetric monoidal categories to (colored) props.

\begin{defn}
  A \defin{(colored) prop} consists of a set $S$, a strict symmetric monoidal category $B$, and a strict symmetric monoidal functor $i \maps PS \to B$ that is bijective on objects.   A morphism of props consists of a function $S\to S'$ and a strict symmetric monoidal functor $B\to B'$ making the evident square commute.  We denote the category of props by $\PROP$.
\end{defn}

\begin{lem}
\label{lem:prop-adj}
  There is an adjunction
  \[
    \begin{tikzcd}[row sep=large,column sep=large]
      \PROP \ar[r,bend left, "F_2"] \ar[r, phantom, "\bot"]  &
      \SSMC. \ar[l,bend left, "U_2"]
    \end{tikzcd}  
  \]
\end{lem}

\begin{proof}
  For a prop $(S,B,i)$ we define $F_2(S,B,i) = B$.  And for a strict symmetric monoidal category $B$, we let $S$ be the set of objects of $B$, so that we have a strict symmetric monoidal functor $PS \to B$.
  Now we factor this functor as a bijective-on-objects functor $i \maps PS \to B'$ followed by a fully faithful one $p \maps B' \to B$.
  Then $B'$ can be given a symmetric strict monoidal structure making both $i$ and $p$ strict symmetric monoidal functors, and we define $U_2(B) = (S,B',i)$.
\end{proof}

Therefore, it will suffice to construct an adjunction between $\SigmaNetCat$ and $\PROP$.
We work with the profunctor representation of \SigmaNet{}s.
Let $U_1 \maps \PROP \to \SigmaNetCat$ be the functor sending $(S,B,i)$ to $(S,N)$ where $N(a_1,a_2) = \hom_B(i(a_2),i(a_1))$.
This is the functor we aim to construct a left adjoint of.  As in~\cite{Master}, we do this ``fiberwise'' for a fixed $A$, then piece the fiberwise adjunctions together.

\begin{lem}
\label{lem:U2-cart}
  We have a commutative triangle
  \[
    \begin{tikzcd}[column sep=0.2cm]
      \SigmaNetCat \ar[dr] & & \PROP \ar[dl] \ar[ll,"U_1"'] \\
      & \Set
    \end{tikzcd}
  \]
  in which the two diagonal functors are split fibrations and $U_1$ is cartesian.
\end{lem}

\begin{proof}
  The two diagonal functors send $(S,N)$ to $S$ and $(S,B,i)$ to $S$, respectively.
  To show the left-hand diagonal functor is a split fibration, let $(S,H)\in \SigmaNetCat$ and $g:S'\to S$; then $(S', N\circ (Pg\times Pg\op))$ is the domain of a cartesian lifting.
  For the right-hand functor, given $(S,B,i)\in \PROP$ and $g \maps S'\to S$, the composite $i\circ g \maps PS' \to B$ may no longer be bijective on objects, but we can factor it as a bijective-on-objects functor $i' \maps PS' \to B'$ followed by a fully faithful one $g' \maps B' \to B$.  These are both again strict symmetric monoidal functors, and the induced map $(S',B',i') \to (S,B,i)$ is cartesian.   Finally, $U_1$ is cartesian by construction, since $g'$ is fully faithful.
\end{proof}

Let $U_{1,S} \maps \SigmaNetCat_S \to \PROP_S$ denote the restriction of $U_1$ to the fibers over a particular set $S$. We will construct a left adjoint $F_{1,S}$ of this functor, then piece these together fiberwise.

Following the proof of~\cite{Master}, we need to decompose the structure of a prop with object set $S$ into the ``monoidal piece'' and the ``composition piece''.
This can be accomplished as follows.
Batanin and Markl~\cite{BataninMarkl} define a \textbf{duoidal category} to be a category $\C$ with two monoidal structures $(\star,J)$ and $(\diamond,I)$ and additional natural morphisms
\begin{gather*}
  I\to J \qquad I\to I\star I \qquad J\diamond J \to J \\
  (A\star B) \diamond (C\star D) \to (A \diamond C) \star (B\diamond D)
\end{gather*}
satisfying axioms that say $(\star,J)$ is a pseudomonoid structure on $(\C,\diamond,I)$ in the 2-category of lax monoidal functors.
It is \textbf{$\diamond$-symmetric} if $\diamond$ is a symmetric monoidal structure and the above maps commute with the symmetry in an evident way.

In a duoidal category, the monoidal structure $\star$ lifts to a monoidal structure on the category of $\diamond$-monoids.
A $\star$-monoid in this monoidal category of $\diamond$-monoids is called a \textbf{duoid}.
Similarly, if the duoidal category is $\diamond$-symmetric, then $\star$ lifts to the category of commutative $\diamond$-monoids, and a monoid therein is called a \textbf{$\diamond$-commutative duoid}.

\begin{lem}\label{lem:ssmcpair-duoidal}
  There is a $\diamond$-symmetric duoidal structure on $\SigmaNetCat_S$ whose category of $\diamond$-commutative duoids is equivalent to $\PROP_S$.
\end{lem}
\begin{proof}
  Note that $\SigmaNetCat = \mathbf{Prof}(PS,PS)$ is the hom-category of $PS$ in the bicategory $\mathbf{Prof}$ of categories and profunctors.
  Since this is an endo-hom-category in a bicategory, it has a monoidal structure given by composition in $\mathbf{Prof}$, which we call $\star$ (thus $J$ is the hom-functor of $A$).   The monoidal structure $\diamond$ is given by convolution:

  {\iflics\small\fi
  \[ (H\diamond K)(x,z) = \int^{a,b,c,d} \! PS(x,ab) \times H(a,c) \times K(b,d) \times PS(cd, z). \]}%
  with $I(x,y) = PS(x,\epsilon) \times PS(\epsilon,y)$.

  A $\star$-monoid is a monad on $PS$ in the bicategory $\mathbf{Prof}$, which is well-known to be equivalent to a category $B$ with a bijective-on-objects functor $PS\to B$.
  Applying the Yoneda lemma, we find that a $\diamond$-monoid structure on such a $B$ consists of morphisms
  \[ B(a,c) \times B(b,d) \to B(a b, c d) \]
  that are suitably compatible.
  This extends the monoidal structure of $PS$ to the arrows of $B$ (it is already defined on the objects of $B$ since they are the same as the objects of $A$).
  Compatibility with the duoidal exchange morphism says that this action is functorial, while compatibility with the map $J\diamond J \to J$ says that it extends the functorial action of the monoidal structure on $PS$.
  The associativity and unitality of a $\diamond$-monoid says $B$ has a strict monoidal structure and the functor $PS\to B$ is strict monoidal.
  Finally, the symmetry of $\diamond$ switches $H$ and $K$ and composes with the symmetry isomorphisms in $PS$ on either side; thus $\diamond$-commutativity of a duoid makes $B$ a \emph{symmetric} strict monoidal category and $PS\to B$ a strict symmetric monoidal functor.
\end{proof}

  In fact, an analogous result holds with $PS$ replaced by any symmetric monoidal category.
  A more abstract construction of this duoidal structure was given by Garner and L\'opez Franco~\cite[Proposition 51]{GarnerLopezFranco}, while the identification of its duoids follows from their Proposition 49 and the remarks after Proposition 54.   Note that the adjective ``commutative'' in~\cite{GarnerLopezFranco} is used with a different meaning than ours; we repeat that for us, ``$\diamond$-commutative'' simply means that the monoid structure with respect to $\diamond$ is commutative in the ordinary sense for a monoid object in a symmetric monoidal category.

Note that both monoidal structures $\star$ and $\diamond$ of $\SigmaNetCat_S$ preserve colimits in each variable.
We can now work at a higher level of abstraction.

\begin{lem}\label{lem:free-duoids}
  For any cocomplete $\diamond$-symmetric duoidal category $\C$ such that $\star$ and $\diamond$ preserve colimits in each variable, the forgetful functor
  \[\diamond \maps \CommDuoid(\C) \to \C\]
  has a left adjoint.
\end{lem}
\begin{proof}
  Recall that free commutative monoids exist in any cocomplete monoidal category whose tensor product preserves colimits in each variable, given by
  \[ F X = \coprod_n X^{\diamond n} / \Sigma_n \]
  where $X^{\diamond n} / \Sigma_n$ denotes the $n^{\mathrm{th}}$ tensor power of $X$ quotiented by the action of the $n^{\mathrm{th}}$ symmetric group.
  Indeed, commutative monoids are monadic over such a category.
  Thus, the category of commutative $\diamond$-monoids in our $\C$ is monadic over $\C$.

  Moreover, since $\diamond$ preserves colimits in each variable, by standard arguments it preserves reflexive coequalizers and sequential colimits in both variables together.
  Thus $X \mapsto X^{\diamond n}$ also preserves reflexive coequalizers and sequential colimits, hence so does the functor $F$ and thus the monad for commutative $\diamond$-monoids.
  It follows that reflexive coequalizers and sequential colimits in the category of commutative $\diamond$-monoids are computed as in $\C$, and therefore are preserved in each variable by the lifted tensor product $\star$.
  Therefore, by~\cite{LackMonoids}, the free $\star$-monoid on a commutative $\diamond$-monoid exists.
  Composing these two free constructions, we find that free $\diamond$-commutative duoids exist.
\end{proof}

\begin{lem}\label{lem:ssmcadj2}
  There is an adjunction
  \[
    \begin{tikzcd}[row sep=large,column sep=large]
      \SigmaNetCat \ar[r,bend left, "F_{1}"] \ar[r, phantom, "\bot"]  &
      \PROP. \ar[l,bend left, "U_{1}"]
    \end{tikzcd}  
  \]
\end{lem}
\begin{proof}
  By \cref{lem:ssmcpair-duoidal,lem:free-duoids}, each fiber functor $U_{1,S}$ has a left adjoint $F_{1,S}$; thus it remains to piece these adjoints together.
  Suppose $(S',N) \in \SigmaNetCat$ and $(S,B,i) \in \PROP$.
  By \cref{lem:U2-cart}, a morphism $(S',N) \to U_1 (S,B,i)$ is equivalently given by a function $g \maps S'\to S$ and a morphism $(S',N) \to U_{1,S'}(S',B', i')$ in $\SigmaNetCat_{S'}$, where $P S' \xrightarrow{i'} B' \xrightarrow{g'} B$ is the factorization of $i\circ P g$ as a bijective-on-objects functor followed by a fully faithful one.
  But the latter morphism is equivalently a morphism $F_{1,S'}(S',N) \to (S',B', i')$ in $\PROP_{S'}$, hence a morphism $F_{1,S'}(S',N) \to (S,B,i)$ in $\PROP$.
  Thus, defining $F_1(S',N) = F_{1,S'}(S',N)$ yields a left adjoint to $U_1$.
  (Note that it is unnecessary to ask whether $F_1$ is cartesian.)
\end{proof}

\begin{proof}[Proof of \cref{thm:SigmaNet_SSMC_adjunction}]
  Combining Lemma \ref{lem:prop-adj} and \ref{lem:ssmcadj2} we obtain the composite adjunction:
  \[
    \begin{tikzcd}[row sep=large,column sep=large]
      \SigmaNetCat \ar[r,bend left, "F_{1}"] \ar[r, phantom, "\bot"]  &
      \PROP \ar[r,bend left, "F_{2}"] \ar[r, phantom, "\bot"] \ar[l,bend left, "U_{1}"] &
      \SSMC. \ar[l,bend left, "U_{2}"]
    \end{tikzcd}
    \qedhere
  \]
\end{proof}

We end this section by considering the commutativity properties of the squares in \cref{eq:adjunctions}.

\begin{prop}\label{thm:left-square}
  There is a natural isomorphism $G_\pre \circ U_{\SigmaNetCat} \cong U_{\PreNet} \circ U_{\StrMC}$.
  Therefore, there is also a natural isomorphism $F_{\StrMC} \circ F_{\PreNet} \cong F_{\SigmaNetCat} \circ F_{\pre}$.
\end{prop}
\begin{proof}
  With our precise definitions, the first isomorphism is actually a strict equality: both functors take a symmetric strict monoidal category $C$ to the pre-net whose places are the objects of $C$ and whose transitions from a word $p$ to a word $q$ are the morphisms in $C$ from the tensor product of $p$ to the tensor product of $q$.  The second isomorphism follows by passage to left adjoints.
\end{proof}

Recalling from \cref{sec:pre-nets} that the composite $F_{\StrMC} \circ F_{\PreNet}$ has been used to give a categorical semantics for pre-nets, we see that this semantics factors through \SigmaNet{}s.

\begin{prop}\label{thm:right-square-1}
  There is a natural isomorphism $G_{\pet} \circ U_{\Petri} \cong U_{\SigmaNetCat} \circ U_{\SSMC}$.
  Therefore, there is also a natural isomorphism $F_{\SSMC} \circ F_{\SigmaNetCat} \cong F_{\Petri} \circ F_{\pet}$.
\end{prop}
\begin{proof}
  Again, the first isomorphism is a strict equality: both functors take a commutative monoidal category $C$ to the \SigmaNet{} whose places are the objects of $C$ and whose transitions from $p$ to $q$ are the morphisms in $C$ from the tensor product of $p$ to the tensor product of $q$, with symmetries acting trivially.
  The second isomorphism follows by passage to left adjoints.
\end{proof}

Though analogous to \cref{thm:left-square}, \cref{thm:right-square-1} does not imply that the categorical semantics of Petri nets factors through \SigmaNet{}s.
However, that is also true:

\begin{prop}
  There is a natural isomorphism $F_{\Petri}\cong F_{\SSMC} \circ F_{\SigmaNetCat} \circ G_{\pet}$.
\end{prop}
\begin{proof}
  Let $N$ be a Petri net and $C$ a commutative monoidal category; since $G_{\pet}$ is fully faithful we have natural isomorphisms
  \begin{align*}
    \Petri(N, U_{\Petri}(C))
    &\cong \SigmaNetCat(G_\pet(N), G_\pet (U_{\Petri}(C)))\\
    &\cong \SigmaNetCat(G_\pet(N), U_{\SigmaNetCat} (U_{\SSMC}(C)))\\
    &\cong \SSMC(F_{\SigmaNetCat} (G_\pet(N)), U_{\SSMC}(C))\\
    &\cong \CMC(F_{\SSMC} (F_{\SigmaNetCat} (G_\pet(N))), C).
  \end{align*}
  Thus $F_{\SSMC} \circ F_{\SigmaNetCat} \circ G_{\pet}$ is left adjoint to $U_{\Petri}$, hence isomorphic to $F_{\Petri}$.
\end{proof}


\section{Relation to whole-grain Petri nets}
\label{sec:whole-grain}

We now clarify the relation of our work to Kock's ``whole-grain Petri nets'' \cite{Kock}.  We show that a whole-grain Petri net can be thought of as a special sort of \SigmaNet{}: one that is free on a pre-net.   We first recall Kock's definition:

\begin{defn}
  A \textbf{whole-grain Petri net} is a diagram
  \[
    \begin{tikzcd}
      S & I \ar[l] \ar[r] & T & O \ar[l] \ar[r] & S
    \end{tikzcd}
  \]
  in which the fibers of the functions $I\to T$ and $O\to T$ are finite.
  A morphism of whole-grain Petri nets, sometimes called an \emph{etale map}, is a diagram
  \[
    \begin{tikzcd}
      S \ar[d] & I \ar[l] \ar[r]\ar[d] \ar[dr,phantom,near start,"\lrcorner"] & T \ar[d] & O \ar[l] \ar[r]\ar[d] \ar[dl,phantom,near start,"\llcorner"] & S\ar[d]\\
      S' & I' \ar[l] \ar[r] & T' & O' \ar[l] \ar[r] & S'
    \end{tikzcd}
  \]
  This defines the category $\WGPet$.
\end{defn}

\begin{thm}\label{thm:wgpet}
  The category $\WGPet$ is equivalent to the full image of $F_\pre \maps \PreNet\to \SigmaNetCat$.
  In other words, there are functors
  \[
    \begin{tikzcd}
      \PreNet \ar[r,"Z_1"] & \WGPet \ar[r,"Z_2"] & \SigmaNetCat
    \end{tikzcd}
  \]
  such that $Z_1$ is essentially surjective, $Z_2$ is fully faithful, and the composite $Z_2 \circ Z_1$ is isomorphic to $F_\pre$.
\end{thm}
\begin{proof}
  Given a pre-net $s,t:T \to S^* \times S^*$, let $I$ be the set of transitions $u\in T$ equipped with a choice of an element of $s(u)$, and define $O$ similarly using $t(u)$.
  There are forgetful functions $I\to T$ and $O\to T$, and maps $I\to S$ (resp.\ $O\to S$) that select the chosen element of $s(u)$ (resp.\ $t(u)$).
  This defines a whole-grain Petri net $Z_1(s,t)$.
  Note that the fibers of $I\to T$ and $O\to T$ are not just finite but equipped with a linear ordering, and the morphisms in the image of $Z_1$ (which is faithful) are precisely those that preserve these orderings.

  To see that $Z_1$ is essentially surjective, given a whole-grain Petri net $N$ we choose linear orderings on each fiber of the maps $I\to T$ and $O\to T$.
  These orderings associated each element of $T$ to two elements of $S^*$, yielding a pre-net whose image under $Z_1$ is isomorphic to $N$.

  Now, given a whole-grain Petri net $S \leftarrow I \to T \leftarrow O \to S$, we define a \SigmaNet{} in the presheaf perspective.
  Its set of places is $S$, and its $(m,n)$-transitions are elements $u\in T$ equipped with a linear ordering on the fibers of $I$ and $O$ over $u$, which we require to have $m$ and $n$ elements respectively.
  These linear orderings enable us to define the source and target maps picking out places, while the permutations act on the linear orderings.
  This defines the functor $Z_2$.
  
  Note that the set $T$ in a whole-grain Petri net $N$ is naturally isomorphic to the set of transition \emph{classes} of $Z_2(N)$.
  Thus, to show that $Z_2$ is fully faithful it remains to show that a morphism $\alpha:Z_2(N) \to Z_2(N')$ uniquely determines the maps $I\to I'$ and $O\to O'$.
  Given $i\in I$ lying over $t\in T$, choose any ordering on the fibers over $t$, in which $i$ appears as the $k^{\mathrm{th}}$ element of its fiber.
  This choice determines a transition $\hat{t}$ of $Z_2(N)$, and hence a transition $\alpha(\hat{t})$ of $Z_2(N')$, which is an element $\check{\alpha}(\hat{t})$ of $T'$ with ordered fibers.
  Then the function $I\to I'$ can and must send $i$ to the $k^{\mathrm{th}}$ element of the $I$-fiber over $\check{\alpha}(\hat{t})$.
  This is independent of the choice of ordering because $\alpha$ commutes with the permutation actions, and it is straightforward to check that it indeed defines a morphism $N\to N'$.

  Finally, the composite $Z_2 \circ Z_1$ preserves the set of places and replaces each $(m,n)$-transition by $m!\, n!$ transitions with free permutation action; but this is the same as $F_\pre$.
\end{proof}

Another construction of the functor $Z_2$ appears in~\cite{Kock}, as a restricted Yoneda embedding or ``nerve''.
Recall the categories $C$ and $D$ from \cref{prop:topos of pre-nets,thm: sigma nets is a topos}.
In fact $D$ is the full image of the composite of the Yoneda embedding $C \hookrightarrow [C\op,\Set] \simeq \PreNet$ with $Z_1:\PreNet \to \WGPet$; we can then define $Z_2$ as the composite $\WGPet \to [\WGPet\op,\Set] \to [D\op,\Set] \simeq \SigmaNetCat$.
In addition, \cref{thm:wgpet} can be viewed as a two-sided version of the relationship between symmetric and nonsymmetric collections, as in~\cite[2.4.4]{Kock3}.

\iflics\else
\section{Open nets}
\label{sec:open-nets}

Various kinds of ``open'' nets have been proposed, which allow one to build nets by gluing together smaller open nets \cite{BBGM, BCEH, BMMS2013, SassoneCongruence}.  \iflics In particular, there is \else In earlier work we introduced \fi a symmetric monoidal double category of open Petri nets, where composing open Petri nets is done by identifying places \cite{BM}. 
For example, \iflics consider the following nets:
\[
\scalebox{0.6}{
\begin{tikzpicture}[baseline=0]
	\begin{pgfonlayer}{nodelayer}
		\node [place,tokens=0] (A) at (-4, 0.5) {$ $};
		\node [place,tokens=0] (B) at (-4, -0.5) {$ $};
		\node [place,tokens=0] (C) at (-1, 0.5) {$ $};
		\node [place,tokens=0] (D) at (-1, -0.5) {$ $};
             \node [style=transition] (a) at (-2.5, 0) {$ $}; 
		
		\node [style=empty] (X) at (-5.1, 1) {$X$};
		\node [style=none] (Xtr) at (-4.75, 0.75) {};
		\node [style=none] (Xbr) at (-4.75, -0.75) {};
		\node [style=none] (Xtl) at (-5.4, 0.75) {};
             \node [style=none] (Xbl) at (-5.4, -0.75) {};
	
		\node [style=inputdot] (1) at (-5, 0.5) {};
		\node [style=empty] at (-5.2, 0.5) {$1$};
		\node [style=inputdot] (2) at (-5, 0) {};
		\node [style=empty] at (-5.2, 0) {$2$};
		\node [style=inputdot] (3) at (-5, -0.5) {};
		\node [style=empty] at (-5.2, -0.5) {$3$};

		\node [style=empty] (Y) at (0.1, 1) {$Y$};
		\node [style=none] (Ytr) at (.4, 0.75) {};
		\node [style=none] (Ytl) at (-.25, 0.75) {};
		\node [style=none] (Ybr) at (.4, -0.75) {};
		\node [style=none] (Ybl) at (-.25, -0.75) {};

		\node [style=inputdot] (4) at (0, 0.5) {};
		\node [style=empty] at (0.2, 0.5) {$4$};
		\node [style=inputdot] (5) at (0, -0.5) {};
		\node [style=empty] at (0.2, -0.5) {$5$};		
		
	\end{pgfonlayer}
	\begin{pgfonlayer}{edgelayer}
		\draw [style=inarrow] (A) to (a);
		\draw [style=inarrow] (B) to (a);
		\draw [style=inarrow] (a) to (C);
		\draw [style=inarrow] (a) to (D);
		\draw [style=inputarrow] (1) to (A);
		\draw [style=inputarrow] (2) to (B);
		\draw [style=inputarrow] (3) to (B);
		\draw [style=inputarrow] (4) to (C);
		\draw [style=inputarrow] (5) to (D);
		\draw [style=simple] (Xtl.center) to (Xtr.center);
		\draw [style=simple] (Xtr.center) to (Xbr.center);
		\draw [style=simple] (Xbr.center) to (Xbl.center);
		\draw [style=simple] (Xbl.center) to (Xtl.center);
		\draw [style=simple] (Ytl.center) to (Ytr.center);
		\draw [style=simple] (Ytr.center) to (Ybr.center);
		\draw [style=simple] (Ybr.center) to (Ybl.center);
		\draw [style=simple] (Ybl.center) to (Ytl.center);
	\end{pgfonlayer}
\end{tikzpicture}
}
\qquad
\scalebox{0.6}{
\begin{tikzpicture}[baseline=0]
	\begin{pgfonlayer}{nodelayer}

		\node [style = transition] (b) at (2.5, 1) {$ $};
		\node [style = transition] (c) at (2.5, -1) {$ $};
		\node [place,tokens=0] (E) at (1, 0) {$ $};
		\node [place,tokens=0] (F) at (4,0) {$ $};

	    \node [style=empty] (Y) at (-0.1, 1) {$Y$};
		\node [style=none] (Ytr) at (.25, 0.75) {};
		\node [style=none] (Ytl) at (-.4, 0.75) {};
		\node [style=none] (Ybr) at (.25, -0.75) {};
		\node [style=none] (Ybl) at (-.4, -0.75) {};

		\node [style=inputdot] (4) at (0, 0.5) {};
		\node [style=empty] at (-0.2, 0.5) {$4$};
		\node [style=inputdot] (5) at (0, -0.5) {};
		\node [style=empty] at (-0.2, -0.5) {$5$};		
		
		\node [style=empty] (Z) at (5, 1) {$Z$};
		\node [style=none] (Ztr) at (4.75, 0.75) {};
		\node [style=none] (Ztl) at (5.4, 0.75) {};
		\node [style=none] (Zbl) at (5.4, -0.75) {};
		\node [style=none] (Zbr) at (4.75, -0.75) {};

		\node [style=inputdot] (6) at (5, 0) {};
		\node [style=empty] at (5.2, 0) {$6$};	
		
	\end{pgfonlayer}
	\begin{pgfonlayer}{edgelayer}
		\draw [style=inarrow, bend left=30, looseness=1.00] (E) to (b);
		\draw [style=inarrow, bend left=30, looseness=1.00] (b) to (F);
		\draw [style=inarrow, bend left=30, looseness=1.00] (c) to (E);
		\draw [style=inarrow, bend left=30, looseness=1.00] (F) to (c);
		\draw [style=inputarrow] (4) to (E);
		\draw [style=inputarrow] (5) to (E);
		\draw [style=inputarrow] (6) to (F);
		\draw [style=simple] (Ytl.center) to (Ytr.center);
		\draw [style=simple] (Ytr.center) to (Ybr.center);
		\draw [style=simple] (Ybr.center) to (Ybl.center);
		\draw [style=simple] (Ybl.center) to (Ytl.center);
		\draw [style=simple] (Ztl.center) to (Ztr.center);
		\draw [style=simple] (Ztr.center) to (Zbr.center);
		\draw [style=simple] (Zbr.center) to (Zbl.center);
		\draw [style=simple] (Zbl.center) to (Ztl.center);
	\end{pgfonlayer}
\end{tikzpicture}
}
\]
In addition to a Petri net, they consists of sets $X$, $Y$, $Z$ and arbitrary functions from these sets into their set of places. These indicate places at which tokens could flow in or out. We may write these open Petri net as $P \maps X \to Y$ and $Q \maps Y \to Z$ for short.
We can compose $P$ with $Q$ and obtain the following open Petri net $Q \circ P \maps X \to Z$:
\else here is an open Petri net with one transition:
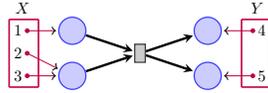
\[
\scalebox{0.6}{
\begin{tikzpicture}
	\begin{pgfonlayer}{nodelayer}
		\node [place,tokens=0] (A) at (-4, 0.5) {$ $};
		\node [place,tokens=0] (B) at (-4, -0.5) {$ $};
		\node [place,tokens=0] (C) at (-1, 0.5) {$ $};
		\node [place,tokens=0] (D) at (-1, -0.5) {$ $};
             \node [style=transition] (a) at (-2.5, 0) {$ $}; 
		
		\node [style=empty] (X) at (-5.1, 1) {$X$};
		\node [style=none] (Xtr) at (-4.75, 0.75) {};
		\node [style=none] (Xbr) at (-4.75, -0.75) {};
		\node [style=none] (Xtl) at (-5.4, 0.75) {};
             \node [style=none] (Xbl) at (-5.4, -0.75) {};
	
		\node [style=inputdot] (1) at (-5, 0.5) {};
		\node [style=empty] at (-5.2, 0.5) {$1$};
		\node [style=inputdot] (2) at (-5, 0) {};
		\node [style=empty] at (-5.2, 0) {$2$};
		\node [style=inputdot] (3) at (-5, -0.5) {};
		\node [style=empty] at (-5.2, -0.5) {$3$};

		\node [style=empty] (Y) at (0.1, 1) {$Y$};
		\node [style=none] (Ytr) at (.4, 0.75) {};
		\node [style=none] (Ytl) at (-.25, 0.75) {};
		\node [style=none] (Ybr) at (.4, -0.75) {};
		\node [style=none] (Ybl) at (-.25, -0.75) {};

		\node [style=inputdot] (4) at (0, 0.5) {};
		\node [style=empty] at (0.2, 0.5) {$4$};
		\node [style=inputdot] (5) at (0, -0.5) {};
		\node [style=empty] at (0.2, -0.5) {$5$};		
		
	\end{pgfonlayer}
	\begin{pgfonlayer}{edgelayer}
		\draw [style=inarrow] (A) to (a);
		\draw [style=inarrow] (B) to (a);
		\draw [style=inarrow] (a) to (C);
		\draw [style=inarrow] (a) to (D);
		\draw [style=inputarrow] (1) to (A);
		\draw [style=inputarrow] (2) to (B);
		\draw [style=inputarrow] (3) to (B);
		\draw [style=inputarrow] (4) to (C);
		\draw [style=inputarrow] (5) to (D);
		\draw [style=simple] (Xtl.center) to (Xtr.center);
		\draw [style=simple] (Xtr.center) to (Xbr.center);
		\draw [style=simple] (Xbr.center) to (Xbl.center);
		\draw [style=simple] (Xbl.center) to (Xtl.center);
		\draw [style=simple] (Ytl.center) to (Ytr.center);
		\draw [style=simple] (Ytr.center) to (Ybr.center);
		\draw [style=simple] (Ybr.center) to (Ybl.center);
		\draw [style=simple] (Ybl.center) to (Ytl.center);
	\end{pgfonlayer}
\end{tikzpicture}
}
\]
In addition to a Petri net, it consists of sets $X$ and $Y$ and arbitrary functions from these sets into the set of places.  These indicate places at which tokens could flow in or out.   We may write this open Petri net as $P \maps X \to Y$ for short.

Given another open Petri net $Q \maps Y \to Z$:
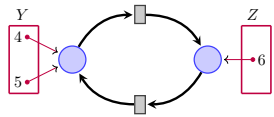
\[
\scalebox{0.6}{
\begin{tikzpicture}
	\begin{pgfonlayer}{nodelayer}

		\node [style = transition] (b) at (2.5, 1) {$ $};
		\node [style = transition] (c) at (2.5, -1) {$ $};
		\node [place,tokens=0] (E) at (1, 0) {$ $};
		\node [place,tokens=0] (F) at (4,0) {$ $};

	    \node [style=empty] (Y) at (-0.1, 1) {$Y$};
		\node [style=none] (Ytr) at (.25, 0.75) {};
		\node [style=none] (Ytl) at (-.4, 0.75) {};
		\node [style=none] (Ybr) at (.25, -0.75) {};
		\node [style=none] (Ybl) at (-.4, -0.75) {};

		\node [style=inputdot] (4) at (0, 0.5) {};
		\node [style=empty] at (-0.2, 0.5) {$4$};
		\node [style=inputdot] (5) at (0, -0.5) {};
		\node [style=empty] at (-0.2, -0.5) {$5$};		
		
		\node [style=empty] (Z) at (5, 1) {$Z$};
		\node [style=none] (Ztr) at (4.75, 0.75) {};
		\node [style=none] (Ztl) at (5.4, 0.75) {};
		\node [style=none] (Zbl) at (5.4, -0.75) {};
		\node [style=none] (Zbr) at (4.75, -0.75) {};

		\node [style=inputdot] (6) at (5, 0) {};
		\node [style=empty] at (5.2, 0) {$6$};	
		
	\end{pgfonlayer}
	\begin{pgfonlayer}{edgelayer}
		\draw [style=inarrow, bend left=30, looseness=1.00] (E) to (b);
		\draw [style=inarrow, bend left=30, looseness=1.00] (b) to (F);
		\draw [style=inarrow, bend left=30, looseness=1.00] (c) to (E);
		\draw [style=inarrow, bend left=30, looseness=1.00] (F) to (c);
		\draw [style=inputarrow] (4) to (E);
		\draw [style=inputarrow] (5) to (E);
		\draw [style=inputarrow] (6) to (F);
		\draw [style=simple] (Ytl.center) to (Ytr.center);
		\draw [style=simple] (Ytr.center) to (Ybr.center);
		\draw [style=simple] (Ybr.center) to (Ybl.center);
		\draw [style=simple] (Ybl.center) to (Ytl.center);
		\draw [style=simple] (Ztl.center) to (Ztr.center);
		\draw [style=simple] (Ztr.center) to (Zbr.center);
		\draw [style=simple] (Zbr.center) to (Zbl.center);
		\draw [style=simple] (Zbl.center) to (Ztl.center);
	\end{pgfonlayer}
\end{tikzpicture}
}
\]
we can compose it with $P$ and obtain the following open Petri net $Q \circ P \maps X \to Z$:\fi
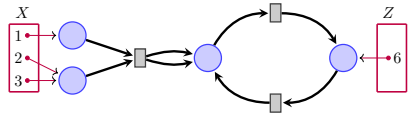
\[
\scalebox{0.6}{
\begin{tikzpicture}
	\begin{pgfonlayer}{nodelayer}
		\node [place,tokens=0] (A) at (-4, 0.5) {$$};
		\node [place,tokens=0] (B) at (-4, -0.5) {$$};;
             \node [style=transition] (a) at (-2.5, 0) {$ $}; 
		\node [place,tokens=0] (E) at (-1, 0) {$$};
		\node [place,tokens=0] (F) at (2,0) {$ $};

	     \node [style = transition] (b) at (.5, 1) {$ $};
		\node [style = transition] (c) at (.5, -1) {$ $};
		
		\node [style=empty] (X) at (-5.1, 1) {$X$};
		\node [style=none] (Xtr) at (-4.75, 0.75) {};
		\node [style=none] (Xbr) at (-4.75, -0.75) {};
		\node [style=none] (Xtl) at (-5.4, 0.75) {};
             \node [style=none] (Xbl) at (-5.4, -0.75) {};
	
		\node [style=inputdot] (1) at (-5, 0.5) {};
		\node [style=empty] at (-5.2, 0.5) {$1$};
		\node [style=inputdot] (2) at (-5, 0) {};
		\node [style=empty] at (-5.2, 0) {$2$};
		\node [style=inputdot] (3) at (-5, -0.5) {};
		\node [style=empty] at (-5.2, -0.5) {$3$};	
		
		\node [style=empty] (Z) at (3, 1) {$Z$};
		\node [style=none] (Ztr) at (2.75, 0.75) {};
		\node [style=none] (Ztl) at (3.4, 0.75) {};
		\node [style=none] (Zbl) at (3.4, -0.75) {};
		\node [style=none] (Zbr) at (2.75, -0.75) {};

		\node [style=inputdot] (6) at (3, 0) {};
		\node [style=empty] at (3.2, 0) {$6$};	
		
	\end{pgfonlayer}
	\begin{pgfonlayer}{edgelayer}
		\draw [style=inarrow] (A) to (a);
		\draw [style=inarrow] (B) to (a);
	     \draw [style=inarrow, bend right=15, looseness=1.00] (a) to (E);
	     \draw [style=inarrow, bend left =15, looseness=1.00] (a) to (E);	
	     	\draw [style=inarrow, bend left=30, looseness=1.00] (E) to (b);
		\draw [style=inarrow, bend left=30, looseness=1.00] (b) to (F);
		\draw [style=inarrow, bend left=30, looseness=1.00] (c) to (E);
		\draw [style=inarrow, bend left=30, looseness=1.00] (F) to (c);	
		\draw [style=inputarrow] (1) to (A);
		\draw [style=inputarrow] (2) to (B);
		\draw [style=inputarrow] (3) to (B);
		\draw [style=inputarrow] (6) to (F);
		\draw [style=simple] (Xtl.center) to (Xtr.center);
		\draw [style=simple] (Xtr.center) to (Xbr.center);
		\draw [style=simple] (Xbr.center) to (Xbl.center);
		\draw [style=simple] (Xbl.center) to (Xtl.center);
		\draw [style=simple] (Ztl.center) to (Ztr.center);
		\draw [style=simple] (Ztr.center) to (Zbr.center);
		\draw [style=simple] (Zbr.center) to (Zbl.center);
		\draw [style=simple] (Zbl.center) to (Ztl.center);
	\end{pgfonlayer}
\end{tikzpicture}
}
\]

\iflics
We can also tensor open Petri nets, putting them side by side ``in parallel''.  This might suggest that open Petri nets should be the morphisms of a symmetric monoidal category.  However, composition of open Petri is not strictly associative, and there are very interesting maps \emph{between} open Petri nets.  To get a feeling for these, there is a morphism between the following nets:
\[
\scalebox{0.6}{
\begin{tikzpicture}
	\begin{pgfonlayer}{nodelayer}

		\node [style = transition] (a) at (2.5, 0.5) {$\alpha$};
		\node [style = transition] (a') at (2.5, -0.5) {$\alpha'$};
		\node [place,tokens=0] (A) at (1, 0.5) {$A$};
		\node [place,tokens=0] (A') at (1, -0.5) {$A'$};
		\node [place,tokens=0] (B) at (4,0) {$B$};
	
		\node [style=empty] (Y) at (-0.1, 1) {$X_1$};
		\node [style=none] (Ytr) at (.25, 0.75) {};
		\node [style=none] (Ytl) at (-.4, 0.75) {};
		\node [style=none] (Ybr) at (.25, -0.75) {};
		\node [style=none] (Ybl) at (-.4, -0.75) {};

		\node [style=inputdot] (4) at (0, 0.5) {};
		\node [style=empty] at (-0.2, 0.5) {$1$};
		\node [style=inputdot] (5) at (0, -0.5) {};
		\node [style=empty] at (-0.2, -0.5) {$1'$};		
		
		\node [style=empty] (Z) at (5.1, 1) {$Y_1$};
		\node [style=none] (Ztr) at (4.75, 0.75) {};
		\node [style=none] (Ztl) at (5.4, 0.75) {};
		\node [style=none] (Zbl) at (5.4, -0.75) {};
		\node [style=none] (Zbr) at (4.75, -0.75) {};

		\node [style=inputdot] (6) at (5, 0) {};
		\node [style=empty] at (5.2, 0) {$2$};	
		
	\end{pgfonlayer}
	\begin{pgfonlayer}{edgelayer}

		\draw [style=inarrow] (A) to (a);
		\draw [style=inarrow, bend left=20, looseness=1.00] (a) to (B);
		\draw [style=inarrow] (A') to (a');
		\draw [style=inarrow, bend right=20, looseness=1.00] (a') to (B);
		\draw [style=inputarrow] (4) to (A);
		\draw [style=inputarrow] (5) to (A');
		\draw [style=inputarrow] (6) to (B);

		\draw [style=simple] (Ytl.center) to (Ytr.center);
		\draw [style=simple] (Ytr.center) to (Ybr.center);
		\draw [style=simple] (Ybr.center) to (Ybl.center);
		\draw [style=simple] (Ybl.center) to (Ytl.center);
		\draw [style=simple] (Ztl.center) to (Ztr.center);
		\draw [style=simple] (Ztr.center) to (Zbr.center);
		\draw [style=simple] (Zbr.center) to (Zbl.center);
		\draw [style=simple] (Zbl.center) to (Ztl.center);
	\end{pgfonlayer}
\end{tikzpicture}
}
\qquad
\scalebox{0.6}{
\begin{tikzpicture}
	\begin{pgfonlayer}{nodelayer}

		\node [style = transition] (a) at (2.5, 0.0) {$\alpha$};
		\node [place,tokens=0] (A) at (1, 0.0) {$A$};
		\node [place,tokens=0] (B) at (4,0) {$B$};

		\node [style=empty] (Y) at (-0.1, 1) {$X_2$};
		\node [style=none] (Ytr) at (.25, 0.75) {};
		\node [style=none] (Ytl) at (-.4, 0.75) {};
		\node [style=none] (Ybr) at (.25, -0.75) {};
		\node [style=none] (Ybl) at (-.4, -0.75) {};

		\node [style=inputdot] (4) at (0, 0.0) {};
		\node [style=empty] at (-0.2, 0.0) {$1$};
		
		\node [style=empty] (Z) at (5.1, 1) {$Y_2$};
		\node [style=none] (Ztr) at (4.75, 0.75) {};
		\node [style=none] (Ztl) at (5.4, 0.75) {};
		\node [style=none] (Zbl) at (5.4, -0.75) {};
		\node [style=none] (Zbr) at (4.75, -0.75) {};

		\node [style=inputdot] (6) at (5, 0) {};
		\node [style=empty] at (5.2, 0) {$2$};	
		
	\end{pgfonlayer}
	\begin{pgfonlayer}{edgelayer}
		\draw [style=inarrow] (A) to (a);
		\draw [style=inarrow] (a) to (B);

		\draw [style=inputarrow] (4) to (A);
		\draw [style=inputarrow] (6) to (B);

		\draw [style=simple] (Ytl.center) to (Ytr.center);
		\draw [style=simple] (Ytr.center) to (Ybr.center);
		\draw [style=simple] (Ybr.center) to (Ybl.center);
		\draw [style=simple] (Ybl.center) to (Ytl.center);
		\draw [style=simple] (Ztl.center) to (Ztr.center);
		\draw [style=simple] (Ztr.center) to (Zbr.center);
		\draw [style=simple] (Zbr.center) to (Zbl.center);
		\draw [style=simple] (Zbl.center) to (Ztl.center);
	\end{pgfonlayer}
\end{tikzpicture}
}
\]
\else
We can also tensor open Petri nets, putting them side by side ``in parallel''.  This might suggest that open Petri nets should be the morphisms of a symmetric monoidal category.  However, composition of open Petri is not strictly associative, and there are very interesting maps \emph{between} open Petri nets.  To get a feeling for these, there is a morphism from this open Petri net:
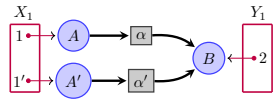
\[
\scalebox{0.6}{
\begin{tikzpicture}
	\begin{pgfonlayer}{nodelayer}

		\node [style = transition] (a) at (2.5, 0.5) {$\alpha$};
		\node [style = transition] (a') at (2.5, -0.5) {$\alpha'$};
		\node [place,tokens=0] (A) at (1, 0.5) {$A$};
		\node [place,tokens=0] (A') at (1, -0.5) {$A'$};
		\node [place,tokens=0] (B) at (4,0) {$B$};
	
		\node [style=empty] (Y) at (-0.1, 1) {$X_1$};
		\node [style=none] (Ytr) at (.25, 0.75) {};
		\node [style=none] (Ytl) at (-.4, 0.75) {};
		\node [style=none] (Ybr) at (.25, -0.75) {};
		\node [style=none] (Ybl) at (-.4, -0.75) {};

		\node [style=inputdot] (4) at (0, 0.5) {};
		\node [style=empty] at (-0.2, 0.5) {$1$};
		\node [style=inputdot] (5) at (0, -0.5) {};
		\node [style=empty] at (-0.2, -0.5) {$1'$};		
		
		\node [style=empty] (Z) at (5.1, 1) {$Y_1$};
		\node [style=none] (Ztr) at (4.75, 0.75) {};
		\node [style=none] (Ztl) at (5.4, 0.75) {};
		\node [style=none] (Zbl) at (5.4, -0.75) {};
		\node [style=none] (Zbr) at (4.75, -0.75) {};

		\node [style=inputdot] (6) at (5, 0) {};
		\node [style=empty] at (5.2, 0) {$2$};	
		
	\end{pgfonlayer}
	\begin{pgfonlayer}{edgelayer}

		\draw [style=inarrow] (A) to (a);
		\draw [style=inarrow, bend left=20, looseness=1.00] (a) to (B);
		\draw [style=inarrow] (A') to (a');
		\draw [style=inarrow, bend right=20, looseness=1.00] (a') to (B);
		\draw [style=inputarrow] (4) to (A);
		\draw [style=inputarrow] (5) to (A');
		\draw [style=inputarrow] (6) to (B);

		\draw [style=simple] (Ytl.center) to (Ytr.center);
		\draw [style=simple] (Ytr.center) to (Ybr.center);
		\draw [style=simple] (Ybr.center) to (Ybl.center);
		\draw [style=simple] (Ybl.center) to (Ytl.center);
		\draw [style=simple] (Ztl.center) to (Ztr.center);
		\draw [style=simple] (Ztr.center) to (Zbr.center);
		\draw [style=simple] (Zbr.center) to (Zbl.center);
		\draw [style=simple] (Zbl.center) to (Ztl.center);
	\end{pgfonlayer}
\end{tikzpicture}
}
\]
to this one:
\[
\scalebox{0.6}
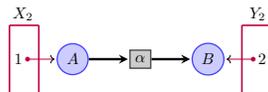
{
\begin{tikzpicture}
	\begin{pgfonlayer}{nodelayer}

		\node [style = transition] (a) at (2.5, 0.0) {$\alpha$};
		\node [place,tokens=0] (A) at (1, 0.0) {$A$};
		\node [place,tokens=0] (B) at (4,0) {$B$};

		\node [style=empty] (Y) at (-0.1, 1) {$X_2$};
		\node [style=none] (Ytr) at (.25, 0.75) {};
		\node [style=none] (Ytl) at (-.4, 0.75) {};
		\node [style=none] (Ybr) at (.25, -0.75) {};
		\node [style=none] (Ybl) at (-.4, -0.75) {};

		\node [style=inputdot] (4) at (0, 0.0) {};
		\node [style=empty] at (-0.2, 0.0) {$1$};
		
		\node [style=empty] (Z) at (5.1, 1) {$Y_2$};
		\node [style=none] (Ztr) at (4.75, 0.75) {};
		\node [style=none] (Ztl) at (5.4, 0.75) {};
		\node [style=none] (Zbl) at (5.4, -0.75) {};
		\node [style=none] (Zbr) at (4.75, -0.75) {};

		\node [style=inputdot] (6) at (5, 0) {};
		\node [style=empty] at (5.2, 0) {$2$};	
		
	\end{pgfonlayer}
	\begin{pgfonlayer}{edgelayer}
		\draw [style=inarrow] (A) to (a);
		\draw [style=inarrow] (a) to (B);

		\draw [style=inputarrow] (4) to (A);
		\draw [style=inputarrow] (6) to (B);

		\draw [style=simple] (Ytl.center) to (Ytr.center);
		\draw [style=simple] (Ytr.center) to (Ybr.center);
		\draw [style=simple] (Ybr.center) to (Ybl.center);
		\draw [style=simple] (Ybl.center) to (Ytl.center);
		\draw [style=simple] (Ztl.center) to (Ztr.center);
		\draw [style=simple] (Ztr.center) to (Zbr.center);
		\draw [style=simple] (Zbr.center) to (Zbl.center);
		\draw [style=simple] (Zbl.center) to (Ztl.center);
	\end{pgfonlayer}
\end{tikzpicture}
}
\]
\fi
mapping both primed and unprimed entities to the corresponding unprimed ones.  This particular morphism maps an open Petri net onto a simpler one.  There are also morphisms that include open Petri nets into more complicated ones.  For example, the above morphism has two right inverses: two ways to include the bottom open Petri net into the top one.

Given that we can compose open Petri net but also compose maps between them, it is natural to formalize them using a symmetric monoidal double category, which we call $\Open\Petri$.  To construct this one can use the following result on `structured cospans' \cite{BC,Courser}, which we state in summary form:

\begin{lem} {\bf \cite[Thm.\ 3.9]{BC}} 
\label{structured}
Let $\X$ be a category with finite colimits. Given a left adjoint $L \maps \Set \to \X$, there is a symmetric monoidal double category $\Open(\X)$ such that:
	\begin{itemize}
		\item objects are sets,
		\item vertical 1-morphisms are functions,
		\item a horizontal 1-cell from $a \in \Set$ to $b\in \Set$ is a cospan in $\X$ of this form:
			\[
			\begin{tikzpicture}[scale=1.5]
			\node (A) at (0,0) {$La$};
			\node (B) at (1,0) {$x$};
			\node (C) at (2,0) {$Lb$.};
			\path[->,font=\scriptsize,>=angle 90]
			(A) edge node[above,left]{$$} (B)
			(C)edge node[above]{$$}(B);
			\end{tikzpicture}
			\]
		\item{a 2-morphism is a commutative diagram in $\X$ of this form:
			\[
			\begin{tikzpicture}[scale=1.5]
			\node (E) at (3,0) {$La$};
			\node (F) at (5,0) {$Lb$};
			\node (G) at (4,0) {$x$};
			\node (E') at (3,-1) {$La'$};
			\node (F') at (5,-1) {$Lb'$.};
			\node (G') at (4,-1) {$x'$};
			\path[->,font=\scriptsize,>=angle 90]
			(F) edge node[above]{$$} (G)
			(E) edge node[left]{$Lf$} (E')
			(F) edge node[right]{$Lg$} (F')
			(G) edge node[left]{$h$} (G')
			(E) edge node[above]{$$} (G)
			(E') edge node[below]{$$} (G')
			(F') edge node[below]{$$} (G');
			\end{tikzpicture}
			\]
			}		
			\end{itemize}
Composition of vertical 1-morphisms is composition of functions.   Composition of horizontal 1-cells is composition of cospans in $\X$ via pushout.   Horizontal composition of 2-morphisms is also done via pushout, vertical composition of 2-morphisms is done via composition in $\X$, and the tensor product and symmetry are defined using chosen coproducts in $\Set$ and $\X$.
\end{lem}

We obtain the symmetric double category $\Open(\Petri)$ by applying this result to the functor $L \maps \Set \to \Petri$ that is left adjoint to the functor sending any 
Petri net to its set of places.  

The same sort of construction gives symmetric monoidal double categories of `open' versions of all our favorite kinds of nets and categories.  Moreover, the construction in \cref{structured} is functorial in the following sense:

\begin{lem} {\bf \cite[Thm.\ 4.2]{BC}} 
\label{thm:structured2}
Suppose $\X$ and $\X'$ have finite colimits and this triangle of finitely cocontinuous functors commutes up to natural isomorphism:
\[
\begin{tikzpicture}[scale=1.5]
\node (A) at (0,-0.5) {$\Set$};
\node (B) at (1.5,0) {$\X$};
\node (D) at (1.5,-1) {$\X'$.};
\node (E) at (1,-0.5) {$\NEarrow \phi$};
\path[->,font=\scriptsize,>=angle 90]
(A) edge node[above]{$L$} (B)
(B) edge node [right]{$F$} (D)
(A) edge node [below]{$L'$} (D);
\end{tikzpicture}
\]
Then there a symmetric monoidal double functor $\Open(F) \maps \Open(\X) \to \Open(\X')$ that acts as follows on 2-morphisms:
\[
\iflics \begin{tikzpicture}[scale=1.2] \else \begin{tikzpicture}[scale=1.5] \fi
\node (A) at (0,0) {$La$};
\node (B) at (1,0) {$x$};
\node (C) at (2,0) {$Lb$};
\node (A') at (0,-1) {$La'$};
\node (B') at (1,-1) {$x'$};
\node (C') at (2,-1) {$Lb'$};
\path[->,font=\scriptsize,>=angle 90]
(A) edge node[above]{$i$} (B)
(C)edge node[above]{$o$}(B)
(A') edge node[below]{$i'$} (B')
(C')edge node[below]{$o'$}(B')
(C)edge node[right]{$Lg$}(C')
(B)edge node[left]{$h$}(B')
(A)edge node[left]{$Lf$}(A');
\end{tikzpicture}
\raisebox{2.75em}{$\quad \mapsto \quad$}
\iflics \begin{tikzpicture}[scale=1.2] \else \begin{tikzpicture}[scale=1.5] \fi
\node (A) at (0,0) {$L'a$};
\node (B) at (1.5,0) {$Fx$};
\node (C) at (3,0) {$L'b$};
\node (A') at (0,-1) {$L'c$};
\node (B') at (1.5,-1) {$Fy$};
\node (C') at (3,-1) {$L'd$.};
\path[->,font=\scriptsize,>=angle 90]
(A) edge node[above]{$ Fi \circ \phi_a$} (B)
(C)edge node[above]{$ Fo \circ \phi_b $}(B)
(A') edge node[below]{$F i'\circ \phi_{a'}$} (B')
(C')edge node[below]{$ F o' \circ \phi_{b'}$}(B')
(A)edge node[left]{$ L'f $}(A')
(B)edge node[left]{$ Fh $}(B')
(C)edge node[right]{$L'g$}(C');
\end{tikzpicture}
\]
\end{lem}

It can be checked that given two composable triangles of the above sort:
\[
\begin{tikzpicture}[scale=1.5]
\node (A) at (0,0) {$\Set$};
\node (B) at (1.5,1) {$\X$};
\node (D) at (1.5,0) {$\X'$};
\node (E) at (1.1,0.3) {$\Narrow \phi$};
\node (F) at (1.1,-0.3) {$\Narrow \phi'$};
\node (G) at (1.5,-1) {$\X''$};
\path[->,font=\scriptsize,>=angle 90]
(A) edge node[above]{$L$} (B)
(B) edge node [right]{$F$} (D)
(A) edge node [above]{$L'$} (D)
(D) edge node [right]{$F'$} (G)
(A) edge node [below]{$L''$} (G);
\end{tikzpicture}
\]
we have 
\[   \Open(F' \circ F) = \Open(F') \circ \Open(F) .\]
Furthermore, isomorphic triangles give isomomorphic symmetric monoidal double functors, so we obtain the following result:

\begin{thm}
\label{thm:symmetric_monoidal_double_functor}
There is a diagram of symmetric monoidal double functors
\[
\adjustbox{scale=0.8}{
\begin{tikzcd}[row sep=large,column sep=4em]
  \Open(\StrMC) \ar[r,shift left,"\Open(F_\StrMC)"]  &
  \Open(\SSMC)  \ar[r,shift left,"\Open(F_\SSMC)"]  &
  \Open(\CMC)   \\
  \\
  \Open(\PreNet) 
  \ar[uu,"\Open(F_\PreNet)"] \ar[r,shift left,"\Open(F_\pre)"] 
  &\Open(\SigmaNetCat) \ar[uu,"\Open(F_\SigmaNetCat)"] \ar[l,shift left,"\Open(G_\pre)"] \ar[r,"\Open(F_\pet)"]
  &\Open(\Petri) \ar[uu,"\Open(F_\Petri)"]
\end{tikzcd}}
\]
where we get each arrow from one of the left adjoints in \cref{eq:adjunctions} using \cref{thm:structured2}, and the squares built from arrows going up or right commute up to 2-isomorphism.
\end{thm}\fi

\section{Conclusion and future work}

In this work we have systematized the theory of Petri nets, their variants, and their categorical semantics.  To this end, we have shown that the notion of \SigmaNet{}, almost absent from standard Petri net literature, is in fact central.  Our framework gives a consistent view of the relations between these interacting notions of net in terms of adjunctions, such that the most important adjunctions present in the literature can be recovered as composites of our fundamental ones.

Our work makes substantial use of tools from homotopy theory and related fields, such as groupoids and fibrations.  We believe this will open up exciting new directions of research in the study of distributed systems and network theory in general.

In fact, the relationships between the various notions of net in our work have analogues in topology. On one hand, a manifold can always be given ``local coordinates'', but it is too restrictive to ask that such coordinates be preserved strictly by maps between manifolds.  Such coordinates can be regarded as analogous to the orderings on sources and targets in a pre-net. On the other hand, when a group acts on a manifold, the quotient topological space may no longer be a manifold, but has singularities at points of non-free action.  This ``coarse moduli space'' can be regarded as analogous to a Petri net, where symmetry information has been lost. Kock's whole-grain Petri nets are analogous to abstract manifolds themselves: they are free of undesirable ``coordinates'', but neither can they have singularities.
Finally, our $\Sigma$-nets play the role of \emph{orbifolds}, coordinate-free manifold-like structures that retain the information of ``isotropy groups'' at singular points, yielding a better-behaved notion of quotient.

\iflics\else
\section*{Acknowledgements}
The second author was supported by the project MIUR PRIN 2017FTXR7S “IT-MaTTerS” and by the \href{https://gitcoin.co/grants/1086/independent-ethvestigator-program}{Independent Ethvestigator Program}.

The fourth author was supported by The United States Air Force Research Laboratory under agreement number FA9550-15-1-0053.  The U.S. Government is authorized to reproduce and distribute reprints for Governmental purposes notwithstanding any copyright notation thereon.  The views and conclusions contained herein are those of the author and should not be interpreted as necessarily representing the official policies or endorsements, either expressed or implied, of the United States Air Force Research Laboratory, the U.S. Government, or Carnegie Mellon University.
\fi

\end{document}